\documentclass[11pt,reqno,a4paper]{amsart}

\usepackage[utf8]{inputenc}
\usepackage[T1]{fontenc}
\usepackage{amsmath,textcomp,  microtype, mathrsfs,csquotes}
\IfFileExists{fbb.sty}{\usepackage{fbb}}{\usepackage{lmodern}}
\usepackage{esint}
\usepackage{amsthm,amssymb,amsfonts}
\usepackage{mathtools}
\usepackage{enumitem,needspace,etoolbox}
\allowdisplaybreaks

\newcommand{\R}{\mathbb{R}}
\newcommand{\Sph}{\mathbb{S}}

\DeclareMathOperator{\Ric}{Ric}
\numberwithin{equation}{section}
\usepackage[top=2.5cm,bottom=2.4cm,left=2.6cm,right=2.6cm,headsep=0.2in]{geometry}
\usepackage{fancyhdr}
\usepackage[pagebackref]{hyperref}
\hypersetup{
colorlinks=true, linkcolor=blue!50!green, citecolor=blue!50!green, urlcolor=blue!50!green, pdfauthor={Mayukh Mukherjee and Utsab Sarkar}, pdftitle={One-bubble stability, reduced coefficients, and threshold compactness for the Escobar functional} }
\usepackage{xcolor}
\newcommand{\dist}{\mathrm{dist}}
\newcommand{\Def}{\mathrm{Def}}
\theoremstyle{plain}
\newtheorem{theorem}{Theorem}[section]
\newtheorem{lemma}[theorem]{Lemma}
\newtheorem{proposition}[theorem]{Proposition}
\newtheorem{corollary}[theorem]{Corollary}
\newtheorem{definition}[theorem]{Definition}

\newcommand{\Scal}{\mathrm{Scal}}

\newcommand{\Esc}{\ensuremath{\mathrm{Esc}}}

\theoremstyle{remark}
\newtheorem{remark}[theorem]{Remark}

\newtheorem{convention}[theorem]{Convention}

\theoremstyle{plain}
\newtheorem{introthm}{Theorem}

\title[One-bubble stability and threshold compactness]{One-bubble stability, reduced coefficients, and threshold compactness for the Escobar functional}
\author{Mayukh Mukherjee and Utsab Sarkar}
\address{MM: Indian Institute of Technology Bombay, Powai, Maharashtra 400076, India}
\email{mathmukherjee@gmail.com, mukherjee@math.iitb.ac.in}
\address{US: Indian Institute of Technology Bombay, Powai, Maharashtra 400076, India}
\email{utsab@math.iitb.ac.in, reachutsab@gmail.com}
\date{}

\begin{document}
\begin{abstract}
Let $(M^n,g)$ be a smooth compact manifold with nonempty boundary and coercive conformal energy. For $n\ge5$, in a conformal representative with $H_g=0$, we prove quantitative stability near a concentrating boundary bubble whose center stays away from the umbilic locus. The norm of the constrained gradient controls the distance to the standard bubble family with sharp exponent $1/2$, and the distance to the Lyapunov--Schmidt corrected family linearly. In the same one-bubble neighborhood, at energies at least the sharp hemisphere value $S_\ast$, the distance to the standard family is also controlled linearly. These estimates follow from the exact reduced quadratic coefficient and the scale-derivative expansion of the reduced functional. The same scale equation shows that sequences of critical points at $S_\ast$ can concentrate only at umbilic points. For $n\ge3$, a profile decomposition without a sign assumption gives Palais--Smale compactness between $S_\ast$ and the first mixed level $(S_\ast^{n-1}+Y^{n-1})^{1/(n-1)}$, where $Y=C^*_{\Esc}(M,g)$, and this endpoint is sharp. Between $S_\ast$ and the two-bubble level, we characterize noncompact levels by the lower critical values. For $n\ge7$, we obtain the fourth-order reduced expansion at umbilic centers. Under local total umbilicity, we evaluate the correction exactly, determine the conformal weight of the reduced coefficient, and obtain the optimal algebraic comparison with the squared ambient Weyl norm.
\end{abstract}
\keywords{Escobar functional, one-bubble stability, Lyapunov--Schmidt reduction, threshold compactness, Weyl curvature}
\subjclass{53C21 (primary); 35B35, 35B44, 35J65 (secondary).}
\maketitle

\section{Introduction}
The scalar-flat boundary Yamabe problem asks whether a compact Riemannian manifold with boundary admits a conformal metric with zero scalar curvature and constant boundary mean curvature. Let $(M^n,g)$ be a smooth compact Riemannian manifold with nonempty boundary, $n\ge3$, and let $q=2(n-1)/(n-2)$ be the critical trace exponent. Denote the scalar curvature by $\Scal_g$, the averaged boundary mean curvature by $H_g$, and the volume and boundary measures by $dV_g$ and $d\sigma_g$. We use the convention for which the Euclidean unit ball has $H_g=1$. For $u\in H^1(M)$, the conformal energy is
\begin{equation}\label{B-recall-I-eq:Esc-energy-weak}
\mathcal N_g^\circ(u)=\int_M\left(|\nabla u|_g^2+\frac{n-2}{4(n-1)}\Scal_gu^2\right)dV_g
+\frac{n-2}{2}\int_{\partial M}H_gu^2\,d\sigma_g.
\end{equation}
Let $\operatorname{Tr}u$ denote the boundary trace. The Escobar quotient and its infimum are
\begin{equation}\label{B-recall-I-eq:Esc-sharp-constant}
\mathcal J_g[u]=\frac{\mathcal N_g^\circ(u)}{\|u\|_{L^q(\partial M)}^2},
\qquad C^*_{\Esc}(M,g)=\inf_{\substack{u\in H^1(M)\\\operatorname{Tr}u\ne0}}\mathcal J_g[u].
\end{equation}
Write $S_\ast$ for this infimum on the round hemisphere, equivalently the sharp half-space trace constant. Rescaled half-space optimizers, called boundary bubbles, have quotient tending to $S_\ast$. Their concentration accounts for the loss of compactness at this value. We study how the norm of the constrained gradient measures proximity to one concentrating bubble, and which geometric terms prevent such a bubble from being stationary.

Throughout the analytic arguments we assume strict coercivity, denoted \textup{(Y$^+_\partial$)}:
\begin{equation}\label{B-recall-I-eq:Yplus}
\mathcal N_g^\circ(u)\ge c_Y\|u\|_{H^1(M)}^2
\qquad(u\in H^1(M)),\qquad c_Y>0.
\end{equation}
Polarization defines the energy inner product. The regularity and geometric bounds are stated with each result. The argument has two parts. A profile decomposition without a sign assumption determines the possible loss of compactness below the first mixed level. At the remaining single-bubble level $S_\ast$, we solve the equation in the complementary directions and study the resulting energy as a function of the center and scale. A Palais--Smale sequence may concentrate at any prescribed boundary point; a concentrating sequence of exact critical points must also make this reduced energy stationary. The quadratic and fourth-order expansions give successive geometric obstructions to this stationarity, under the hypotheses stated below.

Escobar developed the variational approach \cite{EscobarJDG92} and its scalar-flat formulation \cite{EscobarAnnals92,EscobarAddendum94,Escobar1996a}; the half-space optimizers are given in \cite[Eq.~(2)]{Escobar1988}. Write $\mathring{\mathrm{II}}_g$ for the trace-free second fundamental form and $W_g$ for the Weyl tensor. A point is umbilic if $\mathring{\mathrm{II}}_g=0$ there, and the boundary is totally umbilic if this holds everywhere. The strict comparison at a non-umbilic point was proved by Escobar for $n\ge6$ and Marques for $n=4,5$; see \cite[pp.~383--384]{Marques2007} and \cite[pp.~1--2]{Almaraz2010}. Almaraz used fourth-order test functions for umbilic boundaries with $W_g\ne0$ \cite[Section~3]{Almaraz2010}. Under a curvature-jet vanishing hypothesis, Mayer and Ndiaye proved existence for positive conformal invariant and totally umbilic boundary when $n\ge6$ \cite[Theorem~1.1]{MayerNdiaye2017}.

Deficit estimates controlling the squared distance from extremals were proved by Bianchi and Egnell for the Euclidean Sobolev inequality \cite{BianchiEgnell1991} and by Ho for the half-space trace inequality \cite{Ho2022}; see also \cite[Section~2]{BorquezCajuBosch2026}. Fractional Sobolev analogues were proved by Chen, Frank, and Weth \cite{ChenFrankWeth2013}. Zhang, Zhou, and Zou proved stability and profile decomposition for the flat trace problem \cite[Theorems~1.2, 1.6, and~5.1]{ZhangZhouZou2025}. On compact manifolds with nonnegative scalar curvature and conformal invariant below $S_\ast$, Borquez, Caju, and Van Den Bosch obtained estimates near the positive minimizing set \cite[Theorem~1.1]{BorquezCajuBosch2026}. Quantitative stability and threshold compactness for the homogeneous Dirichlet Gagliardo--Nirenberg quotient on manifolds with boundary are studied in \cite{MukherjeeSarkarGN2026}.

For the Euclidean Sobolev equation, sharp quantitative stability estimates were proved by Figalli and Glaudo \cite{FigalliGlaudo2020} and Deng, Sun, and Wei \cite{DengSunWei2025}. Chen and Kim proved sharp estimates on closed manifolds with rates depending on the geometry and the concentrating functions \cite{ChenKim2024}. Quantitative stability estimates for the flat trace problem appear in \cite{ZhangZhouZou2025}; Yu and Zhou treat one bubble for the $p$-Laplacian \cite{YuZhou2026}. We determine how the norm of the constrained gradient controls the distance to the standard and corrected bubble families on a compact manifold with boundary.

Positive-solution compactness was proved under various hypotheses by Almaraz \cite{Almaraz2011CVPDE}, Almaraz, de Queiroz, and Wang \cite{AlmarazQueirozWang2019}, Kim, Musso, and Wei \cite{KimMussoWei2021}, and Ghimenti and Micheletti \cite{GhimentiMicheletti2020NA,GhimentiMicheletti2021CVPDE}. Gong, Kim, Musso, and Wei establish compactness through dimension $14$ for general boundary and $21$ for umbilic boundary, excluding the conformal hemisphere \cite[Theorem~1]{GongKimMussoWei2026}. Our stability estimates concern approximate critical points and distinguish the standard bubble family from the corrected family.

Let $\mathcal S=\{u\in H^1(M):\|u\|_{L^q(\partial M)}=1\}$. Polarization of $\mathcal N_g^\circ$ gives the energy norm $\|\cdot\|_{E,g}$, and we put
\[
r_g(u)=\|\mathsf{grad}_{\mathcal S}\mathcal J_g(u)\|_{E,g}.
\]
A constrained Palais--Smale sequence at $c$ satisfies $\mathcal J_g[u_j]\to c$ and $r_g(u_j)\to0$. Throughout, the smooth compact metric satisfies \textup{(Y$^+_\partial$)}, and constants may depend on $(M,g)$. The finite-regularity results are stated in Theorem~\ref{II-thm:generic-compactness} and proved in Appendix~\ref{app:genericity-proof}.

For the bubble expansions, choose a conformal representative with $H_g=0$. Let $\mathcal U_{x,\varepsilon}\in\mathcal S$ be the normalized boundary bubble centered at $x$ and of scale $\varepsilon$, with the cutoff specified in Section~\ref{B-recall-I-sec:escobar}. Write $Z_0=\varepsilon\partial_\varepsilon\mathcal U_{x,\varepsilon}$ and $Z_a=\varepsilon\partial_{x^a}\mathcal U_{x,\varepsilon}$ in local boundary coordinates, and set
\[
X_{x,\varepsilon}=\left\{h\in H^1(M):
\int_{\partial M}\mathcal U_{x,\varepsilon}^{q-1}h\,d\sigma_g=0,
\quad \langle h,Z_j\rangle_{E,g}=0\quad(0\le j\le n-1)\right\}.
\]
The correction solves the equation transverse to changes of center and scale. For $h\in X_{x,\varepsilon}$ put $\mathsf{ret}_{\mathcal U}(h)=(\mathcal U+h)/\|\mathcal U+h\|_q$, the normalization map (a local retraction onto $\mathcal S$). When $n\ge5$, there is a unique small correction $w_{x,\varepsilon}\in X_{x,\varepsilon}$ for which the differential of $\mathcal J_g\circ\mathsf{ret}_{\mathcal U_{x,\varepsilon}}$ vanishes on $X_{x,\varepsilon}$. Define
\[
\Phi(x,\varepsilon)=\mathsf{ret}_{\mathcal U_{x,\varepsilon}}(w_{x,\varepsilon}),
\qquad \mathcal J_1(x,\varepsilon)=\mathcal J_g[\Phi(x,\varepsilon)].
\]
Let $\mathcal B_0$ and $\mathcal B_{\mathrm{LS}}$ be the standard and Lyapunov--Schmidt (LS) corrected families at sufficiently small scales, including both signs. All distances below are measured in the energy norm.

\begin{introthm}[Quantitative stability and sharpness]\label{I-thm:main-residual}
Assume $n\ge5$ and $H_g=0$. Let $K\subset\partial M$ be compact with $\inf_K|\mathring{\mathrm{II}}_g|>0$. There are $\rho,\varepsilon_0,C>0$ such that every $u\in\mathcal S$ whose normalized modulation gives
\[
u=\sigma\mathsf{ret}_{\mathcal U_{x,\varepsilon}}(v),\qquad
v\in X_{x,\varepsilon},\quad \|v\|_{E,g}<\rho,\quad
x\in K,\quad 0<\varepsilon<\varepsilon_0,\quad \sigma\in\{\pm1\},
\]
satisfies
\begin{equation}\label{I-eq:intro-residual-estimates}
\begin{gathered}
\|u-\sigma\Phi(x,\varepsilon)\|_{E,g}\le Cr_g(u),
\qquad \varepsilon^2\le Cr_g(u),\\
\dist_E(u,\mathcal B_0)\le Cr_g(u)^{1/2},
\qquad \dist_E(u,\mathcal B_{\mathrm{LS}})\le Cr_g(u).
\end{gathered}
\end{equation}
If $\mathcal J_g[u]\ge S_\ast$, then
\[
\dist_E(u,\mathcal B_0)\le Cr_g(u).
\]
Whenever $\mathring{\mathrm{II}}_g$ is not identically zero on a connected boundary component, there are concentrating sequences showing that the powers $1/2$ in the standard-family estimate and $1$ in its improvement above $S_\ast$ cannot be increased. Furthermore, these statements hold independently of any sign assumptions on the underlying functions.
\end{introthm}

Theorem~\ref{I-thm:main-residual} is proved in Theorem~\ref{II-thm:residual-stability} and Proposition~\ref{II-prop:residual-sharpness}. The complementary correction has order $\varepsilon$, whereas the logarithmic scale derivative of the reduced functional has order $\varepsilon^2$ away from the umbilic locus. The negative energy term in the next theorem accounts both for the square-root estimate and for its linear improvement above $S_\ast$.

\begin{introthm}[The exact quadratic coefficient]\label{I-thm:main-quadratic}
Assume $n\ge5$ and $H_g=0$. Uniformly in $x\in\partial M$,
\begin{align*}
\mathcal J_1(x,\varepsilon)/S_\ast
&=1+\kappa_3^{\mathrm{red}}(n)|\mathring{\mathrm{II}}_g(x)|^2\varepsilon^2
+o(\varepsilon^2),\\
S_\ast^{-1}\partial_\varepsilon\mathcal J_1(x,\varepsilon)
&=2\kappa_3^{\mathrm{red}}(n)|\mathring{\mathrm{II}}_g(x)|^2\varepsilon
+o(\varepsilon).
\end{align*}
The reduced coefficient is strictly negative. Its exact value and the complementary correction are
\begin{align}
\delta\kappa_3^{\mathrm{LS}}(n)
&=\frac{(2n-1)(n+2)}{2(n-1)^2(n+1)(n-3)(n-4)},\nonumber\\
\kappa_3^{\mathrm{red}}(n)
&=-\frac{(n-2)(n^2-2n-2)}{2(n-1)^2(n+1)(n-3)(n-4)}<0.
\label{I-eq:intro-reduced-coefficient}
\end{align}
The standard-bubble coefficient is $\kappa_3^{\mathrm{bare}}(n)=(6-n)/(2(n-1)(n-3)(n-4))$, and $\kappa_3^{\mathrm{red}}=\kappa_3^{\mathrm{bare}}-\delta\kappa_3^{\mathrm{LS}}$. In particular, $\delta\kappa_3^{\mathrm{LS}}(5)=21/128$ and $\kappa_3^{\mathrm{red}}(5)=-13/128$.
\end{introthm}

The quadratic form involving the inverse linearized operator is evaluated by an explicit particular solution and a harmonic Robin correction, using spherical harmonics as in \cite[Proposition~3.3]{GongKimMussoWei2026}. Its evaluation and differentiated expansion are proved in Proposition~\ref{I-prop:quadratic-exact} and Proposition~\ref{II-prop:onebubble-scale-derivative-II}.

The compactness statements use a profile decomposition without a sign assumption, with a remainder converging strongly to zero in $H^1$, based on \cite[Lemma~2.4]{Almaraz2014}. Put
\[
Y=C^*_{\Esc}(M,g)>0,\qquad
c_{\mathrm{mix}}=(S_\ast^{n-1}+Y^{n-1})^{1/(n-1)},\qquad
\Lambda_2=2^{1/(n-1)}S_\ast.
\]

\begin{introthm}[Palais--Smale levels and compactness of critical points]\label{I-thm:main-compactness}\leavevmode\begin{enumerate}[label=\textup{(\roman*)},leftmargin=2em]
\item Assume $n\ge3$. Every constrained Palais--Smale sequence with level $c<S_\ast$ or $S_\ast<c<c_{\mathrm{mix}}$ is precompact in $H^1(M)$. The endpoint $c_{\mathrm{mix}}$ admits a noncompact Palais--Smale sequence.

For $S_\ast<c<\Lambda_2$, a noncompact constrained Palais--Smale sequence at $c$ exists if and only if
\[
d=(c^{n-1}-S_\ast^{n-1})^{1/(n-1)}
\]
is a critical value of $\mathcal J_g|_{\mathcal S}$.

\item Assume $n\ge5$. In a representative with $H_g=0$, let $(u_j)\subset\mathcal S$ be constrained critical points with $\mathcal J_g[u_j]\to S_\ast$ that are not precompact. After extraction there are $\sigma\in\{\pm1\}$, $\varepsilon_j\downarrow0$, and $x_j\to p\in\partial M$ such that
\[
u_j=\sigma\Phi(x_j,\varepsilon_j),\qquad
\mathring{\mathrm{II}}_g(p)=0.
\]
The umbilicity conclusion is conformally invariant. In any conformal representative, if the boundary is nowhere umbilic, then for every $c_0<c_{\mathrm{mix}}$ the set
\[
\{u\in\mathcal S:\mathsf{grad}_{\mathcal S}\mathcal J_g(u)=0,
\quad \mathcal J_g[u]\le c_0\}
\]
is compact in $H^1(M)$.
\end{enumerate}
\end{introthm}

The Palais--Smale assertions are Theorems~\ref{I-thm:compactness-escobar} and~\ref{I-thm:signfree-first-gap-compactness}. A stronger mass bound for sign-changing profiles was proved by Almaraz and Wang \cite[Theorem~1.4]{AlmarazWang2022}. The critical-point conclusions follow from Theorem~\ref{II-thm:threshold-Sstar-onebubble} and Corollary~\ref{II-cor:threshold-compactness}. For $n\ge5$, the nowhere-umbilic condition is open and dense among coercive $C^m$ metrics, $m\ge3$; see Theorem~\ref{II-thm:generic-compactness}.

The preceding theorem leaves umbilic points as possible blow-up points of sequences of critical points. At an umbilic center, the order-$\varepsilon$ coefficient of the complementary first variation vanishes. The next coefficient can occur at order $\varepsilon^2$, and the correction obtained by completing the square contributes at order $\varepsilon^4$. We call these leading covectors the geometric sources; their precise definitions are given in Remark~\ref{I-rem:geometric-source}. The resulting quartic coefficient has an explicit curvature formula under local total umbilicity.

For the fourth-order expansion assume $n\ge7$ and take $R_\beta(\varepsilon)=\varepsilon^{-\beta}$ with $4/(n-2)<\beta<1$, and let $\mathcal J_{1,\beta}$ be the corresponding corrected one-bubble energy. Write $E_g=W_g(\nu,\cdot,\nu,\cdot)$ on $T\partial M$, where $\nu$ is a unit normal, and let $\overline W_g$ be the intrinsic Weyl tensor of the boundary metric.

\begin{introthm}[The fourth-order coefficient]\label{I-thm:main-quartic}
Assume $n\ge7$ and $H_g=0$.
\begin{enumerate}[label=\textup{(\roman*)},leftmargin=2em]
\item At every umbilic point $p$ there is a coefficient $\Xi_g^{\mathrm{red}}(p)$, independent of $\beta$, such that
\begin{align*}
\mathcal J_{1,\beta}(p,\varepsilon)/S_\ast
&=1+\Xi_g^{\mathrm{red}}(p)\varepsilon^4+o(\varepsilon^4),\\
S_\ast^{-1}\partial_\varepsilon\mathcal J_{1,\beta}(p,\varepsilon)
&=4\Xi_g^{\mathrm{red}}(p)\varepsilon^3+o(\varepsilon^3).
\end{align*}
The remainders are uniform on compact subsets of the umbilic locus.

\item If the boundary is totally umbilic near $p$, then
\begin{equation}\label{I-eq:intro-quartic-intrinsic}
\Xi_g^{\mathrm{red}}(p)
=-e_n\frac{(n-2)^2}{(n-3)^2}|E_g(p)|_g^2
-b_n|\overline W_g(p)|_g^2,
\end{equation}
where
\[
e_n=\frac{3n^3-18n^2+21n+14}
{3(n-2)(n-1)^2(n+1)(n-4)(n-5)(n-6)},
\quad
b_n=\frac{n-2}{24(n-1)(n-3)(n-5)(n-6)}.
\]
Both constants are positive. The optimal algebraic comparison is
\[
\Xi_g^{\mathrm{red}}(p)\le-c(n)|W_g(p)|_g^2,
\qquad c(n)=\frac{n-2}{4(n-3)}e_n,
\]
with equality precisely when $\overline W_g(p)=0$. In particular, $\Xi_g^{\mathrm{red}}(p)=0$ if and only if $W_g(p)=0$. For another representative $\widetilde g=\psi^{4/(n-2)}g$ with $H_{\widetilde g}=0$,
\[
\Xi_{\widetilde g}^{\mathrm{red}}(p)
=\psi(p)^{-8/(n-2)}\Xi_g^{\mathrm{red}}(p).
\]
\item Suppose that the whole boundary is totally umbilic. If a normalized sequence of constrained critical points at values tending to $S_\ast$ is not precompact, then every limiting blow-up center $p$ satisfies $W_g(p)=0$. If $\inf_{\partial M}|W_g|>0$, then for every $c_0<c_{\mathrm{mix}}$ the normalized constrained critical points with values at most $c_0$ form a compact subset of $H^1(M)$.
\end{enumerate}
\end{introthm}

The expansion is Proposition~\ref{II-prop:quartic-reduced-umbilic}; the curvature formula, conformal weight, and compactness consequence are Theorem~\ref{II-thm:quartic-Weyl-sign} and Corollaries~\ref{II-cor:quartic-conformal-weight} and~\ref{II-cor:quartic-subcriticality-compactness}.

The two coefficient calculations use the same flat complementary Hessian. If $\mathscr S$ and $\mathscr L$ are the complementary differential and Hessian of the functional, respectively, solving the complementary equation changes the energy by
\[
-\frac12\langle\mathscr L^{-1}\mathscr S,\mathscr S\rangle
\]
to leading order. The source is $O(\varepsilon)$ in the quadratic calculation and $O(\varepsilon^2)$ at umbilic centers with the cutoff used in Theorem~\ref{I-thm:main-quartic}. We control the remainder and its weighted scale derivative. The explicit curvature calculation in Theorem~\ref{I-thm:main-quartic} uses the same harmonic Robin pairing as Theorem~\ref{I-thm:main-quadratic}.

We use $\nabla_\partial$ for the gradient of the induced boundary metric. The dimensional restrictions reflect the corresponding integrability: the profile decomposition holds for $n\ge3$, the first geometric source has finite norm in the dual of the energy space for $n\ge5$, and the fourth-order analysis requires $n\ge7$. The standard-bubble calculations used below are stated in Proposition~\ref{I-prop:standard-bubble-inputs}, with their sources in the companion paper \cite{EscobarGeometry}.

After proving coercivity and modulation, we establish the profile decomposition without a sign assumption. We then construct the corrected family and compute its quadratic coefficient. The differentiated expansion gives the stability estimates and the compactness results for critical points. The fourth-order calculation continues this argument at umbilic centers and concludes the main text. Supplementary deficit estimates, hemisphere rigidity, the finite-regularity genericity proof, and further consequences are given in the appendices.

\section{Analytic and geometric preliminaries}
\label{B-recall-I-sec:escobar}
Throughout, $M$ is a compact smooth manifold of dimension $n\ge3$, with nonempty smooth boundary, and $g$ is a smooth Riemannian metric. Let $\bar g=g|_{T\partial M}$ be the induced boundary metric and $\nu$ the inward unit normal. Our conventions are
\[
\mathrm{II}(X,Y)=-\langle\nabla_X\nu,Y\rangle_g,
\qquad K_g=\operatorname{tr}_{\bar g}\mathrm{II},
\qquad H_g=\frac{K_g}{n-1},
\qquad \mathring{\mathrm{II}}=\mathrm{II}-H_g\bar g.
\]
Thus the Euclidean unit ball has $H_g=1$; see \cite[Remark~3.12]{EscobarGeometry}. 

\subsection{Energy and geometric assumptions}
\label{B-recall-I-subsec:escobar-notation}
Retain $q$ from the introduction, also denoted $2^*_\partial$, and put $a_n=4(n-1)/(n-2)$. The conformal operators are
\begin{equation}\label{B-recall-I-eq:canonical-pair}
L_g^\circ=-\Delta_g+a_n^{-1}\Scal_g,
\qquad B_g^\circ=-\partial_\nu+\frac{n-2}{2}H_g.
\end{equation}
The polarization of \eqref{B-recall-I-eq:Esc-energy-weak} is
\begin{equation}\label{B-recall-I-eq:Esc-energy-bilinear}
\mathcal N_g^\circ(u,v)=\int_M\bigl(\langle\nabla u,\nabla v\rangle_g+a_n^{-1}\Scal_guv\bigr)\,dV_g
+\frac{n-2}{2}\int_{\partial M}H_guv\,d\sigma_g.
\end{equation}
We suppress the subscript $g$ when the metric is fixed. Green's identity gives, for smooth $u$,
\begin{equation}\label{B-recall-I-eq:bridge-firstorder}
\int_M uL_g^\circ u\,dV_g+\int_{\partial M}uB_g^\circ u\,d\sigma_g
=\mathcal N_g^\circ(u),
\qquad
\mathcal Q_g(u)=a_n\mathcal N_g^\circ(u),
\end{equation}
where $\mathcal Q_g(u)=a_n\int_M|\nabla u|^2+\int_M\Scal_gu^2+2\int_{\partial M}K_gu^2$. For $H^1$ functions, the operator pairing in \eqref{B-recall-I-eq:bridge-firstorder} is defined by the energy bilinear form.

We retain the strict coercivity hypothesis \textup{(Y$^+_\partial$)} of \eqref{B-recall-I-eq:Yplus}.
The bounded coefficients and the trace inequality also give
\begin{equation}\label{B-recall-I-eq:Esc-energy-continuity}
\mathcal N_g^\circ(u)\le C\|u\|_{H^1(M)}^2.
\end{equation}
Thus $\langle u,v\rangle_{E,g}=\mathcal N_g^\circ(u,v)$ is an inner product, and $\|u\|_{E,g}=\sqrt{\mathcal N_g^\circ(u)}$ is equivalent to the $H^1$ norm.

The quotient and sharp constant are defined in \eqref{B-recall-I-eq:Esc-sharp-constant}.

If $C_{\mathrm{tr}}$ is the norm of the trace map $H^1(M)\to L^q(\partial M)$, then
\begin{equation}\label{B-recall-I-eq:Esc-strict-positive}
C^*_{\Esc}(M,g)\ge c_YC_{\mathrm{tr}}^{-2}>0.
\end{equation}
We write $\Def_g^{\Esc}(u)=\mathcal J_g[u]-C^*_{\Esc}(M,g)$. The hypothesis \textup{(Y$^+_\partial$)} is standard in the boundary Yamabe problem; see \cite{EscobarAnnals92,EscobarJDG92,HanLi1999}.

Write $\operatorname{Rm}_g$ and $\operatorname{Ric}_g$ for the Riemann and Ricci curvature tensors, respectively. The bounded-geometry assumption \textup{(BG$^{k+}$)} means that there is a uniform radius $r_0>0$ for boundary normal charts, interior geodesic charts, and the inward Fermi collar, together with bounds
\[
\sup_M|\nabla^j\operatorname{Rm}_g|\le C_j\quad(0\le j\le k-2),
\qquad
\sup_{\partial M}|\nabla_\top^j\mathrm{II}|\le\widetilde C_j\quad(0\le j\le k-1).
\]
More precisely, $\operatorname{inj}(\partial M,\bar g)\ge r_0$, a center at distance at least $r\le r_0$ from the boundary has an interior chart of radius $r$, and the normal exponential map is a diffeomorphism from $\partial M\times[0,r_0)$ onto the collar. Here $\nabla_\top$ is the boundary Levi-Civita connection. The trace and mixed coordinate estimates are given in \cite[Lemmas~3.1--3.2]{EscobarGeometry}. For a fixed smooth compact manifold their constants may depend on $(M,g)$ and the chart radii. For sharp manifold trace inequalities with a lower-order boundary term, see \cite{LiZhu1997}.

\smallskip\noindent \emph{Curvature conventions}.\label{B-recall-I-rem:intro-H-convention}
We take \[R(X,Y)Z=\nabla_X\nabla_YZ-\nabla_Y\nabla_XZ-\nabla_{[X,Y]}Z\] and $R_{ijkl}=\langle R(\partial_k,\partial_l)\partial_j,\partial_i\rangle_g$. The round sphere therefore has $R_{ijkl}=g_{ik}g_{jl}-g_{il}g_{jk}$. Let $\Sigma_t$ be the parallel hypersurface at inward distance $t$, and differentiate the Fermi components of its second fundamental form. At $t=0$,
\[
\partial_t(\mathrm{II}_{\Sigma_t})_{ab}
=-(\mathrm{II}^2)_{ab}+R_{anbn},
\qquad
(\nabla_\nu\mathrm{II})_{ab}=(\mathrm{II}^2)_{ab}+R_{anbn}.
\]
The second identity uses the covariant normal derivative of the same tensor family. The Gauss equation is
\[
\Scal_{\bar g}=\Scal_g|_{\partial M}-2\Ric_g(\nu,\nu)+K_g^2-|\mathrm{II}|^2.
\]
Tangential indices $a,b$ range from $1$ to $n-1$, repeated indices are summed, and tangential indices are raised with $\bar g$.

\subsection{Model profiles and boundary bubbles}
We use the closed half-space $\R^n_+=\R^{n-1}\times[0,\infty)$ and the round hemisphere $S^n_+=\{x\in S^n:x_{n+1}\ge0\}$ as manifolds with boundary. Equations in either domain refer to its interior. Write $y=(y',t)$, also denoted $Y$ in rescaled calculations, $B_R^+=B_R\cap\R^n_+$, and $\nabla_{\tan}=(\partial_{y_a})_{a=1}^{n-1}$. Flat volume and boundary measures are denoted by $dy$ and $d\sigma$.

The space $\mathcal D^{1,2}(\R^n_+)=\dot H^1(\R^n_+)$ is the completion, in the norm $\|\nabla w\|_{L^2}$, of restrictions of $C_c^\infty(\R^n)$ to $\R^n_+$. Its elements belong to $L^{2n/(n-2)}(\R^n_+)$ and have traces in $L^q(\partial\R^n_+)$. We use the analogous definition on $\R^n$.

The hemisphere constant $S_\ast$ also equals the sharp flat half-space constant \cite{Escobar1988,Beckner1993}. With $\Gamma$ denoting Euler's gamma function, we fix the optimizer \cite[Eq.~(2)]{Escobar1988}
\[
U_+(y',t)=A_n\bigl(|y'|^2+(1+t)^2\bigr)^{-(n-2)/2},
\qquad
A_n^2=\frac{2^{n-2}\Gamma(n/2)}{(n-2)\pi^{n/2}}.
\]
This choice gives $\|\nabla U_+\|_2=1$ and $\|U_+(\cdot,0)\|_q^2=S_\ast^{-1}$. Put $t_\ast=S_\ast^{1/2}$ and
\begin{equation}\label{B-recall-I-eq:mass-one-Uplus}
U_\ast=t_\ast U_+.
\end{equation}
Then $\|U_\ast(\cdot,0)\|_q=1$, $\|\nabla U_\ast\|_2^2=S_\ast$, and
\[
-\Delta U_\ast=0,
\qquad -\partial_tU_\ast=S_\ast U_\ast^{q-1}\quad\text{on }\partial\R^n_+.
\]
For $U_+$ the boundary equation coefficient is instead $\lambda_{\mathrm{flat}}=S_\ast^{q/2}$. All profile moments below use this fixed $U_+$.

\begin{convention}[Cutoff function]\label{II-conv:cutoff-choice}
Fix $\eta\in C_c^\infty([0,\infty);[0,1])$ equal to $1$ on $[0,1]$ and $0$ on $[2,\infty)$. Set $\chi_R(y)=\eta(|y|/R)^{n-2}$. Since $q-2=2/(n-2)$,
\[
\chi_R^{q-2}=\eta_R^2,\qquad
\chi_R^{q-1}=\eta_R^n,\qquad
\chi_R^q=\eta_R^{2(n-1)},\qquad \eta_R(y)=\eta(|y|/R).
\]
Put $\gamma=q-2$, $\alpha=(n-2)/2$, and $R=\varepsilon^{-\beta}$. We use logarithmic scale derivatives $\mathbf D=\varepsilon\partial_\varepsilon$. For $\rho=|Y|/R$, differentiation at fixed model coordinates gives
\[
\mathbf D(\chi_R^\gamma)=\beta\rho(\eta^2)'(\rho),
\]
whereas at fixed coordinates on $M$ the annular factor is $(\beta-1)\rho(\eta^2)'(\rho)$. Their first two logarithmic derivatives are bounded. After multiplication by $U_+^\gamma$, each annular coefficient has $L^{n-1}(\partial\R^n_+)$ norm $O(R^{-1})$. We call $\mathcal U^\gamma$ a boundary potential. The estimate
\begin{equation}\label{I-eq:boundary-potential-form-bound}
\left|\int_{\partial M}Vh\phi\right|
\le\|V\|_{n-1}\|h\|_q\|\phi\|_q
\le c_Y^{-1}C_{\mathrm{tr}}^2\|V\|_{n-1}\|h\|_E\|\phi\|_E
\end{equation}
explains its action on functions in the energy space. Ordinary $\partial_\varepsilon$ derivatives introduce an additional factor $\varepsilon^{-1}$.
\end{convention}

Let $\Phi_x$ be a Fermi chart centered at $x\in\partial M$. For $w$ supported in $B_{2R}^+$ and $2\varepsilon R<r_0$, define its rescaling in Fermi coordinates by
\[
(\mathcal T_{x,\varepsilon}w)(z)
=\varepsilon^{-(n-2)/2}w\bigl(\Phi_x^{-1}(z)/\varepsilon\bigr)
\quad\text{in the chart},
\]
and extend it by zero outside. For general $w$ we first insert $\chi_R$.

\begin{remark}[Cutoff radius]\label{B-recall-I-rem:bubble-cutoff}
Unless another exponent is specified, we take $R(\varepsilon)=\varepsilon^{-3/4}$. The support radius on $M$ is $2\varepsilon R(\varepsilon)=2\varepsilon^{1/4}$ which tends to zero as $\varepsilon\to0$.
\end{remark}

\begin{definition}[Bubbles and their parameter families]\label{B-recall-I-def:bubbles}
The boundary bubble and its boundary normalization are
\[
U_{+,x,\varepsilon}=\mathcal T_{x,\varepsilon}(\chi_{R(\varepsilon)}U_+),
\qquad
\mathcal U_{x,\varepsilon}=\frac{U_{+,x,\varepsilon}}{\|U_{+,x,\varepsilon}\|_{L^q(\partial M)}}.
\]
Allowing positive scalar multiples gives
\[
\mathcal C_\partial=\{aU_{+,x,\varepsilon}:a>0\}.
\]
Only scales whose supports lie in the stated charts are included. A concentrating subfamily will be specified in each modulation statement. For $\alpha\in(0,1)$, $\mathcal U^{(\alpha)}_{x,\varepsilon}$ denotes the same normalized construction with $R(\varepsilon)=\varepsilon^{-\alpha}$.
\end{definition}

\begin{convention}[Coefficients and diagonal limits]\label{B-recall-I-conv:cutoff-limits}
The geometric coefficients are obtained by expanding at fixed $R$ and then letting $R\to\infty$. Under the stated dimension assumptions, the required moments converge. We evaluate the expansions along $R(\varepsilon)=\varepsilon^{-3/4}$ unless otherwise specified; the stated remainders include the resulting truncation errors. The same choice is admissible for finitely many scales $k\varepsilon$, $k=1,2,3$. Fixed-radius statements display their dependence on $R$ explicitly.
\end{convention}

The modulation directions are $\partial_\varepsilon U_{+,x,\varepsilon}$ and $\partial_{x^a}U_{+,x,\varepsilon}$. We impose energy orthogonality to them, together with the boundary constraint described below.

Higher-order coordinate remainders use the stated \textup{(BG)} assumptions.

\subsection{The boundary constraint}
On $\mathcal S=\{u\in H^1(M):\|u\|_{L^q(\partial M)}=1\}$, the tangent space and normalization map are
\[
T_u\mathcal S=\left\{h:\int_{\partial M}|u|^{q-2}uh\,d\sigma_g=0\right\},
\qquad
\mathsf{ret}_u(h)=\frac{u+h}{\|u+h\|_{L^q(\partial M)}}.
\]
The constrained gradient is the energy Riesz representative of $D(\mathcal J_g|_{\mathcal S})[u][h]=2\mathcal N_g^\circ(u,h)$. Thus a Palais--Smale sequence at level $c$ satisfies $\mathcal J_g[u_k]\to c$ and $\|\mathsf{grad}_{\mathcal S}\mathcal J_g(u_k)\|_{E,g}\to0$. We use $H^{-1}_E(M)$ for the energy-norm dual; boundary distributions act through the trace. For $n\ge5$, put $\gamma=q-2=2/(n-2)\in(0,1)$.

\subsection{Standard-bubble expansions}
The standard-bubble expansion supplies the geometric terms in the reduction.

\begin{proposition}[Standard-bubble expansions]
\label{I-prop:standard-bubble-inputs}
Let $\mathcal C\subset\partial M$ be an open boundary region with a fixed Fermi collar, and suppose that $H_g=0$ on $\mathcal C$. Use the normalized bubbles of Definition~\ref{B-recall-I-def:bubbles}.
\begin{enumerate}[label=\textup{(\roman*)},leftmargin=2em]
\item Assume $n\ge5$ and \textup{(BG$^{3+}$)}, and use the cutoff $R(\varepsilon)=\varepsilon^{-3/4}$. Put $\kappa_3^{\mathrm{bare}}(n)=(6-n)/[2(n-1)(n-3)(n-4)]$. Uniformly for $x$ in compact subsets of $\mathcal C$,
\begin{equation}\label{I-eq:standard-bubble-quadratic-input}
\begin{aligned}
&\frac{\mathcal J_g[\mathcal U_{x,\varepsilon}]}{S_\ast}
=1+\kappa_3^{\mathrm{bare}}(n)|\mathring{\mathrm{II}}_g(x)|^2\varepsilon^2
+e_2(x,\varepsilon),\\&\qquad
|e_2(x,\varepsilon)|+|\varepsilon\partial_\varepsilon e_2(x,\varepsilon)|
=o(\varepsilon^2).
\end{aligned}
\end{equation}
In particular,
\begin{equation}\label{I-eq:standard-bubble-scale-input}
\frac1{S_\ast}\partial_\varepsilon\mathcal J_g[\mathcal U_{x,\varepsilon}]
=2\kappa_3^{\mathrm{bare}}(n)|\mathring{\mathrm{II}}_g(x)|^2\varepsilon
+o(\varepsilon).
\end{equation}
\item Assume $n\ge6$ and \textup{(BG$^{4+}$)}. At an umbilic point $p\in\mathcal C$, both the quadratic coefficient and the geometric cubic coefficient $\Theta_g(p)$ vanish. Here $\Theta_g$ is the cubic coefficient of $\mathcal J_g/S_\ast$, obtained from the fixed-$R$ expansion followed by $R\to\infty$. Consequently,
\[
\mathcal J_g[\mathcal U_{p,\varepsilon}]=S_\ast+o(\varepsilon^3)
\quad(n\ge7).
\]
For $n=6$, the same conclusion holds for $\mathcal U^{(\alpha)}_{p,\varepsilon}$ with each fixed $3/4<\alpha<1$.
\end{enumerate}
\end{proposition}

\begin{proof}
The quadratic coefficient is \cite[Proposition~4.7]{EscobarGeometry}. Put $\mathbf D=\varepsilon\partial_\varepsilon$. In the rescaled Fermi integrals the first-order term vanishes by $H_g=0$ and angular integration. The flat tail, the second-moment tail, and the Taylor error, including one $\mathbf D$ derivative, are bounded by
\[
C\left(R^{2-n}+\varepsilon^2R^{4-n}
+\varepsilon^3 L_n(R)\right),
\]
where $L_5(R)=\log(2R)$ and $L_n(R)=1$ for $n\ge6$. The radial derivative of the metric remainder has the same order under \textup{(BG$^{3+}$)}; a derivative of $\chi_R$ is supported on the same annulus and has the same decay. Scalar and boundary-density errors obey these bounds as well. The boundary norm is bounded away from zero. Thus division preserves the estimates, which are $o(\varepsilon^2)$ for $R=\varepsilon^{-3/4}$. This proves \eqref{I-eq:standard-bubble-quadratic-input} and its ordinary scale-derivative consequence.

Every cubic contraction contains an undifferentiated $\mathring{\mathrm{II}}$ by \cite[Equation~(B.21), Corollary~B.10]{EscobarGeometry}. The quadratic and cubic coefficients therefore vanish at $p$. The expansion of \cite[Lemma~4.11 and Equations~(4.25), (4.27)]{EscobarGeometry} leaves the bound
\[
C\left(R^{2-n}+\varepsilon^2R^{4-n}
+\varepsilon^3R^{5-n}+\varepsilon^4 K_n(R)\right),
\]
where $K_6(R)=\log(2R)$ and $K_n(R)=1$ for $n\ge7$. This is $o(\varepsilon^3)$ for the cutoff $R(\varepsilon)=\varepsilon^{-3/4}$ when $n\ge7$. For $n=6$, take $R=\varepsilon^{-\alpha}$ with $\alpha>3/4$, so the flat tail $\varepsilon^{4\alpha}$ is also $o(\varepsilon^3)$.
\end{proof}

\section{Coercivity and modulation}
\label{I-sec:coercivity-first}

\subsection{Coercivity on the hemisphere and conformal invariance}

The one-bubble reduction rests on a spectral gap transverse to the optimizer orbit and a uniform choice of bubble parameters. The gap is the linearized estimate underlying half-space trace stability: Ho's theorem \cite{Ho2022} is recalled in \cite[Section~2]{BorquezCajuBosch2026}, and a refined spectral estimate appears in \cite[Proposition~7.1]{ZhangZhouZou2025}. Supplementary deficit estimates are given in Appendix~\ref{A-sec:supplementary-stability}.

Fix the Cayley transform $\mathcal C:S^n_+\setminus\{p_\infty\}\to\R^n_+$, where $p_\infty\in\partial S^n_+$. Its conformal pullback $\mathcal C^\ast U_+$ extends smoothly across $p_\infty$. Let $F_\ast$ be its normalization to unit boundary $L^q$ norm, let $\mathcal M_S$ be its positive normalized optimizer orbit under boundary-preserving conformal maps, and let $E=T_{F_\ast}\mathcal M_S$. On the hemisphere the constraint is
\begin{equation}\label{I-eq:Esc-constraint-slice}
\mathcal S=\{u\in H^1(S^n_+):\|u\|_{L^q(\partial S^n_+)}=1\}.
\end{equation}

\begin{lemma}[Equivalence and conformal invariance of the energy norm]\label{I-lem:graph-norm}
The energy norm satisfies
\begin{equation}\label{I-eq:graph-norm-explicit}
\|u\|_{E,g}^2=\mathcal N_g^\circ(u).
\end{equation}
It is equivalent to the $H^1$ norm on the round hemisphere, and on $(M,g)$ under \textup{(BG$^{2+}$)} and \textup{(Y$^+_\partial$)}. If $\widehat g=\phi^{4/(n-2)}g$, multiplication by $\phi$ preserves the energy norm and the Escobar deficit. In particular, the Cayley pullback is an isometry from $\dot H^1(\R^n_+)$ to the hemisphere energy space and preserves the deficit and energy distances.
\end{lemma}

\begin{proof}
Green's identity gives
\begin{equation}\label{I-eq:green-identity}
\int_MuL_g^\circ v\,dV_g+\int_{\partial M}uB_g^\circ v\,d\sigma_g
=\mathcal N_g^\circ(u,v).
\end{equation}
On the round hemisphere, $\Scal_g>0$ and $H_g=0$, so the energy and $H^1$ norms are equivalent. On $M$, the bounds follow from \eqref{B-recall-I-eq:Yplus} and \eqref{B-recall-I-eq:Esc-energy-continuity}. Conformal covariance gives
\[
\mathcal N_{\widehat g}^\circ(v,w)=\mathcal N_g^\circ(\phi v,\phi w),
\qquad
\int_{\partial M}|v|^q\,d\sigma_{\widehat g}
=\int_{\partial M}|\phi v|^q\,d\sigma_g;
\]
see \cite[Lemma~3.13]{EscobarGeometry}. For the Cayley map these identities first hold for functions vanishing near the pole. Such functions are dense because a cutoff at radius $\rho$ has squared Dirichlet norm $O(\rho^{n-2})$ and squared $L^2$ norm $O(\rho^n)$, with $n\ge3$. Density in the completed half-space Dirichlet space then proves the assertion.
\end{proof}

\begin{lemma}[Constrained second variation at $F_\ast$]\label{I-lem:hessian-sphere}\label{I-prop:robin-gap}
For the normalized optimizer $F_\ast$, one has
\[
L_{g_{\mathrm{round}}}^\circ F_\ast=0\ \text{ in }S^n_+,
\qquad
B_{g_{\mathrm{round}}}^\circ F_\ast=\lambda_\ast F_\ast^{q-1}\ \text{ on }\partial S^n_+,
\]
with $\lambda_\ast=\mathcal N^\circ_{g_{\mathrm{round}}}(F_\ast)=\mathcal J_{g_{\mathrm{round}}}[F_\ast]=S_\ast$. For $\phi\in T_{F_\ast}\mathcal S=\{\phi:\int_{\partial S^n_+}F_\ast^{q-1}\phi=0\}$, the constrained Jacobi form, equal to half the Hessian, is
\begin{equation}\label{I-eq:Hessian-sphere}
Q[\phi]
= \mathcal N^\circ_{g_{\mathrm{round}}}(\phi,\phi)-(q-1)\lambda_\ast\int_{\partial S^n_+} F_\ast^{q-2}\phi^2d\sigma.
\end{equation}
We denote by
\[
Q[\phi,\eta]\coloneqq
\mathcal N^\circ_{g_{\mathrm{round}}}(\phi,\eta)
-(q-1)\lambda_\ast
\int_{\partial S^n_+}F_\ast^{q-2}\phi\eta\,d\sigma
\]
the associated symmetric bilinear form, so that $Q[\phi]=Q[\phi,\phi]$. Since $H_{g_{\mathrm{round}}}=0$, the boundary mean-curvature term is absent. Moreover, $\ker Q=E$ on $T_{F_\ast}\mathcal S$, and there exists $c_\ast>0$ such that $Q[\phi]\ge c_\ast\|\phi\|_{H^1(S^n_+)}^2$ for all $\phi\in T_{F_\ast}\mathcal S\cap E^\perp$.

The coercivity constant relative to the energy norm can be chosen as $2/(n+2)$, uniformly along $\mathcal M_S$. In particular,
\[
L_{g_{\mathrm{round}}}^\circ\phi=0,
\qquad B_{g_{\mathrm{round}}}^\circ\phi=(q-1)\lambda_\ast F_\ast^{q-2}\phi
\]
has no nonzero solution in $T_{F_\ast}\mathcal S\cap E^\perp$.
\end{lemma}

\begin{proof}
Write $\mathcal B(u,v)=\mathcal N^\circ_{g_{\mathrm{round}}}(u,v)$ and $\mathcal G(u)=\int_{\partial S^n_+}|u|^q\,d\sigma$. On $\mathcal G=1$, the quotient equals $\mathcal B(u,u)$. The Euler--Lagrange equation is
\begin{equation}\label{I-eq:EL-sphere-weak}
\mathcal B(F_\ast,\eta)
=\lambda_\ast\int_{\partial S^n_+}F_\ast^{q-1}\eta d\sigma
\qquad\forall\eta\in H^1(S^n_+).
\end{equation}
Testing with $F_\ast$ gives $\lambda_\ast=\mathcal B(F_\ast,F_\ast)=S_\ast$.

Differentiating $\int|u_t|^q=1$ twice along a constrained curve with $\dot u_0=\phi$ gives
\[
\int_{\partial S^n_+}F_\ast^{q-1}\ddot u_0
=-(q-1)\int_{\partial S^n_+}F_\ast^{q-2}\phi^2.
\]
Substitution in the second derivative of $\mathcal B(u_t,u_t)$ proves \eqref{I-eq:Hessian-sphere}.

For $Z\in E$, choose a curve of normalized optimizers $F_s$ with $F_0=F_\ast$ and $\dot F_0=Z$. Each satisfies \eqref{I-eq:EL-sphere-weak} with the same multiplier $S_\ast$. Differentiation gives $Q[Z,\eta]=0$ for every $\eta\in H^1(S^n_+)$, so $E\subset\ker Q$. For the reverse inclusion and the gap, use the conformal map from the half-space to the Euclidean unit ball $B^n$ whose conformal factor is $\bigl(2/(|y'|^2+(1+y_n)^2)\bigr)^{(n-2)/2}$. It carries $U_+$ to a constant, so the corresponding map from the hemisphere carries $F_\ast$ to the boundary-normalized constant optimizer on $B^n$. Conformal covariance preserves the energy, constraint, and Jacobi form. Put $\alpha=(n-2)/2$. On the ball,
\[
\|\phi\|_{E,B}^2
=\int_{B^n}|\nabla\phi|^2
+\alpha\int_{\partial B^n}\phi^2,
\qquad
Q_B[\phi]=\int_{B^n}|\nabla\phi|^2
-\int_{\partial B^n}\phi^2,
\]
because the constant optimizer satisfies $\lambda_\ast F_\ast^{q-2}=\alpha$ and $\alpha(q-2)=1$.

Write $\phi=h+w$, where $h$ is the harmonic extension of its trace and $w\in H^1_0(B^n)$. Let $h_\ell$ denote the trace component in spherical harmonics of degree $\ell$. For $H_\ell(r,\theta)=r^\ell h_\ell(\theta)$, Green's identity, with $\nu_{\mathrm{out}}$ the outward unit normal to the ball, gives
\[
\int_{B^n}\nabla H_\ell\cdot\nabla H_m
=\int_{\partial B^n}h_\ell\partial_{\nu_{\mathrm{out}}} H_m
=\delta_{\ell m}\ell\|h_\ell\|_2^2,
\qquad
\int_{B^n}\nabla H_\ell\cdot\nabla w=0.
\]
Together with $\|\operatorname{Tr}\phi\|_2^2=\sum_\ell\|h_\ell\|_2^2$, this yields
\begin{align*}
\|\phi\|_{E,B}^2
&=\|\nabla w\|_{L^2(B^n)}^2
+\sum_{\ell\ge0}(\ell+\alpha)
\|h_\ell\|_{L^2(\partial B^n)}^2,\\
Q_B[\phi]
&=\|\nabla w\|_{L^2(B^n)}^2
+\sum_{\ell\ge0}(\ell-1)
\|h_\ell\|_{L^2(\partial B^n)}^2.
\end{align*}
The tangency condition requires the trace to have zero mean, so $h_0=0$. The kernel is exactly the degree-one harmonic space, generated by the coordinate functions; these are the derivatives of the boundary-preserving conformal motions of the constant optimizer. Under the conformal change of variables, this space gives $E$. Furthermore, $\langle\phi,x_i\rangle_{E,B}=(1+\alpha)
\int_{\partial B^n}\phi x_i$, so orthogonality to this kernel in the energy inner product
also removes $h_1$. Since
\[
\frac{\ell-1}{\ell+\alpha}\ge\frac1{2+\alpha}
=\frac2{n+2}\qquad(\ell\ge2),
\]
one obtains $Q_B[\phi]\ge\frac2{n+2}\|\phi\|_{E,B}^2$ on the constrained complement. Applying the inverse conformal change and using equivalence of the energy and $H^1$ norms proves the stated kernel identity and coercivity. Testing the Jacobi system with $\phi$ gives uniqueness on the complement; conformal covariance gives uniformity along the optimizer orbit.
\end{proof}

\begin{remark}[Modulation on a general manifold]\label{I-rem:modulation-general}
We choose the center, scale, and amplitude by the boundary constraint and energy orthogonality. The corresponding statement on the hemisphere is Lemma~\ref{I-lem:modulation-sphere}.
\end{remark}

\begin{lemma}[Modulation near $\mathcal C_\partial$]\label{I-lem:modulation}
Assume $n\ge5$, \textup{(Y$^+_{\partial}$)}, and \textup{(BG$^{2+}$)}, and fix a Fermi collar. There are sufficiently small constants $\varepsilon_1,\delta_0,\delta_1>0$, with $\delta_0<\delta_1$, for the following assertion. Define the concentrating family
\[
\mathcal C_\partial
=\big\{aU_{+,x,\varepsilon}:a>0,\ x\in\partial M,
\ 0<\varepsilon\le\varepsilon_1/2\big\}.
\]
If $u\in H^1(M)\setminus\{0\}$ satisfies $\operatorname{dist}_{H^1}(u,\mathcal C_\partial)
\le\delta_0\|u\|_{H^1(M)}$, there are parameters
$(a^\sharp,x^\sharp,\varepsilon^\sharp)$ with $a^\sharp>0$ and $0<\varepsilon^\sharp\le\varepsilon_1$ such that, with $v^\sharp=a^\sharp U_{+,x^\sharp,\varepsilon^\sharp}$ and $r=u-v^\sharp$,
\[
\langle r,\partial_{\varepsilon}v^\sharp\rangle_{E,g}=0,
\qquad
\langle r,\partial_{x^\alpha}v^\sharp\rangle_{E,g}=0\quad(\alpha=1,\dots,n-1).
\]
The remaining orthogonality condition is $
\int_{\partial M}(v^\sharp)^{q-1}r\,d\sigma_g=0$.
This solution is unique among triples satisfying these orthogonality conditions, $0<\varepsilon\le\varepsilon_1$, and $\|u-aU_{+,x,\varepsilon}\|_{H^1}\le\delta_1\|u\|_{H^1}$. The parameter map is $C^1$, and
\[
\|r\|_{H^1(M)}\le C\operatorname{dist}_{H^1}(u,\mathcal C_\partial).
\]
The modulation constants depend on the fixed metric and, in the finite-regularity case, on its $C^3$ bounds.
\end{lemma}

\begin{proof}
Write $\Psi(a,\xi,\varepsilon)=aU_{+,x(\xi),\varepsilon}$ in a boundary chart. The equations are
\[
\int_{\partial M}\Psi^{q-1}(u-\Psi)=0,
\qquad \langle u-\Psi,\partial_\varepsilon\Psi\rangle_E=0,
\qquad \langle u-\Psi,\partial_{\xi^a}\Psi\rangle_E=0.
\]
Amplitude homogeneity permits $\|u\|_{H^1}=1$; the relevant amplitudes then lie in a fixed compact subset of $(0,\infty)$. Near $(a_0,\xi_0,\varepsilon_0)$ use parameters
$\bigl(a,(\xi-\xi_0)/\varepsilon_0,\log(\varepsilon/\varepsilon_0)\bigr)$ and multiply the center and scale equations by $\varepsilon_0$ and $\varepsilon$, respectively. The profile has uniformly bounded first two derivatives in these parameters. The case of $C^3$ metrics follows from Lemma~\ref{I-lem:finite-regularity-preparation}. In the flat limit the Jacobian has a nonzero amplitude entry and Jacobi block $-a_0^2(n-2)I/4$; the mixed constraint entries vanish because translation and dilation preserve the boundary mass. Fermi rescaling therefore makes the Jacobian uniformly invertible. Terms involving second profile derivatives are multiplied by $u-\Psi$. The quantitative inverse function theorem gives a unique $C^1$ solution and a uniform local distance estimate. Choosing the logarithmic scale radius smaller than $\log2$ ensures $\varepsilon^\sharp\le\varepsilon_1$.

For uniqueness throughout the stated neighborhood, when the relative scale and center parameters leave every compact set, the flat energy inner product of the corresponding bubbles tends to zero. Fermi rescaling and cutoff tails preserve the resulting positive separation. On a compact set of relative parameters, continuity and injectivity give positive separation outside the local inversion neighborhood. Thus sufficiently small $\varepsilon_1$ and $\delta_1$ place any two candidate triples in one such neighborhood, where they coincide.

Apply the uniform local estimate at approximate nearest bubbles to obtain
$\|u-v^\sharp\|_{H^1}\le C\operatorname{dist}_{H^1}(u,\mathcal C_\partial)$.
Choose $C\delta_0<\delta_1$ and rescale the amplitude.
\end{proof}

\begin{remark}[Amplitude under the boundary normalization]\label{I-rem:amplitude-bound}
For $u\in\mathcal S$, write $w=u-v^\sharp$. The boundary orthogonality and Taylor expansion give
\[
a^\sharp\|U_{+,x^\sharp,\varepsilon^\sharp}\|_q=1+O(\|w\|_{H^1}^2),
\qquad
a^\sharp=S_\ast^{1/2}+O(\varepsilon^\sharp)+O(\|w\|_{H^1}^2).
\]
Thus $a^\sharp\asymp S_\ast^{1/2}$ uniformly near the concentrating family. For a boundary bubble $U$, the homogeneity identities
\[
D\mathcal J_g[a^\sharp U][w]=(a^\sharp)^{-1}D\mathcal J_g[U][w],
\qquad
D^2\mathcal J_g[a^\sharp U][w,w]=(a^\sharp)^{-2}D^2\mathcal J_g[U][w,w]
\]
give uniform constants in the first-variation and coercivity estimates.
\end{remark}

\subsection{Coercivity near a boundary bubble on \texorpdfstring{$(M,g)$}{(M, g)}}

\begin{theorem}[Coercivity near a boundary bubble on $(M,g)$]\label{I-thm:intrinsic-single-gap}
Assume $n\ge5$, \textup{(Y$^+_{\partial}$)}, and \textup{(BG$^{2+}$)}. Use the energy inner product of \eqref{B-recall-I-eq:Esc-energy-bilinear}. Then there exist $\varepsilon_0>0$ and $c_\ast>0$, depending only on $n$ and the \textup{(BG$^{2+}$)} and \textup{(Y$^+_{\partial}$)} data of $(M,g)$, such that for every $x\in\partial M$, every $\varepsilon\in(0,\varepsilon_0]$, and every $w\in H^1(M)$ satisfying
\begin{equation}\label{I-eq:esc-orthog-gauge-final}
\langle w,\partial_\varepsilon U_{+,x,\varepsilon}\rangle_{E,g}
=\langle w,\nabla_x U_{+,x,\varepsilon}\rangle_{E,g}=0,
\qquad
\int_{\partial M}(U_{+,x,\varepsilon})^{q-1}w\,d\sigma_g=0,
\quad q=\tfrac{2(n-1)}{n-2},
\end{equation}
one has
\begin{equation}\label{I-eq:intrinsic-gap-final}
\big\langle D^2\mathcal J_g[U_{+,x,\varepsilon}]w,w\big\rangle
\ge c_\ast\|w\|_{E,g}^2
\asymp\|w\|_{H^1(M)}^2.
\end{equation}
\end{theorem}

\begin{proof}
Suppose the uniform gap fails. Then there are $x_k\in\partial M$, $\varepsilon_k\downarrow0$, and $w_k$ such that
\begin{equation}\label{I-eq:contr-setup}
\|w_k\|_{E,g}=1,\qquad
\eqref{I-eq:esc-orthog-gauge-final}\text{ holds at }(x_k,\varepsilon_k),
\qquad D^2\mathcal J_g[U_{+,x_k,\varepsilon_k}][w_k,w_k]\le k^{-1}.
\end{equation}
Write $w_k\rightharpoonup W$ in $H^1(M)$ and $v_k=w_k-W$. Strong $L^2$ convergence on $M$ and $\partial M$ implies
\[
1=\|W\|_E^2+\|v_k\|_E^2+o(1),
\qquad \|v_k\|_E^2=\int_M|\nabla v_k|^2+o(1).
\]
For the critical rescaling $\widetilde w_k(y)=\varepsilon_k^{(n-2)/2}w_k(\Phi_{x_k}(\varepsilon_ky))$, the fixed function $W$ has vanishing gradient energy on each shrinking ball in $M$. Thus the coordinate bounds give, for fixed $R$,
\begin{equation}\label{I-eq:per-ball-bound}
\int_{B_R^+}|\nabla\widetilde w_k|^2
\le(1+O(\varepsilon_kR))\int_M|\nabla v_k|^2+o(1)
\le1+o(1).
\end{equation}
On the full expanding chart, uniform coordinate comparison gives
\begin{equation}\label{I-eq:global-grad-bound}
\int_{B_{r_0/\varepsilon_k}^+}|\nabla\widetilde w_k|^2
\le C\|w_k\|_{H^1(M)}^2\le C.
\end{equation}
The global Sobolev inequality also bounds $\|\widetilde w_k\|_{2n/(n-2)}$. The trace inequality on the expanding half-ball gives
\begin{equation}\label{I-eq:global-trace-bound}
\|\widetilde w_k(\cdot,0)\|_{L^q(|y'|<r_0/\varepsilon_k)}
\le C\left(\|\nabla\widetilde w_k\|_2+
\frac{\varepsilon_k}{r_0}\|\widetilde w_k\|_2\right)\le C,
\end{equation}
because $\varepsilon_k\|\widetilde w_k\|_2\le C\|w_k\|_{L^2(M)}$. A diagonal subsequence converges weakly in local $H^1$ and strongly in local boundary $L^2$ to $\widetilde w$. By \eqref{I-eq:per-ball-bound},
\begin{equation}\label{I-eq:wlimit-gradient-bound}
\int_{\R^n_+}|\nabla\widetilde w|^2\le1.
\end{equation}
The bound in $L^{2n/(n-2)}$ gives $\widetilde w\in L^{2n/(n-2)}$, and hence $\widetilde w\in\mathcal D^{1,2}(\R^n_+)$.

For tangent directions, differentiation of the quotient gives
\[
D^2\mathcal J_g[U_k][w_k,w_k]
=2\|U_k\|_q^{-2}\left(1-(q-1)\lambda_k
\int_{\partial M}U_k^{q-2}w_k^2\right),
\quad
\lambda_k=\frac{\mathcal N_g^\circ(U_k)}{\|U_k\|_q^q}\to\lambda_{\mathrm{flat}},
\]
where $U_k=U_{+,x_k,\varepsilon_k}$. The weighted integrals are bounded by \eqref{I-eq:global-trace-bound}, so a further subsequence satisfies
\begin{equation}\label{I-eq:bdry-mass}
(q-1)\lambda_{\mathrm{flat}}\int_{\partial M}U_k^{q-2}w_k^2\longrightarrow m\ge1.
\end{equation}
Local trace convergence passes this integral to the limit. Indeed, $U_+^{q-2}(\cdot,0)\in L^{n-1}$ and
\[
\int_{|y'|>R}U_+^{q-2}\widetilde w_k^2
\le C\|U_+^{q-2}\|_{L^{n-1}(|y'|>R)}\longrightarrow0
\]
uniformly in $k$. The same tail estimates hold after insertion of the cutoff and the boundary Jacobian, both of which converge to $1$ on fixed balls. Define
\[
Q_{\mathrm{flat}}[v]=\int_{\R^n_+}|\nabla v|^2
-(q-1)\lambda_{\mathrm{flat}}\int_{\partial\R^n_+}U_+^{q-2}v^2.
\]
Then
\[
(q-1)\lambda_{\mathrm{flat}}\int_{\partial\R^n_+}U_+^{q-2}\widetilde w^2=m,
\qquad Q_{\mathrm{flat}}[\widetilde w]\le1-m\le0.
\]
In particular $\widetilde w\ne0$. The boundary constraint passes to the limit by the critical trace tail estimate. The normalized scale and center orthogonality conditions pass by local weak convergence, the Dirichlet tails of the flat Jacobi fields, and the vanishing lower-order terms. Thus $\widetilde w$ is constrained and orthogonal to every flat translation and dilation field. Conformal invariance and Lemma~\ref{I-prop:robin-gap} give $Q_{\mathrm{flat}}[\widetilde w]>0$, a contradiction.
\end{proof}

\begin{lemma}[Estimate for the first variation]\label{I-lem:first-order-drift}
Assume $n\ge5$, \textup{(BG$^{2+}$)}, and \textup{(Y$^+_{\partial}$)}. There exist $\varepsilon_0>0$ and $C_{\mathrm{geo}}>0$, depending only on $n$ and the \textup{(BG$^{2+}$)} and \textup{(Y$^+_{\partial}$)} data, such that the following holds.

For every $x\in\partial M$, every $\varepsilon\in(0,\varepsilon_0]$, and every $a>0$, set
\[
v_{x,\varepsilon}\coloneqq aU_{+,x,\varepsilon},
\]
where $U_{+,x,\varepsilon}=\mathcal T_{x,\varepsilon}(\chi_{R(\varepsilon)}U_+)$ is the boundary bubble of Definition~\ref{B-recall-I-def:bubbles}, with the cutoff chosen as in Remark~\ref{B-recall-I-rem:bubble-cutoff}. Then
\[
\big|L[v_{x,\varepsilon}](\psi)\big|
\coloneqq\big|D\mathcal J_g[v_{x,\varepsilon}](\psi)\big|
\le C_{\mathrm{geo}}\varepsilon a^{-1}\|\psi\|_{E,g}
\]
for every $\psi\in T_{v_{x,\varepsilon}}\mathcal S_{v_{x,\varepsilon}}$, where
\[
\mathcal S_{v_{x,\varepsilon}}\coloneqq\Big\{F\in H^1(M):\ \|F\|_{L^q(\partial M)}=\|v_{x,\varepsilon}\|_{L^q(\partial M)}\Big\},
\qquad q=2^*_\partial.
\]
\end{lemma}

\begin{proof}
Write $v=v_{x,\varepsilon}$ and $q=2^*_\partial$. For $\psi\in T_v\mathcal S_v$, the first variation is
\begin{equation}\label{I-eq:L-as-N}
L[v](\psi)=\frac{2}{\|v\|_{L^q(\partial M)}^2}\mathcal N_g^\circ(v,\psi),
\end{equation}
since the derivative of the denominator vanishes on $T_v\mathcal S_v$. By the Green identity,
$$
\mathcal N_g^\circ(v,\psi)=\int_M\psi L_g^\circ v\,dV_g+\int_{\partial M}\psi B_g^\circ v\,d\sigma_g.
$$

It suffices by homogeneity to take $a=1$: $D\mathcal J_g[aU](\psi)=a^{-1}D\mathcal J_g[U](\psi)$, and the tangent condition is unchanged by this amplitude. Write $\widetilde\psi$ for the critically rescaled test function. After cutoff in the fixed collar, its Dirichlet and critical trace norms are bounded by $C\|\psi\|_{E,g}$. Write $P_R=\chi_RU_+$. Testing the flat equation by $\chi_R\widetilde\psi$ gives the exact cancellation
\begin{equation}\label{I-eq:cutoff-first-variation-identity}
\begin{aligned}
&\int\nabla P_R\cdot\nabla\widetilde\psi
-\lambda_{\mathrm{flat}}\int_{\partial\R^n_+}P_R^{q-1}\widetilde\psi\\
&\quad=\int U_+\nabla\chi_R\cdot\nabla\widetilde\psi
-\int\widetilde\psi\nabla U_+\cdot\nabla\chi_R
+\lambda_{\mathrm{flat}}\int_{\partial\R^n_+}
(\chi_R-\chi_R^{q-1})U_+^{q-1}\widetilde\psi.
\end{aligned}
\end{equation}
The integrals on the right are supported on $R<|Y|<2R$. Cauchy--Schwarz, Hardy's inequality, and the trace inequality bound their sum by $CR^{-\alpha}\|\widetilde\psi\|_{\mathcal D^{1,2}}$. In the rescaled energy pairing, the principal coefficient error is $O(\varepsilon|Y|)$; its pairing with the bubble is bounded by
\[
C\varepsilon\bigl\||Y|\nabla U_+\bigr\|_{L^2}
\|\widetilde\psi\|_{\mathcal D^{1,2}}.
\]
The weighted norm is finite for $n\ge5$. The boundary mean-curvature term is bounded by $C\varepsilon\|U_+\|_{L^{q'}}\|\widetilde\psi\|_{L^q}$, $q'=q/(q-1)$; $U_+(\cdot,0)\in L^{q'}$ in the same dimension range. The scalar term is controlled by Hardy's inequality, with bound
\[
C\varepsilon^2\bigl\||Y|U_+\bigr\|_{L^2(B^+_{2R})}
\|\widetilde\psi\|_{\mathcal D^{1,2}}
\le C\varepsilon^2\ell'_n(R)\|\psi\|_{E,g},
\]
where $\ell'_5(R)=R^{1/2}$, $\ell'_6(R)=(1+\log R)^{1/2}$, and $\ell'_n(R)=1$ for $n\ge7$. The boundary density error is $O(\varepsilon^2)$, since $|Y'|^2U_+^{q-1}\in L^{q'}(\mathbb R^{n-1})$ for $n\ge5$. The flat boundary term vanishes by tangency, up to the boundary-density error just estimated. The annular terms in \eqref{I-eq:cutoff-first-variation-identity} are bounded by
\[
C\omega(R)\|\psi\|_{E,g},\qquad
\omega(R)=O(R^{-(n-2)/2}).
\]
For $R=\varepsilon^{-3/4}$ both $\varepsilon^2\ell'_n(R)$ and $\omega(R)$ are $o(\varepsilon)$. The quotient denominator is bounded above and below independently of the center and scale. Equation~\eqref{I-eq:L-as-N} now gives the claimed estimate, including the factor $a^{-1}$.\qedhere

\end{proof}

\section{Compactness and Palais--Smale decompositions}
\label{I-sec:compactness}

Palais--Smale sequences are precompact below $S_\ast$ and admit a finite profile decomposition at every level. The profile masses determine the first mixed level and characterize the noncompact levels between one and two bubbles. For the concentration-compactness and profile-decomposition background, see \cite{Lions85a,Lions85b,Gerard1998}.

\subsection{Profile decomposition without a sign assumption}\label{I-sec:signfree-compactness}

\begin{proposition}[Boundary profile decomposition without a sign assumption]
\label{I-prop:signfree-boundary-profile-decomposition}
Assume \textup{(Y$^+_\partial$)} and \textup{(BG$^{2+}$)}. Let $(u_k)\subset\mathcal S$ be a constrained Palais--Smale sequence at level $c>0$, with no sign assumption.
Then $c\ge C_{\Esc}^*(M,g)>0$.  After passing to a subsequence, there exist a weak limit $u_\infty\in H^1(M)$, an integer $J\ge0$, boundary points $x_k^j\in\partial M$, scales $\varepsilon_k^j\downarrow0$, nonzero profiles $V^{(j)}\in\dot H^1(\R^n_+)$, $1\le j\le J$, and cutoff radii $R_k^j\to\infty$ with $\varepsilon_k^jR_k^j\to0$, such that
\begin{equation}\label{I-eq:signfree-profile-decomp}
u_k
=
u_\infty
+
\sum_{j=1}^J
\mathcal T_{x_k^j,\varepsilon_k^j}(\chi_{R_k^j}V^{(j)})
+
r_k,
\qquad
\|r_k\|_{H^1(M)}\to0.
\end{equation}
Write $d_\partial$ for the distance of the induced boundary metric. The centers and scales satisfy the separation condition
\begin{equation}\label{I-eq:signfree-orthogonality}
\frac{\varepsilon_k^i}{\varepsilon_k^j}
+
\frac{\varepsilon_k^j}{\varepsilon_k^i}
+
\frac{d_\partial(x_k^i,x_k^j)^2}{\varepsilon_k^i\varepsilon_k^j}
\longrightarrow\infty
\qquad(i\ne j).
\end{equation}
The weak limit solves
\[
L_g^\circ u_\infty=0
\quad\text{in }M,
\qquad
B_g^\circ u_\infty
=
c|u_\infty|^{q-2}u_\infty
\quad\text{on }\partial M,
\]
and each profile solves the critical boundary problem
\[
-\Delta V^{(j)}=0\quad\text{in }\R^n_+,
\qquad
-\partial_{y_n}V^{(j)}
=
c|V^{(j)}|^{q-2}V^{(j)}
\quad\text{on }\partial\R^n_+.
\]

The boundary mass and energy split:
\begin{align}
1
&=
\int_{\partial M}|u_\infty|^q\,d\sigma_g
+
\sum_{j=1}^J m_j
+
o(1),
\label{I-eq:signfree-mass-split}\\
c
&=
\mathcal N_g^\circ(u_\infty)
+
\sum_{j=1}^J E_j
+
\mathcal N_g^\circ(r_k)
+
o(1),
\label{I-eq:signfree-energy-split}
\end{align}
where $m_j\coloneqq \int_{\partial\R^n_+}|V^{(j)}|^q\,dy'$ and $E_j\coloneqq \int_{\R^n_+}|\nabla V^{(j)}|^2\,dy$.  Testing the profile equations by their respective solutions gives
\[
\mathcal N_g^\circ(u_\infty)
=
c\int_{\partial M}|u_\infty|^q\,d\sigma_g,
\qquad
E_j=cm_j.
\]
Each nonzero profile has
\begin{equation}\label{I-eq:profile-mass-lower-bound}
m_j\ge\left(\frac{S_\ast}{c}\right)^{n-1},
\end{equation}
so only finitely many profiles can occur.
\end{proposition}

\begin{proof}
Since $C^*_{\Esc}(M,g)=\inf_{\operatorname{Tr}u\ne0}\mathcal J_g[u]$ and $\mathcal J_g[u_k]\to c$, the level satisfies $c\ge C^*_{\Esc}(M,g)>0$.

The profile extraction follows Struwe's method \cite{Struwe1984} and the boundary extraction lemma \cite[Lemma~2.4 and (2-11)]{Almaraz2014}. For a related Struwe decomposition for polyharmonic operators on compact Riemannian manifolds, with or without boundary, see \cite{Mazumdar2017}. After normalizing the boundary coefficient to one, the limiting equation has the weak form
\[
\int_{\R^n_+}\nabla V\cdot\nabla\phi
=\int_{\partial\R^n_+}|V|^{q-2}V\phi,
\quad\text{equivalently }-\partial_{y_n}V=|V|^{q-2}V.
\]
The extraction applies to sequences without a sign assumption and gives a nonzero profile. The corresponding flat decomposition is proved in \cite[Theorem~5.1]{ZhangZhouZou2025}.

\smallskip\noindent(i) \emph{Profile extraction.}
For $u\in\mathcal S$ set $\ell_u(\phi)=\int_{\partial M}|u|^{q-2}u\phi$ and $\psi=\phi-\ell_u(\phi)u$. Then $\psi\in T_u\mathcal S$ and, on energy-bounded sets, $\|\psi\|_E\le C\|\phi\|_E$. The exact identity
\begin{equation}\label{I-eq:exact-PS-multiplier}
\mathcal N_g^\circ(u,\phi)-\mathcal J_g[u]\ell_u(\phi)
=\mathcal N_g^\circ(u,\psi)
=\tfrac12D(\mathcal J_g|_{\mathcal S})(u)[\psi]
\end{equation}
therefore gives $\lambda_k=\mathcal J_g[u_k]\to c$ and
\[
\mathcal N_g^\circ(u_k,\phi)
-\lambda_k\int_{\partial M}|u_k|^{q-2}u_k\phi\,d\sigma_g
=o(1)\|\phi\|_{H^1(M)}.
\]
Since $\lambda_k\to c$, the sequence is Palais--Smale for
\[
\mathscr I_c(u)=\frac12\mathcal N_g^\circ(u)
-\frac cq\int_{\partial M}|u|^q\,d\sigma_g.
\]
After weak convergence and almost-everywhere trace convergence, the limit $u_\infty$ solves the stated equation. The absolute-value estimates in the proof of \cite[Lemma~2.2]{Almaraz2014} give the same derivative splitting for the nonlinearity $|u|^{q-2}u$. Hence $u_k^{(0)}=u_k-u_\infty$ is Palais--Smale for $\mathscr I_c$. The scalar-curvature and boundary mean-curvature forms are compact from $H^1(M)$ to its dual: for example,
\[
\left|\int_M\Scal_g u_k^{(0)}\phi\right|
+\left|\int_{\partial M}H_g u_k^{(0)}\phi\right|
\le C\left(\|u_k^{(0)}\|_{L^2(M)}
+\|u_k^{(0)}\|_{L^2(\partial M)}\right)\|\phi\|_{H^1}
=o(1)\|\phi\|_{H^1}.
\]
Thus $u_k^{(0)}=u_k-u_\infty$ converges weakly to zero and is a Palais--Smale sequence for
\[
\mathscr I_c^0(u)=\frac12\int_M|\nabla u|^2\,dV_g
-\frac cq\int_{\partial M}|u|^q\,d\sigma_g.
\]
Multiplication by $c^{1/(q-2)}$ reduces the boundary coefficient to one. We may therefore apply \cite[Lemma~2.4]{Almaraz2014} successively, reversing this fixed rescaling after each extraction. The extraction preserves the Palais--Smale condition after subtraction and splits the energy. Each remainder, rescaled at a previously extracted center and scale, converges weakly to zero; comparison of rescalings gives \eqref{I-eq:signfree-orthogonality}. The Dirichlet tails of each fixed profile permit a diagonal choice $R_k^j\to\infty$ with $\varepsilon_k^jR_k^j\to0$. The Brezis--Lieb lemma \cite[Theorem~1]{BrezisLieb1983} and energy orthogonality give \eqref{I-eq:signfree-mass-split}--\eqref{I-eq:signfree-energy-split}; the compact lower-order forms vanish on the rescaled profiles. An interior profile would be a finite-energy harmonic function on $\R^n$ and hence zero.

\smallskip\noindent(ii) \emph{Energy estimate for each profile.}
Testing each profile equation against $V^{(j)}$ gives $E_j=cm_j$.  The sharp Sobolev trace inequality bounds $\int_{\R^n_+}|\nabla V|^2$ from below by $S_\ast\big(\int_{\partial\R^n_+}|V|^q\big)^{2/q}$ for every $V\in\dot H^1(\R^n_+)$, including a sign-changing $V$. Thus $E_j\ge S_\ast m_j^{2/q}$.  Together with $E_j=cm_j$ and $m_j>0$ this gives \eqref{I-eq:profile-mass-lower-bound},
\[
m_j\ge\Big(\frac{S_\ast}{c}\Big)^{n-1}.
\]
Each profile therefore has boundary mass at least $(S_\ast/c)^{n-1}>0$, and hence $J\le(c/S_\ast)^{n-1}$. After the last extraction the concentration function of the remainder $r_k$ tends to zero. The vanishing case gives $\|r_k\|_{L^q(\partial M)}\to0$. The mass and energy identities, together with $E_j=cm_j$ and the corresponding identity for $u_\infty$, then give $\mathcal N_g^\circ(r_k)\to0$. Coercivity also yields $\|r_k\|_{H^1(M)}\to0$.
\end{proof}

\subsection{Compactness below the hemisphere threshold}

\begin{corollary}[Concentration measure and bubble energy]\label{I-lem:signfree-CC-bubble-cost}
Under the hypotheses of Proposition~\ref{I-prop:signfree-boundary-profile-decomposition}, pass to limiting profile centers $a_j=\lim_kx_k^j$. The boundary measures satisfy
\begin{equation}\label{I-eq:boundary-measure-splitting}
|u_k|^q\,d\sigma_g\overset{*}{\rightharpoonup}
|u_\infty|^q\,d\sigma_g+\sum_{a\in\mathcal A}\nu_a\delta_a,
\qquad \nu_a=\sum_{j:a_j=a}m_j,
\end{equation}
where $\mathcal A$ is the finite set of limiting centers. Every atom satisfies
\begin{equation}\label{I-eq:atom-cost-sharp}
c\nu_a\ge S_\ast\nu_a^{2/q},\qquad 0<\nu_a\le1.
\end{equation}
In particular, $c<S_\ast$ excludes all profiles and gives strong $H^1(M)$ convergence.
\end{corollary}
\begin{proof}
The strong decomposition, parameter orthogonality, and Brezis--Lieb splitting imply \eqref{I-eq:boundary-measure-splitting}: each continuous weight tends to its value at $a_j$ on the shrinking support of the $j$th profile. Since $E_j=cm_j\ge S_\ast m_j^{2/q}$ and $2/q<1$,
\[
c\nu_a=\sum_{j:a_j=a}E_j
\ge S_\ast\sum_{j:a_j=a}m_j^{2/q}
\ge S_\ast\nu_a^{2/q}.
\]
For $c<S_\ast$, the lower mass bound \eqref{I-eq:profile-mass-lower-bound} exceeds one, so there are no profiles. This recovers the measure formulation of \cite[Lemma~2.3, p.~61]{Lions85b}.
\end{proof}

\begin{theorem}[Palais--Smale compactness below the hemisphere threshold]
\label{I-thm:compactness-escobar}
Assume \textup{(Y$^+_{\partial}$)} and \textup{(BG$^{2+}$)}.  Let $q=2^*_\partial$ and let $(u_k)$ be a Palais--Smale sequence for $\mathcal J_g$ subject to
\[
\|u_k\|_{L^q(\partial M)}=1
\]
at a level $c<S_\ast=C^*_{\Esc}(S^n_+)$. Then $(u_k)$ is precompact in $H^1(M)$.

In particular, if $C^*_{\Esc}(M,g)<S_\ast$, then every Palais--Smale sequence at level $C^*_{\Esc}(M,g)$ is precompact.
\end{theorem}

\begin{proof}
The preceding proof gives $J=0$ in \eqref{I-eq:signfree-profile-decomp}; the remainder converges strongly to zero in $H^1$, proving the assertion.
\end{proof}

\begin{corollary}[Compactness of minimizing sequences below the hemisphere threshold]
\label{I-cor:escobar-precompact-cov}
Let $(u_k)$ be a minimizing sequence for $\mathcal J_g$, so that
\[
\|u_k\|_{L^q(\partial M)}=1,
\qquad
\Def^{\Esc}_g(u_k)\to0.
\]
If $C^*_{\Esc}(M,g)<S_\ast$, then $(u_k)$ is precompact in $H^1(M)$.  Along a subsequence,
\[
u_k\to u_\infty\quad\text{strongly in }H^1(M),
\qquad
\|u_\infty\|_{L^q(\partial M)}=1,
\]
and $u_\infty$ is an extremal.
\end{corollary}

\begin{proof}
The boundary constraint is closed in $H^1(M)$ and the quotient is bounded below there. Ekeland's principle \cite[Theorem~1]{Ekeland1979} gives a Palais--Smale minimizing sequence $(\widetilde u_k)$ with $\|u_k-\widetilde u_k\|_{H^1}\to0$. Theorem~\ref{I-thm:compactness-escobar} makes $(\widetilde u_k)$, and hence $(u_k)$, precompact. Along a strongly convergent subsequence, continuity of the trace and energy gives $\|u_\infty\|_q=1$ and $\mathcal J_g[u_\infty]=C^*_{\Esc}(M,g)$.
\end{proof}

\begin{lemma}[Mass bounds for half-space profiles]
\label{I-lem:signfree-mass-gap}
Assume \textup{(Y$^+_\partial$)} and \textup{(BG$^{2+}$)}, let $q=2^*_\partial$, and let $(u_k)\subset\mathcal S$ be a constrained Palais--Smale sequence for $\mathcal J_g$ at level $c>0$, with profile decomposition
\eqref{I-eq:signfree-profile-decomp}.  Write
\[
\beta(c)\coloneqq \Big(\frac{S_\ast}{c}\Big)^{n-1},
\qquad
V^{(j)}_\pm\coloneqq \max(\pm V^{(j)},0),
\qquad
m_j^\pm\coloneqq \int_{\partial\R^n_+}\bigl(V^{(j)}_\pm\bigr)^q\,dy'.
\]

\textup{(i)} \emph{Mass of a profile.} If $V^{(j)}$ does not change sign, then $|V^{(j)}|$ is a half-space
optimizer up to translation and dilation, and
\[
m_j=\beta(c).
\]
If $V^{(j)}$ changes sign, then
\[
m_j^+>\beta(c),\qquad m_j^->\beta(c),
\qquad\text{hence}\qquad m_j>2\beta(c).
\]
In particular no profile has boundary mass in the interval $\bigl(\beta(c),2\beta(c)\bigr]$. For the critical boundary equation on the half-space, Almaraz--Wang~\cite[Theorem~1.4]{AlmarazWang2022} proved that the energy of every sign-changing solution is uniformly greater than twice the bubble energy.

\textup{(ii)} \emph{Zero weak limit.} Suppose $u_\infty\equiv0$, and let
$N_1$ and $N_\pm$ be the numbers of profiles that do not change sign and that change sign, respectively. Then
\[
\bigl(N_1+2N_\pm\bigr)\beta(c)\le1,
\]
with equality if and only if $N_\pm=0$, in which case $c=N_1^{1/(n-1)}S_\ast$ and no profile changes sign, with boundary mass $1/N_1$ for each profile.

\textup{(iii)} \emph{Range of possible levels.} Write
$\Lambda_k\coloneqq k^{1/(n-1)}S_\ast$. There is no Palais--Smale sequence with zero weak limit at any level
\[
c\in\bigl(S_\ast,\Lambda_2\bigr),
\qquad
\Lambda_2=2^{1/(n-1)}S_\ast .
\]
At the level $c=\Lambda_2$, every such sequence has exactly two half-space bubble profiles, each of a fixed sign and boundary mass $\tfrac12$; if moreover
\[
\liminf_{k\to\infty}\int_{\partial M}(u_k^+)^q\,d\sigma_g>0
\qquad\text{and}\qquad
\liminf_{k\to\infty}\int_{\partial M}(u_k^-)^q\,d\sigma_g>0,
\]
the two bubbles have opposite signs. A sign-changing profile requires $c>\Lambda_2$; among the levels defined above, $\Lambda_3$ is the first at which such a profile is not excluded by the mass bound.
\end{lemma}

\begin{proof}
Fix a profile $V=V^{(j)}$ and put $m^\pm=m_j^\pm$.

(i) If $V$ has a fixed sign, the strong maximum principle and Hopf lemma give $V>0$ or $V<0$. The Li--Zhu classification \cite[Theorem~1.2]{LiZhu1995} identifies $|V|$ as a half-space optimizer. Hence
\[
cm_j=E_j=S_\ast m_j^{2/q},
\qquad m_j=(S_\ast/c)^{n-1}=\beta(c).
\]
If $V$ changes sign, testing its equation with $V_+$ and $-V_-$ gives
\[
\int_{\R^n_+}|\nabla V_\pm|^2=cm^\pm
\ge S_\ast(m^\pm)^{2/q},
\qquad m^\pm\ge\beta(c).
\]
Equality for $m^+$ would make $V_+$ an optimizer with strictly positive boundary trace. Since $V_+V_-=0$, this would force $m^-=0$, a contradiction. Thus both inequalities are strict.

(ii) When $u_\infty=0$, the exact mass splitting gives
\[
1=\sum_jm_j\ge(N_1+2N_\pm)\beta(c),
\]
with equality precisely when $N_\pm=0$. Then $N_1\beta(c)=1$, so each mass is $1/N_1$ and $c=\Lambda_{N_1}$.

(iii) For $S_\ast<c<\Lambda_2$, one has $1/2<\beta(c)<1$. A sign-changing profile would exceed the total mass, while (ii) would require the integer $N_1=1/\beta(c)$ to lie in $(1,2)$. At $c=\Lambda_2$, one has $\beta(c)=1/2$, so there are exactly two profiles, each of a fixed sign and mass $1/2$. If both signed boundary masses have positive limits, these profiles must have opposite signs. Finally, $m_j>2\beta(c)$ permits a sign-changing profile only when $c>\Lambda_2$; the first level $\Lambda_k$ not excluded by this bound is $\Lambda_3$.
\end{proof}

\begin{lemma}[Single-bubble blow-up at the threshold]
\label{I-lem:extremal-profile-onesigned}
Assume $n\ge3$, \textup{(Y$^+_\partial$)}, and \textup{(BG$^{2+}$)}. Let $(u_k)\subset\mathcal S$ be a Palais--Smale sequence at level $S_\ast$ that is not precompact in $H^1(M)$. After passing to a subsequence, there are $\sigma\in\{\pm1\}$, centers $x_k\in\partial M$, and scales $\varepsilon_k\downarrow0$ such that
\begin{equation}\label{I-eq:threshold-one-bubble-decomposition}
u_k
=
\sigma\mathcal T_{x_k,\varepsilon_k}
\bigl(\chi_{R(\varepsilon_k)}U_\ast\bigr)
+
r_k,
\qquad
\|r_k\|_{H^1(M)}\to0.
\end{equation}
Here $U_\ast$ is the boundary-normalized optimizer \eqref{B-recall-I-eq:mass-one-Uplus}. The boundary concentration is a single atom of full mass, and the weak limit is zero.
\end{lemma}

\begin{proof}
By \eqref{I-eq:exact-PS-multiplier}, take $\lambda_k=\mathcal J_g[u_k]\to S_\ast$ to obtain
\begin{equation}\label{I-eq:threshold-PS-multiplier}
\mathcal N_g^\circ(u_k,\phi)-\lambda_k
\int_{\partial M}|u_k|^{q-2}u_k\phi\,d\sigma_g
=o(1)\|\phi\|_{H^1}.
\end{equation}
Pass first to a subsequence having no strongly convergent subsequence, and apply Proposition~\ref{I-prop:signfree-boundary-profile-decomposition}. Since the remainder converges strongly to zero, there must be at least one profile. At $c=S_\ast$, every profile has mass at least one, while the total boundary mass is one. Thus there is exactly one profile, with mass one, and the weak limit has zero trace. Its equation and \textup{(Y$^+_\partial$)} imply $u_\infty=0$. Lemma~\ref{I-lem:signfree-mass-gap}(i) identifies the profile as a signed translate and dilation of $U_\ast$.

Absorb its fixed translation and dilation into the center and scale. We may also replace the extraction cutoff by $R(\varepsilon)=\varepsilon^{-3/4}$. The squared Dirichlet norm of the truncation error is bounded by $CR^{2-n}$. For an arbitrary admissible extraction radius, its squared $L^2(M)$ norm is bounded by
\[
C\varepsilon^2
\begin{cases}
R,&n=3,\\
1+\log R,&n=4,\\
1,&n\ge5,
\end{cases}
\]
which tends to zero when $R\to\infty$ and $\varepsilon R\to0$. The fixed parameter change also has vanishing rescaled coordinate error. These changes are therefore $o(1)$ in $H^1(M)$ and give \eqref{I-eq:threshold-one-bubble-decomposition}.
\end{proof}

\begin{theorem}[The first mixed Palais--Smale level]
\label{I-thm:signfree-first-gap-compactness}
Assume $n\ge3$, \textup{(Y$^+_\partial$)}, and \textup{(BG$^{2+}$)}, and put
\begin{equation}\label{I-eq:first-mixed-level}
Y=C^*_{\Esc}(M,g),\qquad
c_{\mathrm{mix}}=(S_\ast^{n-1}+Y^{n-1})^{1/(n-1)},\qquad
\Lambda_2=2^{1/(n-1)}S_\ast.
\end{equation}
Every constrained Palais--Smale sequence at a level $S_\ast<c<c_{\mathrm{mix}}$ is precompact in $H^1(M)$, without a sign assumption. The endpoint $c_{\mathrm{mix}}$ admits a noncompact Palais--Smale sequence.

More precisely, if $S_\ast<c<\Lambda_2$, a noncompact constrained Palais--Smale sequence at $c$ exists if and only if
\begin{equation}\label{I-eq:mixed-critical-value}
d=(c^{n-1}-S_\ast^{n-1})^{1/(n-1)}
\end{equation}
is a critical value of $\mathcal J_g|_{\mathcal S}$.
\end{theorem}

\begin{proof}
We have $0<Y\le S_\ast$, so $S_\ast<c_{\mathrm{mix}}\le\Lambda_2$. Apply Proposition~\ref{I-prop:signfree-boundary-profile-decomposition} to a noncompact sequence at a level $S_\ast<c<\Lambda_2$, passing to a subsequence having no strongly convergent subsequence. Put $\beta=(S_\ast/c)^{n-1}\in(1/2,1)$. Lemma~\ref{I-lem:signfree-mass-gap} excludes sign-changing profiles, and each fixed-sign profile has mass $\beta$. Thus exactly one profile occurs. The weak limit $u_\infty$ has boundary mass $1-\beta>0$ and satisfies the equation with multiplier $c$.

Write $u_\infty=av$, where $a=(1-\beta)^{1/q}$ and $v\in\mathcal S$. Then $v$ is a constrained critical point with value
\[
\mathcal J_g[v]=ca^{q-2}=c(1-\beta)^{1/(n-1)}=d.
\]
Since $d\ge Y$, this implies $c\ge c_{\mathrm{mix}}$ and proves compactness in the stated interval. It also proves the necessity of \eqref{I-eq:mixed-critical-value}.

Conversely, let $v\in\mathcal S$ be critical with value $0<d<S_\ast$. Set
\[
c=(d^{n-1}+S_\ast^{n-1})^{1/(n-1)},\qquad
a_0=(d/c)^{(n-2)/2},\qquad
a_1=(S_\ast/c)^{(n-2)/2}.
\]
Fix $x\in\partial M$ and put
\[
z_\varepsilon=a_0v+a_1\mathcal U_{x,\varepsilon},\qquad
u_\varepsilon=\frac{z_\varepsilon}{\|z_\varepsilon\|_q}.
\]
The bubble converges weakly to zero in $H^1$, tends to zero almost everywhere, and has energy tending to $S_\ast$. The Brezis--Lieb identity and energy polarization give
\[
\|z_\varepsilon\|_q^q=a_0^q+a_1^q+o(1)=1+o(1),\qquad
\mathcal N_g^\circ(z_\varepsilon)=a_0^2d+a_1^2S_\ast+o(1)=c+o(1).
\]
For the derivative, critical rescaling bounds the principal metric error by $C\varepsilon R$ and the Dirichlet cutoff error by $CR^{-(n-2)/2}$. Both tend to zero for our cutoff, in every dimension $n\ge3$. The flat bubble equation therefore gives
\[
\left\|\mathcal N_g^\circ(\mathcal U_{x,\varepsilon},\cdot)
-S_\ast\int_{\partial M}\mathcal U_{x,\varepsilon}^{q-1}(\cdot)\,d\sigma_g\right\|_{H^{-1}_E}\longrightarrow0.
\]
The scalar-curvature and mean-curvature terms tend to zero by compactness of the $L^2$ embeddings. Nonlinear Brezis--Lieb splitting, as in the proof of Proposition~\ref{I-prop:signfree-boundary-profile-decomposition}, gives
\[
\bigl\||z_\varepsilon|^{q-2}z_\varepsilon
-a_0^{q-1}|v|^{q-2}v-a_1^{q-1}\mathcal U_{x,\varepsilon}^{q-1}\bigr\|_{q/(q-1)}\longrightarrow0.
\]
The identities $a_0d=ca_0^{q-1}$ and $a_1S_\ast=ca_1^{q-1}$ now show that $u_\varepsilon$ is Palais--Smale at $c$. Its nonzero concentrating mass $a_1^q$ prevents strong convergence. This proves sufficiency.

If $Y<S_\ast$, Corollary~\ref{I-cor:escobar-precompact-cov} supplies a minimizer of value $d=Y$, and the preceding construction proves noncompactness at $c_{\mathrm{mix}}$. If $Y=S_\ast$, take two bubbles at distinct fixed boundary points, each with amplitude $2^{-1/q}$, and normalize their sum. The same calculation gives a noncompact Palais--Smale sequence at $\Lambda_2=c_{\mathrm{mix}}$.
\end{proof}

\section{Lyapunov--Schmidt reduction and the strict inequality}
\label{I-sec:reduced}

The complementary correction lowers the quadratic coefficient of the standard bubble. We identify it through the flat linearized equation and evaluate the resulting quadratic form exactly. The resulting coefficient is negative at every non-umbilic point in all dimensions $n\ge5$.

\subsection{The leading term of the constrained gradient}

Write $\mathrm{Sym}^2_0(\R^{n-1})$ for the trace-free symmetric matrices. For $S$ in this space, put $|S|^2=\operatorname{tr}(S^2)$ and
\[
r=|y'|,\qquad s=r^2+(1+t)^2,\qquad Q_S(y')=S_{ab}y_ay_b.
\]
Tangential spherical harmonics decompose functions into angular degrees. Radial functions have degree $0$, translation fields have degree $1$, and $Q_S(y')h(r,t)$ has degree $2$. The flat dilation and translation Jacobi fields are
\[
Z^\flat_0=-\frac{n-2}{2}U_+-y\cdot\nabla U_+,
\qquad Z^\flat_a=-\partial_aU_+.
\]
Define the flat constrained complement and its Jacobi form by
\begin{align*}
X^\flat&=\left\{\phi\in\mathcal D^{1,2}:\int_{\partial\R^n_+}U_+^{q-1}\phi=0,
\quad\int_{\R^n_+}\nabla\phi\cdot\nabla Z^\flat_j=0\ (0\le j\le n-1)\right\},\\
\mathcal B_n(\phi,\psi)&=\int_{\R^n_+}\nabla\phi\cdot\nabla\psi
-(q-1)\lambda_{\mathrm{flat}}\int_{\partial\R^n_+}U_+^{q-2}\phi\psi.
\end{align*}
Let $X_2$ be the degree-$2$ subspace, and $X_{2,S}$ its closed subspace generated by $Q_Sh(r,t)$ for a fixed $S$. Angular orthogonality places both in $X^\flat$. We use the cutoff bound $\omega(R)=R^{-(n-2)/2}$ and the metric remainder bound
\[
\eta_n(\varepsilon)=
\begin{cases}
\varepsilon^{13/8},&n=5,\\
\varepsilon^2\sqrt{1+\log(1/\varepsilon)},&n=6,\\
\varepsilon^2,&n\ge7.
\end{cases}
\]
With the choice $R(\varepsilon)=\varepsilon^{-3/4}$, $\eta_n(\varepsilon)+\omega(R(\varepsilon))=o(\varepsilon)$ for $n\ge5$ and is $O(\varepsilon^2)$ for $n\ge8$.

\begin{lemma}[Constrained gradient of a bubble: degree-$2$ term orthogonal to the kernel]\label{I-lem:residual-kernel-pure}\label{II-lem:residual-kernel-pure}
Assume $n\ge5$, \textup{(Y$^+_{\partial}$)}, and \textup{(BG$^{3+}$)}. Use the normalized bubble of Definition~\ref{B-recall-I-def:bubbles} with the cutoff $R(\varepsilon)=\varepsilon^{-3/4}$, and write $t(x,\varepsilon)=\|U_{+,x,\varepsilon}\|_q^{-1}$. Its normalized parameter derivatives are
\[
Z_0=\varepsilon\partial_\varepsilon\mathcal U_{x,\varepsilon},
\qquad Z_\alpha=\varepsilon\partial_{x^\alpha}\mathcal U_{x,\varepsilon}
\quad(1\le\alpha\le n-1).
\]
Set $\mathcal K_{x,\varepsilon}=\operatorname{span}\{Z_0,\ldots,Z_{n-1}\}\subset T_{\mathcal U_{x,\varepsilon}}\mathcal S$, and let $\Pi_{x,\varepsilon}$ be the energy-orthogonal projection onto $\mathcal K_{x,\varepsilon}^\perp\cap T_{\mathcal U_{x,\varepsilon}}\mathcal S$.

Then, uniformly for $x\in\partial M$ as $\varepsilon\downarrow0$, there exist coefficients $a_\beta(x,\varepsilon)\in\mathbb R$ ($\beta=0,\dots,n-1$) such that
\begin{equation}\label{I-eq:residual-kernel-decomp}
\mathsf{grad}_{\mathcal S}\mathcal J_g(\mathcal U_{x,\varepsilon})
=\sum_{\beta=0}^{n-1} a_\beta(x,\varepsilon)Z_\beta(x,\varepsilon)+\mathcal R_{x,\varepsilon},
\qquad
\mathcal R_{x,\varepsilon}=\Pi_{x,\varepsilon}\mathsf{grad}_{\mathcal S}\mathcal J_g(\mathcal U_{x,\varepsilon})
\in\mathcal K_{x,\varepsilon}^{\perp}\cap T_{\mathcal U_{x,\varepsilon}}\mathcal S,
\end{equation}
and
\begin{equation}\label{I-eq:residual-perp-bound}
\|\mathcal R_{x,\varepsilon}\|_{H^1(M)}\le C\bigl(\varepsilon|H_g(x)|+\varepsilon|\mathring{\mathrm{II}}(x)|+\varepsilon^2+\omega(R(\varepsilon))\bigr),
\end{equation}
with $\omega$ as above, and $|a_\beta(x,\varepsilon)|\le C\varepsilon$. In particular, the constrained gradient is $O(\varepsilon)$; when $H_g(x)=0$, the first term in the bound disappears. At an umbilic point with zero mean curvature the bound is $O(\varepsilon^2+\omega(R(\varepsilon)))$, which is $O(\varepsilon^2)$ for $n\ge8$ and $o(\varepsilon)$ for $n=5,6,7$.
\end{lemma}

\begin{proof}
The normalized scale and center fields have a uniformly invertible Gram matrix by flat nondegeneracy and Fermi rescaling. Their energy projection gives \eqref{I-eq:residual-kernel-decomp}. Indeed, with $G_{\alpha\beta}=\langle Z_\alpha,Z_\beta\rangle_E$,
\[
a_\beta=\sum_\alpha(G^{-1})_{\beta\alpha}
\langle\mathsf{grad}_{\mathcal S}\mathcal J_g(\mathcal U),Z_\alpha\rangle_E,
\qquad |a_\beta|\le C\varepsilon
\]
by the bounded normalized frame, its inverse Gram matrix, and Lemma~\ref{I-lem:first-order-drift}. It remains to estimate the differential on the complement.

The flat equation cancels the leading term. The first-order Fermi and Robin corrections give
\begin{equation}\label{I-eq:first-residual-spin-split}
D\mathcal J_g[\mathcal U_{x,\varepsilon}](\psi)
=\varepsilon H_g(x)\ell_0(\psi)
+\varepsilon\mathring{\mathrm{II}}_{ab}(x)\ell_{ab}^{(2)}(\psi)
+O\bigl((\eta_n(\varepsilon)+\omega(R(\varepsilon)))\|\psi\|_E\bigr).
\end{equation}
Here $\ell_0$ is the radial functional from the trace of $\mathrm{II}$ and the boundary mean curvature; the trace-free contraction is the degree-$2$ functional generated by $t\mathring{\mathrm{II}}_{ab}\partial_{ab}U_+$. These functionals are bounded because $|y|\nabla U_+\in L^2$ for $n\ge5$. The metric remainder follows from
\[
\varepsilon^2\bigl\||y|^2\nabla U_+\bigr\|_{L^2(B_{2R}^+)}
\le C\varepsilon^2
\begin{cases}R^{1/2},&n=5,\\(1+\log R)^{1/2},&n=6,\\1,&n\ge7,\end{cases}
\]
with the scalar and trace terms estimated as in Lemma~\ref{I-lem:first-order-drift}. For $R=\varepsilon^{-3/4}$ this is $C\eta_n(\varepsilon)\le C(\varepsilon^2+\omega(R))$.

The radial source vanishes when $H_g(x)=0$. The trace-free source is orthogonal to the dilation and translation fields, since
\[
\int_{\Sph^{n-2}}\mathring{\mathrm{II}}_{ab}\omega_a\omega_b\,d\omega=0,
\qquad
\int_{\Sph^{n-2}}\mathring{\mathrm{II}}_{ab}\omega_a\omega_b\omega_c\,d\omega=0.
\]
Consequently, on the complement,
\begin{equation}\label{I-eq:residual-dual-test}
|D\mathcal J_g[\mathcal U_{x,\varepsilon}](\psi)|
\le C\bigl(\varepsilon|H_g(x)|+\varepsilon|\mathring{\mathrm{II}}(x)|+\varepsilon^2+\omega(R(\varepsilon))\bigr)\|\psi\|_E.
\end{equation}
Riesz representation and energy--$H^1$ equivalence give \eqref{I-eq:residual-perp-bound}. Setting $H_g(x)=0$ and then $\mathring{\mathrm{II}}(x)=0$ gives the stated refinements.
\end{proof}

\begin{lemma}[Leading degree-$2$ term of the projected gradient]
\label{I-lem:spin2-profile}
Assume $n\ge5$ and retain $X^\flat$ and $X_2$ from the beginning of this subsection. For $S\in\mathrm{Sym}^2_0(\mathbb R^{n-1})$, define
\begin{equation}\label{I-eq:spin2-source-functional}
\mathscr L_S(\phi)
\coloneqq
-2\int_{\mathbb R^n_+}
y_n S^{ab}\partial_aU_+\partial_b\phi\,dy
=
2\int_{\mathbb R^n_+}
y_n S^{ab}\partial_{ab}^2U_+\phi\,dy,
\end{equation}
where the second identity is tangential integration by parts ($y_n=0$ on $\partial\mathbb R^n_+$). Let $\mathcal E_S^\flat\in X^\flat$ be the $\dot H^1$-Riesz representative of $\mathscr L_S$.

In the setting of Lemma~\ref{I-lem:residual-kernel-pure}, assume $H_g(x)=0$. Set $X_{x,\varepsilon}\coloneqq\mathcal K_{x,\varepsilon}^{\perp_E}
\cap T_{\mathcal U_{x,\varepsilon}}\mathcal S$, and let
$\mathscr L_S^{x,\varepsilon}\in X_{x,\varepsilon}^\ast$ denote the covector obtained from $\mathscr L_S$ by the rescaled Fermi pullback, localized with the same cutoff convention as the boundary bubble. Define the restricted differential $\mathscr R_{x,\varepsilon}\coloneqq D(\mathcal J_g\circ\mathsf{ret}_{\mathcal U_{x,\varepsilon}})(0)|_{X_{x,\varepsilon}}$. Then, in the norm dual to the energy norm on $X_{x,\varepsilon}$,
\begin{equation}\label{I-eq:spin2-profile-strong}
\left\|
\mathscr R_{x,\varepsilon}
+2\varepsilon\,t(x,\varepsilon)\,
\mathscr L_{\mathring{\mathrm{II}}(x)}^{x,\varepsilon}
\right\|_{X_{x,\varepsilon}^\ast}
\le
C\big(\eta_n(\varepsilon)+\omega(R(\varepsilon))\big),
\end{equation}
where $\eta_n$ and $\omega$ are defined above. The error is $o(\varepsilon)$ for $n\ge5$ and $O(\varepsilon^2)$ for $n\ge8$. The first variation of the flat cutoff bubble is radial, so its degree-$2$ projection is zero. The leading source is therefore $-2\varepsilon t_\ast\mathscr L_{\mathring{\mathrm{II}}(x)}^{x,\varepsilon}$ for $\mathcal J_g$, or $-2\varepsilon t_\ast^{-1}\mathscr L_{\mathring{\mathrm{II}}(x)}^{x,\varepsilon}$ for $\mathcal J_g/S_\ast$. The map $S\mapsto \mathcal E_S^\flat$ is $O(n-1)$-equivariant and injective on $\mathrm{Sym}^2_0(\mathbb R^{n-1})$.
\end{lemma}

\begin{proof}
At the unnormalized bubble, differentiating the first geometric energy correction gives
\[
4\varepsilon\int t\mathrm{II}_{ab}(x)\partial_aU_+\partial_b\phi
=-2\varepsilon\mathscr L_{\mathrm{II}(x)}(\phi).
\]
Under the unit boundary $L^q$-norm constraint the first variation is linear in the base function. At $\mathcal U_{x,\varepsilon}=t(x,\varepsilon)U_{+,x,\varepsilon}$ the source is therefore multiplied by $t(x,\varepsilon)$. The degree-$2$ term is unchanged by flat projection. The separate error bounds are
\[
\begin{array}{ll}
\text{second metric jet:}&C\varepsilon^2\ell_{2,n}(R),\\
\text{flat cutoff:}&CR^{-\alpha},\\
\text{tail of the first-order source:}&C\varepsilon R^{-(n-4)/2},\\
\text{projection:}&C(\varepsilon^2+\varepsilon R^{-\alpha}),
\end{array}
\]
where $\ell_{2,5}(R)=R^{1/2}$, $\ell_{2,6}(R)=\sqrt{\log(2R)}$, and $\ell_{2,n}(R)=1$ for $n\ge7$. The first two estimates follow from the weighted moment in Lemma~\ref{I-lem:residual-kernel-pure} and \eqref{I-eq:cutoff-first-variation-identity}. The source tail is bounded by $CR^{-\alpha}$ because $\varepsilon R\le1$. For the projection, the normalized frame differs from the critically rescaled flat frame by $O_E(\varepsilon+R^{-\alpha})$. This follows from the Fermi expansion, the explicit profile derivatives, and $\partial_xd_{\bar g}(x,z')^2=-2\exp_x^{-1}(z')$. The constraint normal is handled by the exact Riesz identity in Remark~\ref{B-rem:inverse-Gram}. Angular orthogonality removes the flat moments, and the inverse Gram formula \eqref{II-eq:projection-inverse-Gram-formula-onebubble} gives the stated projection error. All four terms are bounded by the right-hand side of \eqref{I-eq:spin2-profile-strong}.

The boundary normalization satisfies $t(x,\varepsilon)=t_\ast+O(\varepsilon^2)$ by \cite[Lemma~3.8]{EscobarGeometry}, which gives the stated leading coefficient.

Rotational symmetry gives equivariance. The explicit identity
\[
S_{ab}\partial_{ab}U_+=n(n-2)A_nQ_Ss^{-(n+2)/2}
\]
shows that $\mathscr L_S$ is nonzero on degree-$2$ tests whenever $S\ne0$. Hence its Riesz representative is nonzero, proving injectivity.
\end{proof}

\begin{remark}[The geometric source]\label{I-rem:geometric-source}
For $F_{x,\varepsilon}(v)=\mathcal J_g(\mathsf{ret}_{\mathcal U_{x,\varepsilon}}v)$ on $X_{x,\varepsilon}$, the complementary source is the covector $D_XF_{x,\varepsilon}(0)\in X_{x,\varepsilon}^*$. In fixed model coordinates its leading coefficient in a minimal boundary region is
\[
\mathscr S_{1,x}=-2t_\ast\mathscr L_{\mathring{\mathrm{II}}(x)},
\qquad D_XF_{x,\varepsilon}(0)=\varepsilon\mathscr S_{1,x}+o(\varepsilon).
\]
For $n\ge7$, at an umbilic center with the refined cutoff, the corresponding expansion starts with $\varepsilon^2\mathscr S_{2,x}$; see \eqref{II-eq:second-source-explicit}. These leading covectors are the geometric sources. The full complementary source also includes cutoff and higher-order terms. Its energy Riesz representative is a vector, defined by pairing with the energy inner product. For $\mathcal J_g/S_\ast$, all source covectors are divided by $S_\ast$.
\end{remark}

\begin{proposition}[Lyapunov--Schmidt correction for the degree-$2$ source]
\label{I-prop:spin2-schur-identification}
Assume $n\ge5$. Let $\mathcal B_n$ be the flat Jacobi form on the constrained complement $X^\flat$ and let $\mathscr L_S$ be the degree-$2$ source functional from \eqref{I-eq:spin2-source-functional}. Define $\delta\kappa_3^{\mathrm{LS}}(n)$ by the first equality in \eqref{I-eq:spin2-schur-general} with $|S|=1$. Its value is independent of the choice of unit $S$, is strictly positive, and for every $S\in\mathrm{Sym}^2_0(\mathbb R^{n-1})$ one has
\begin{equation}\label{I-eq:spin2-schur-general}
\delta\kappa_3^{\mathrm{LS}}(n)|S|^2
=
\sup_{\phi\in X^\flat}
\left\{
2\mathscr L_S(\phi)-\mathcal B_n(\phi,\phi)
\right\}
=
\sup_{\phi\in X^\flat\setminus\{0\}}
\frac{\mathscr L_S(\phi)^2}{\mathcal B_n(\phi,\phi)}.
\end{equation}
This is the Lyapunov--Schmidt correction entering the reduced decomposition \eqref{I-eq:Rg-red-decomp}. Its exact value is computed in Proposition~\ref{I-prop:quadratic-exact}.
\end{proposition}

\begin{proof}
The trace inequality makes $\mathcal B_n$ continuous on $X^\flat$, and the flat spectral gap from Lemma~\ref{I-lem:hessian-sphere} gives $\mathcal B_n(\phi,\phi)\ge\gamma_\ast\|\phi\|_{\dot H^1}^2$. Moreover,
\[
|\mathscr L_S(\phi)|\le 2|S|\|y_n\nabla_{\tan}U_+\|_{L^2}\|\phi\|_{\dot H^1}\le C_n|S|\|\phi\|_{\dot H^1},
\]
since $y_n\nabla_{\tan}U_+\in L^2(\mathbb R^n_+)$ for $n\ge5$. Thus Lax--Milgram gives a unique $w_S\in X^\flat$ satisfying $\mathcal B_n(w_S,\psi)=\mathscr L_S(\psi)$ for every $\psi\in X^\flat$. Completing the square yields
\[
2\mathscr L_S(\phi)-\mathcal B_n(\phi,\phi)=\mathcal B_n(w_S,w_S)-\mathcal B_n(\phi-w_S,\phi-w_S).
\]
Hence the first supremum is attained uniquely at $w_S$ and equals $\mathcal B_n(w_S,w_S)$. For $\phi\ne0$, maximizing $2c\mathscr L_S(\phi)-c^2\mathcal B_n(\phi,\phi)$ over $c\in\mathbb R$ gives $c=\mathscr L_S(\phi)/\mathcal B_n(\phi,\phi)$ and maximum $\mathscr L_S(\phi)^2/\mathcal B_n(\phi,\phi)$. Taking the supremum over these lines proves the second equality in \eqref{I-eq:spin2-schur-general}.

The map $S\mapsto w_S$ is linear and $O(n-1)$-equivariant, since $X^\flat$ and $\mathcal B_n$ are invariant under tangential rotations. Therefore $S\mapsto\mathcal B_n(w_S,w_S)$ is an orthogonally invariant quadratic form on trace-free symmetric matrices. Diagonalizing $S$ and permuting its eigenvalues shows that it is a linear combination of $\operatorname{tr}(S^2)$ and $(\operatorname{tr}S)^2$. The latter is zero, so the form equals $\delta\kappa_3^{\mathrm{LS}}(n)|S|^2$. Lemma~\ref{I-lem:spin2-profile} gives $\mathscr L_S\ne0$ for $S\ne0$, proving $\delta\kappa_3^{\mathrm{LS}}(n)>0$; when $S=0$, $w_S=0$ and both suprema are zero. The coefficients of $\mathcal B_n$ depend on $y'$ only through $|y'|$. Thus distinct tangential spherical harmonics are $\mathcal B_n$-orthogonal. Since $\mathscr L_S$ pairs only with the harmonic $S_{ab}y_ay_b$, uniqueness places $w_S$ in the closed subspace generated by $S_{ab}y_ay_bh(|y'|,y_n)$. For $S\ne0$, both suprema can therefore be restricted to this subspace.

To identify the normalization in the geometric expansion, put $S=\mathring{\mathrm{II}}(x)$ in a boundary region where $H_g\equiv0$ and rescale $\phi\in X^\flat$ in Fermi coordinates, multiply it by the base bubble factor $t(x,\varepsilon)$, and project it onto $X_{x,\varepsilon}$; denote this direction by $\phi_{x,\varepsilon}$. Lemma~\ref{I-lem:spin2-profile} gives the leading covector $-2\varepsilon t_\ast^{-1}\mathscr L_S$ for $\mathcal J_g/S_\ast$ at $t_\ast U_+$, while the limiting quadratic form is $t_\ast^{-2}\mathcal B_n$, where $t_\ast^2=S_\ast$. The factor $t_\ast$ in this rescaled direction cancels both normalization factors. Consequently, for fixed $c\in\mathbb R$,
\[
\frac{\mathcal J_g(\mathcal U_{x,\varepsilon}+c\varepsilon\phi_{x,\varepsilon})}{S_\ast}
=1+\varepsilon^2\left(\kappa_3^{\mathrm{bare}}(n)|S|^2-2c\mathscr L_S(\phi)+c^2\mathcal B_n(\phi,\phi)\right)+o(\varepsilon^2).
\]
The preceding minimization lowers the quadratic coefficient by $\delta\kappa_3^{\mathrm{LS}}(n)|S|^2$. The passage from this flat minimum to the parameter-dependent Lyapunov--Schmidt correction is established in Lemma~\ref{I-lem:resolvent-limit} and the proof of \eqref{I-eq:Rg-red-decomp} below.
\end{proof}

\subsection{Construction of the correction}

The construction uses the Lyapunov--Schmidt scheme \cite[Chapter~2]{AmbrosettiMalchiodi2006}; a boundary Yamabe application with a geometric bubble correction is given in \cite[Section~3]{GhimentiMichelettiPistoia2018}.

\begin{remark}[Orthogonal projection onto the complement]\label{B-rem:inverse-Gram}
Let $P_{x,\varepsilon}=\Pi_{x,\varepsilon}$ be the energy-orthogonal projection onto $X_{x,\varepsilon}$, and choose a basis $\{\zeta_\beta\}$ for its orthogonal complement, consisting of the constraint normal and normalized modulation directions as in Remark~\ref{I-rem:banach-trivialization}. Put $G_{\beta\gamma}=\langle\zeta_\beta,\zeta_\gamma\rangle_{E,g}$. Then
\begin{equation}\label{II-eq:projection-inverse-Gram-formula-onebubble}
(I-P_{x,\varepsilon})R
=\sum_{\beta,\gamma}\zeta_\beta(G^{-1})_{\beta\gamma}
\langle R,\zeta_\gamma\rangle_{E,g}.
\end{equation}
Differentiating this projection formula uses $D(G^{-1})=-G^{-1}(DG)G^{-1}$. The constraint normal also has the required flat limit. For a boundary-normalized bubble $\mathcal U$, the quotient identity gives
\[
\mathcal N_g^\circ(\mathcal U,\phi)
=\mathcal J_g[\mathcal U]
\int_{\partial M}\mathcal U^{q-1}\phi\,d\sigma_g
+\tfrac12D\mathcal J_g[\mathcal U][\phi].
\]
Thus the energy Riesz representative of the covector $\phi\mapsto q\int\mathcal U^{q-1}\phi$ is
\[
\frac{q}{\mathcal J_g[\mathcal U]}\mathcal U
-\frac{q}{2\mathcal J_g[\mathcal U]}
R_E D\mathcal J_g[\mathcal U],
\]
where $R_E$ is the energy Riesz map. The first-variation bound makes its second term $O(\varepsilon)$, or $O(\varepsilon^2)$ at an umbilic point with the improved cutoff. Thus the normal converges in the global energy norm; the corresponding estimate for the parameter derivatives of the first variation gives the differentiated convergence.
\end{remark}

\begin{remark}[Local trivialization of the complements]\label{I-rem:banach-trivialization}
Write $\mathbf p=(x,\varepsilon)$, fix $X_0=X_{\mathbf p_0}$, and let $\Pi_0$ be its energy projection. For nearby parameters, $\Pi_0|_{X_{\mathbf p}}$ is an isomorphism and
\[
\mathcal I_{\mathbf p}=(\Pi_0|_{X_{\mathbf p}})^{-1}:X_0\longrightarrow X_{\mathbf p}
\]
identifies the complements. Their orthogonal complements have a frame consisting of the energy Riesz representative of $q\int\mathcal U_{\mathbf p}^{q-1}\phi$ and the normalized parameter fields. The explicit cutoff profiles have the required $H^1$ parameter derivatives at positive scales; the powered cutoff in Convention~\ref{II-conv:cutoff-choice} gives the same regularity to $\mathcal U^{q-1}$. The trace theorem and fixed Riesz map transfer it to the constraint normal. Formula~\eqref{II-eq:projection-inverse-Gram-formula-onebubble} therefore makes $\Pi_{\mathbf p}$ and $\mathcal I_{\mathbf p}$ $C^1$ in all parameters and locally $C^2$ in scale. The normalized frame, its logarithmic scale derivative, and its inverse Gram matrix are uniformly bounded on the local logarithmic-scale neighborhoods. Hence $\varepsilon\partial_\varepsilon\Pi_{\mathbf p}$ and $\varepsilon\partial_\varepsilon\mathcal I_{\mathbf p}$ are uniformly bounded there. Lemma~\ref{I-lem:finite-regularity-preparation} covers the case of $C^3$ metrics.
\end{remark}

\begin{lemma}[Uniform $C^1$ Lyapunov--Schmidt inversion]\label{I-lem:LS-C1-inversion}
In the local coordinates of Remark~\ref{I-rem:banach-trivialization}, put
\[
F_{\mathbf p}(w)=D_w\bigl[\mathcal J_g(\mathcal U_{\mathbf p}+\mathcal I_{\mathbf p}w)\bigr],
\qquad R_{\mathbf p}=F_{\mathbf p}(0),\qquad L_{\mathbf p}=D_wF_{\mathbf p}(0).
\]
Suppose $F_{\mathbf p}$ is $C^1$ jointly in $(\mathbf p,w)$, its $w$ derivative is uniformly continuous on a fixed ball, and $\|L_{\mathbf p}^{-1}\|\le C_L$. For sufficiently small $\|R_{\mathbf p}\|$, the equation $F_{\mathbf p}(w)=0$ has a unique solution in a fixed small ball, with
\begin{equation}\label{I-eq:LS-C1-size}
\|w_{\mathbf p}\|\le2C_L\|R_{\mathbf p}\|.
\end{equation}
The solution is $C^1$ in $\mathbf p$ and
\[
D_{\mathbf p}w_{\mathbf p}
=-\bigl(D_wF_{\mathbf p}(w_{\mathbf p})\bigr)^{-1}
D_{\mathbf p}F_{\mathbf p}(w_{\mathbf p}).
\]
For the local coordinate representations of the quotient, the required uniform continuity follows from the $C^2$ boundary trace functional, since $q>2$.
\end{lemma}
\begin{proof}
Choose $r_0$ so that $\|L_{\mathbf p}^{-1}(D_wF_{\mathbf p}(w)-L_{\mathbf p})\|\le1/2$ for $\|w\|\le r_0$. If $2C_L\|R_{\mathbf p}\|\le r_0$, the map $w\mapsto w-L_{\mathbf p}^{-1}F_{\mathbf p}(w)$ is a contraction of that ball into itself. Its fixed point satisfies \eqref{I-eq:LS-C1-size}; the parameter-dependent implicit-function theorem gives the derivative formula.
\end{proof}

\begin{lemma}[Expansion of the Lyapunov--Schmidt functional]
\label{II-lem:holder-quadratic-expansion}
Fix an admissible parameter $\mathbf p$ for a family of normalized approximate functions $\mathcal U_{\mathbf p}$. Let $\mathcal T_{\mathbf p}:X_0\to X_{\mathbf p}$ identify a fixed Hilbert space $X_0$ with each complementary space $X_{\mathbf p}$. Assume that the local coordinate representations of the quotient below have uniform $C^{2,\gamma}$ bounds on a fixed small ball, where $0<\gamma\le1$ (for $n\ge5$, $\gamma=q-2$). In the one-bubble case, use $\mathcal T_{\mathbf p}=\mathcal I_{\mathbf p}$ from Remark~\ref{I-rem:banach-trivialization}. Write
\[
\mathscr H_{\mathbf p}(\tilde w)=\mathcal J_g\bigl(\mathsf{ret}_{\mathcal U_{\mathbf p}}(\mathcal T_{\mathbf p}\tilde w)\bigr),\quad
\mathcal S_{\mathbf p}=D_{\tilde w}\mathscr H_{\mathbf p}(0),
\quad \mathcal L_{\mathbf p}=D_{\tilde w}^2\mathscr H_{\mathbf p}(0).
\]
Assume $\mathcal L_{\mathbf p}^{-1}$ is uniformly bounded. If
\[
D_{\tilde w}\mathscr H_{\mathbf p}(\tilde w_{\mathbf p})=0,
\qquad
\|\tilde w_{\mathbf p}\|
\le C\|\mathcal S_{\mathbf p}\|,
\]
then
\begin{equation}\label{II-eq:holder-w-quadratic}
\tilde w_{\mathbf p}
=-\mathcal L_{\mathbf p}^{-1}\mathcal S_{\mathbf p}
+O\bigl(\|\mathcal S_{\mathbf p}\|^{1+\gamma}\bigr),
\end{equation}
and
\begin{equation}\label{II-eq:holder-energy-quadratic}
\mathscr H_{\mathbf p}(\tilde w_{\mathbf p})
=\mathscr H_{\mathbf p}(0)
-\frac12
\left\langle
\mathcal L_{\mathbf p}^{-1}\mathcal S_{\mathbf p},
\mathcal S_{\mathbf p}
\right\rangle
+O\bigl(\|\mathcal S_{\mathbf p}\|^{2+\gamma}\bigr).
\end{equation}
\end{lemma}

\begin{proof}
The inequality $\bigl||a|^\gamma-|b|^\gamma\bigr|\le|a-b|^\gamma$ and the trace embedding give, on a fixed small ball,
\[
\|D^2\mathscr H_{\mathbf p}(\tilde w_1)
-D^2\mathscr H_{\mathbf p}(\tilde w_0)\|
\le C\|\tilde w_1-\tilde w_0\|^\gamma.
\]
Hence
\begin{align*}
D\mathscr H_{\mathbf p}(\tilde w)
&=\mathcal S_{\mathbf p}+\mathcal L_{\mathbf p}\tilde w+O(\|\tilde w\|^{1+\gamma}),\\
\mathscr H_{\mathbf p}(\tilde w)
&=\mathscr H_{\mathbf p}(0)+\langle\mathcal S_{\mathbf p},\tilde w\rangle
+\frac12\langle\mathcal L_{\mathbf p}\tilde w,\tilde w\rangle
+O(\|\tilde w\|^{2+\gamma}).
\end{align*}
Applying the first identity at $\tilde w_{\mathbf p}$ gives
\eqref{II-eq:holder-w-quadratic}. Substitution into the second identity gives
\eqref{II-eq:holder-energy-quadratic}.
\end{proof}

\begin{lemma}[Lyapunov--Schmidt reduction and estimates for the correction]\label{I-lem:LS}
Assume $n\ge5$, \textup{(Y$^+_\partial$)}, and \textup{(BG$^{3+}$)}. For the normalized bubble $\mathcal U=\mathcal U_{x,\varepsilon}$ and complement $X=X_{x,\varepsilon}$, there is a unique small solution $w_{x,\varepsilon}\in X$, with $\|w_{x,\varepsilon}\|\le c_1\varepsilon$, of
\begin{equation}\label{I-eq:LS-projected-grad}
D(\mathcal J_g\circ\mathsf{ret}_{\mathcal U})[w_{x,\varepsilon}][\phi]=0
\qquad(\phi\in X).
\end{equation}
It satisfies
\begin{equation}\label{I-eq:LS-w-bound-general}
\|w_{x,\varepsilon}\|_{H^1(M)}\le C(M,g)\varepsilon.
\end{equation}
If additionally
\begin{equation}\label{I-eq:LS-refined-hyp}
H_g(x)=0,
\end{equation}
then
\begin{equation}\label{I-eq:LS-w-bound-refined}
\|w_{x,\varepsilon}\|_{H^1(M)}
\le C(M,g)\bigl(\varepsilon|\mathring{\mathrm{II}}(x)|+\varepsilon^2+\omega(R(\varepsilon))\bigr).
\end{equation}
At an umbilic center this is $O(\varepsilon^2)$ for $n\ge8$ and $o(\varepsilon)$ for $n=5,6,7$.

The Hessian $D^2(\mathcal J_g\circ\mathsf{ret}_{\mathcal U})(0)$ is uniformly coercive on $X$ and its operator has uniformly bounded inverse. The solution is $C^1$ in $(x,\varepsilon)$ for positive scales, using Remark~\ref{I-rem:banach-trivialization}. On a boundary region with $H_g\equiv0$, the reduced functional has the second-order expansion of Proposition~\ref{I-rem:reduced-vs-bare}.
\end{lemma}

\begin{proof}
The boundary coordinate expansion gives
\[
\|U_{+,x,\varepsilon}\|_q^q
=\|U_+\|_q^q+O(\varepsilon^2+R^{-(n-1)}).
\]
The weighted trace moment is finite for $n\ge5$, and the cutoff $R(\varepsilon)=\varepsilon^{-3/4}$ has $R^{-(n-1)}=O(\varepsilon^2)$. Hence
\begin{equation}\label{I-eq:t-asymp}
t(x,\varepsilon)=t_\ast+O(\varepsilon^2),\qquad 0<c\le t(x,\varepsilon)\le C.
\end{equation}
The supports lie in the fixed collar because their radius is $2\varepsilon^{1/4}$.

Homogeneity gives $\mathcal J_g(\mathsf{ret}_{\mathcal U}h)=\mathcal J_g(\mathcal U+h)$. For tangent $h,\phi$, direct differentiation gives, also at an approximate bubble,
\begin{equation}\label{I-eq:Jacobi-form}
D^2(\mathcal J_g\circ\mathsf{ret}_{\mathcal U})(0)[h,\phi]
=2\mathcal N_g^\circ(h,\phi)-2(q-1)\mathcal J_g[\mathcal U]
\int_{\partial M}\mathcal U^{q-2}h\phi.
\end{equation}
To transfer Theorem~\ref{I-thm:intrinsic-single-gap} to the normalized complement, put $D_0=\varepsilon\partial_\varepsilon$ and $D_a=\varepsilon\partial_{x^a}$ for $1\le a\le n-1$, and write $U_0=U_{+,x,\varepsilon}$ and $Z_j=D_j(tU_0)$ for the weighted parameter fields. If $h\in X$, then
\[
\langle h,D_jU_0\rangle_E
=-(D_j\log t)\langle h,U_0\rangle_E,
\qquad |\langle h,U_0\rangle_E|\le C\varepsilon\|h\|_E.
\]
The second estimate follows from tangency and Lemma~\ref{I-lem:first-order-drift}. The weighted derivatives of $t$ are bounded. Use the full frame consisting of the common constraint normal and the Jacobi fields of $U_0$. The energy projection determined by its uniformly invertible Gram matrix maps $h$ into the complement associated with $U_0$ with error at most $C\varepsilon\|h\|_E$, preserving tangency. Since
\[
D^2\mathcal J_g[tU_0]=t^{-2}D^2\mathcal J_g[U_0],
\]
the uniformly bounded Hessian loses at most $C\varepsilon\|h\|_E^2$, so its positive gap persists. Use the energy Riesz map to represent this Hessian and the complementary differential by $\mathcal L$ and $\mathcal R_{\mathcal U}$. Define
\[
\mathcal N_{\mathcal U}(w)
=R_E\bigl(D_XF(w)-D_XF(0)-D_X^2F(0)w\bigr),
\qquad F=\mathcal J_g\circ\mathsf{ret}_{\mathcal U}|_X.
\]
The projected equation is
\[
\mathcal Lw=-\mathcal R_{\mathcal U}-\mathcal N_{\mathcal U}(w).
\]
Lemma~\ref{I-lem:first-order-drift} gives
\begin{equation}\label{I-eq:LS-general-residual-bound-referenced}
\|\mathcal R_{\mathcal U}\|_{H^1(M)}\le C(M,g)\varepsilon;
\end{equation}
Lemma~\ref{I-lem:residual-kernel-pure} gives the refined bound when $H_g(x)=0$.

Since $f(u)=|u|^{q-2}u$ has a $\gamma$-H\"older derivative, the trace inequality and the uniformly bounded normalization factors imply
\begin{align*}
\|\mathcal N_{\mathcal U}(w)\|_E&\le C\|w\|_{H^1}^{1+\gamma},\\
\|\mathcal N_{\mathcal U}(w)-\mathcal N_{\mathcal U}(v)\|_E
&\le C(\|w\|_{H^1}+\|v\|_{H^1})^\gamma\|w-v\|_{H^1}.
\end{align*}
The multiplier terms obey the same bounds. Uniform $C^1$ inversion (Lemma~\ref{I-lem:LS-C1-inversion}) therefore gives a unique solution in a fixed small ball and
\[
\|w\|_{H^1}\le C\|\mathcal R_{\mathcal U}\|_E.
\]
This proves both estimates for the correction and the parameter regularity. The energy expansion follows from Lemma~\ref{II-lem:holder-quadratic-expansion}; cutoff contributions are $O(\omega(R)^2+\varepsilon\omega(R))=o(\varepsilon^2)$.
\end{proof}

\begin{definition}[Reduced functional]\label{I-def:reduced}
For the correction of Lemma~\ref{I-lem:LS}, put $\Phi(x,\varepsilon)=\mathsf{ret}_{\mathcal U_{x,\varepsilon}}(w_{x,\varepsilon})\in\mathcal S$. The reduced functional is
\begin{equation}\label{I-eq:Jk}
\mathcal J_1(x,\varepsilon)
\coloneqq  \mathcal J_g\big(\Phi(x,\varepsilon)\big)
= \mathcal J_g\big(\mathcal U_{x,\varepsilon}+w_{x,\varepsilon}\big),
\end{equation}
the second equality following from homogeneity of the quotient.
\end{definition}

\subsection{The reduced quadratic coefficient}

\begin{lemma}[Resolvent limit for the degree-$2$ source]
\label{I-lem:resolvent-limit}
Assume $n\ge5$, \textup{(Y$^+_\partial$)}, \textup{(BG$^{3+}$)}, and $H_g\equiv0$ near $x$; write $S\coloneqq \mathring{\mathrm{II}}(x)$. Let $\mathcal B_{x,\varepsilon}$ be the pulled-back Jacobi form on $X_{x,\varepsilon}$, normalized by $D^2(\widehat{\mathcal J}\circ\mathsf{ret}_{\mathcal U_{x,\varepsilon}})(0)[v,\phi] =2\mathcal B_{x,\varepsilon}(v,\phi)$ for $v,\phi\in X_{x,\varepsilon}$, where $\widehat{\mathcal J}=\mathcal J_g/S_\ast$, with Hessian operator $\mathcal L_{x,\varepsilon}=2\mathcal B_{x,\varepsilon}$. Write $\mathcal L_{\mathrm{flat}}=2\mathcal B_n$ on $X^\flat$. Model directions are rescaled without an amplitude factor, so the limiting form of $\mathcal B_{x,\varepsilon}$ is $S_\ast^{-1}\mathcal B_n$. Define the degree-$2$ source $F_1(x)$ by
\begin{equation}\label{I-eq:F1-definition}
\langle F_1(x),\phi\rangle
\coloneqq -2\,\mathscr L_{S}(\phi)
=4\int_{\R^n_+}y_n\,S^{ab}\,\partial_aU_+\,\partial_b\phi\,dy ,
\end{equation}
with $\mathscr L_S$ as in \eqref{I-eq:spin2-source-functional}, and set $\widehat F_1(x)\coloneqq t_\ast^{-1}F_1(x)=S_\ast^{-1/2}F_1(x)$, where $t_\ast=S_\ast^{1/2}$. We use the same notation for the corresponding covector on $X_{x,\varepsilon}$ under the rescaled Fermi identification. Then:
\begin{enumerate}[label=\textup{(\alph*)},leftmargin=2em]
\item $F_1(x)\in\bigl(\dot H^1(\R^n_+)\bigr)^\ast$ for $n\ge5$, and the differential of the
normalized functional $\widehat{\mathcal J}\coloneqq \mathcal J_g/S_\ast$, restricted to the complement and denoted by $\widehat R_{x,\varepsilon}$, satisfies $\|\varepsilon^{-1}\widehat R_{x,\varepsilon}-\widehat F_1(x)\|_{X_{x,\varepsilon}^\ast}\to0$, uniformly for $x$ in compact subsets of a boundary region on which $H_g\equiv0$.
\item $\langle \widehat F_1(x),\mathcal L_{x,\varepsilon}^{-1}\widehat F_1(x)\rangle
\to\langle F_1(x),\mathcal L_{\mathrm{flat}}^{-1}F_1(x)\rangle
=2\,\delta\kappa_3^{\mathrm{LS}}(n)|S|^2$, uniformly on compact subsets.
\end{enumerate}
\end{lemma}

\begin{proof}
For (a), let $f_1^b=4tS_{ab}\partial_aU_+$ and $f_1^n=0$. Then $|f_1(y)|\le C(1+|y|)^{2-n}$, so $f_1\in L^2$ for $n\ge5$ and $F_1=-\operatorname{div}f_1$ belongs to $(\mathcal D^{1,2})^\ast$. Its tail satisfies
\begin{equation}\label{I-eq:F1-tail}
\sup_{\|\nabla\phi\|_2\le1}\left|\int_{\R^n_+\setminus B_R^+}f_1\cdot\nabla\phi\right|
\le\|f_1\|_{L^2(\R^n_+\setminus B_R^+)}=:\eta(R)\longrightarrow0.
\end{equation}
Lemma~\ref{I-lem:spin2-profile} gives the error in the first-variation expansion $C(\eta_n(\varepsilon)+\omega(R(\varepsilon)))=o(\varepsilon)$. Equation~\eqref{I-eq:F1-tail} also removes the localization of the source, proving (a).

For (b), solve
\begin{equation}\label{I-eq:resolvent-eq}
2\mathcal B_{x,\varepsilon}(v_\varepsilon,\phi)
=\langle\widehat F_1(x),\phi\rangle\qquad(\phi\in X_{x,\varepsilon}).
\end{equation}
Uniform coercivity bounds $v_\varepsilon$ in $H^1(M)$. Its Fermi rescaling has bounded gradient norm and bounded $L^{2n/(n-2)}$ norm, as in the proof of Theorem~\ref{I-thm:intrinsic-single-gap}. Extract a subsequence converging weakly on every half-ball of integer radius to $V\in\mathcal D^{1,2}(\R^n_+)$. Let $L$ be a fixed test radius and $R_c(\varepsilon)=\varepsilon^{-3/4}$ the moving cutoff. Eventually $R_c(\varepsilon)>2L$, so local weak convergence and \eqref{I-eq:F1-tail} give
\[
\limsup_{\varepsilon\to0}
|\langle\widehat F_1,v_\varepsilon\rangle-\langle\widehat F_1,V\rangle|
\le C\|f_1\|_{L^2(|Y|>L)}
\le CL^{-(n-4)/2}.
\]
Letting $L\to\infty$ proves convergence of the pairing. The normalized translation and dilation fields converge strongly to the listed $Z_j^\flat$ after rescaling. The energy Riesz representative of the boundary constraint converges by Remark~\ref{B-rem:inverse-Gram}. Their orthogonality conditions pass to $V$, so $V\in X^\flat$. For $\phi\in X^\flat\cap C_c^\infty(\R^n_+)$, subtract an $o(1)$ finite-dimensional projection from its Fermi rescaling to make it belong to $X_{x,\varepsilon}$. On its fixed support the forms converge to $S_\ast^{-1}\mathcal B_n$. Passing to the limit gives
\[
2S_\ast^{-1}\mathcal B_n(V,\phi)=\langle\widehat F_1(x),\phi\rangle.
\]
Compactly supported tests in the complement are dense. To see this, approximate in the Dirichlet norm and subtract a finite linear combination of fixed compactly supported functions to make the defining moments vanish. The flat gap therefore identifies the unique limit
\[
V=S_\ast\mathcal L_{\mathrm{flat}}^{-1}\widehat F_1
=t_\ast\mathcal L_{\mathrm{flat}}^{-1}F_1,
\qquad \mathcal L_{\mathrm{flat}}=2\mathcal B_n.
\]
Consequently the normalization factors cancel, and
\[
\langle\widehat F_1,V\rangle
=\langle F_1,\mathcal L_{\mathrm{flat}}^{-1}F_1\rangle=2\mathscr L_S(\mathcal B_n^{-1}\mathscr L_S)
=2\delta\kappa_3^{\mathrm{LS}}(n)|S|^2.
\]
Uniqueness makes every subsequence have this limit. The same argument for $x_j\to x_\infty$ in convergent Fermi frames proves uniformity on compact sets.
\end{proof}

\begin{proposition}[Reduced second-order coefficient and Lyapunov--Schmidt correction]\label{I-rem:reduced-vs-bare}
Assume $n\ge5$, \textup{(Y$^+_{\partial}$)}, and \textup{(BG$^{3+}$)}, and use the reduction of Lemma~\ref{I-lem:LS} and Definition~\ref{I-def:reduced}.

Suppose $H_g\equiv0$ on a boundary region. Let $\delta\kappa_3^{\mathrm{LS}}(n)>0$ be the Lyapunov--Schmidt correction on the degree-$2$ subspace defined in
\eqref{I-eq:spin2-schur-general} and set
$\kappa_3^{\mathrm{red}}\coloneqq \kappa_3^{\mathrm{bare}}-\delta\kappa_3^{\mathrm{LS}}$. Then
\begin{equation}\label{I-eq:Rg-red-decomp}
\frac{\mathcal J_1(x,\varepsilon)}{S_\ast}
=1+\mathfrak R_g^{\mathrm{red}}(x)\varepsilon^2+o(\varepsilon^2),
\qquad
\mathfrak R_g^{\mathrm{red}}(x)
=\big(\kappa_3^{\mathrm{bare}}(n)-\delta\kappa_3^{\mathrm{LS}}(n)\big)|\mathring{\mathrm{II}}(x)|^2
\end{equation}
uniformly for $x$ on compact subsets of that region, by Proposition~\ref{I-prop:standard-bubble-inputs}. If the ambient metric varies in a $C^m$ family with $m\ge3$, this coefficient depends $C^{m-3}$-smoothly on $g$.
\end{proposition}

\begin{proof}
Put $A=\mathring{\mathrm{II}}(x)$ and use the normalized functional $\widehat{\mathcal J}=\mathcal J_g/S_\ast$. Lemma~\ref{II-lem:holder-quadratic-expansion} gives
\[
\frac{\mathcal J_1}{S_\ast}
=\widehat{\mathcal J}(\mathcal U)
-\frac12\langle\widehat R_{x,\varepsilon},\mathcal L_{x,\varepsilon}^{-1}\widehat R_{x,\varepsilon}\rangle
+O(\varepsilon^{2+\gamma}),
\]
where $\widehat R$ and $\mathcal L$ are the differential and Hessian of the same normalized functional, restricted to the complement. Lemma~\ref{I-lem:resolvent-limit} and uniform inversion give
\begin{align*}
\widehat R_{x,\varepsilon}&=\varepsilon\widehat F_1(x)+o_{X^\ast}(\varepsilon),\\
\langle\widehat R_{x,\varepsilon},\mathcal L_{x,\varepsilon}^{-1}\widehat R_{x,\varepsilon}\rangle
&=2\varepsilon^2\delta\kappa_3^{\mathrm{LS}}(n)|A|^2+o(\varepsilon^2).
\end{align*}
The factor $1/2$ therefore subtracts exactly $\delta\kappa_3^{\mathrm{LS}}(n)|A|^2$ from the standard-bubble coefficient. Proposition~\ref{I-prop:standard-bubble-inputs} gives the standard-bubble term and proves \eqref{I-eq:Rg-red-decomp}. Uniformity follows from the uniform source and resolvent limits. The displayed formula gives the asserted metric regularity. Finally, Cauchy--Schwarz in $\mathcal B_n$ gives
\begin{equation}\label{I-eq:reduced-expansion-clean}
\frac{\mathscr L_S(\phi)^2}{\mathcal B_n(\phi,\phi)}
\le\delta\kappa_3^{\mathrm{LS}}(n)|S|^2
\qquad(0\ne\phi\in X^\flat),
\end{equation}
so the exact coefficient also implies every upper bound obtained from a fixed test function.
\end{proof}

\subsection{Exact evaluation of the quadratic coefficient}
\label{I-subsec:quadratic-exact}

\begin{convention}[Coefficients before and after reduction]\label{I-conv:bare-vs-reduced}
In the representative with minimal boundary, $\mathfrak R_g^{\mathrm{bare}}=\kappa_3^{\mathrm{bare}}|\mathring{\mathrm{II}}|^2$ is the standard-bubble coefficient, while $\mathfrak R_g^{\mathrm{red}}$ is the corrected coefficient in \eqref{I-eq:Rg-red-decomp}.
\end{convention}

We evaluate the quadratic form in \eqref{I-eq:spin2-schur-general}. Fix $S\in\mathrm{Sym}^2_0(\R^{n-1})$ with $|S|=1$, and write
\[
a=1+t,\qquad s=|y'|^2+a^2,\qquad Q=Q_S(y').
\]
Since $(q-1)\lambda_{\mathrm{flat}}U_+^{q-2}=n/(1+|y'|^2)$ on the boundary, the boundary operator associated with $\mathcal B_n$ is
\[
\mathcal Rv=-\partial_t v-\frac{n}{1+|y'|^2}v.
\]
All functions in the following calculation have tangential spherical harmonic degree two, so the constraint and Jacobi orthogonality conditions hold automatically.

\begin{lemma}[A harmonic Robin pairing]\label{I-lem:Robin-pairing-exact}
For $n\ge5$, put $b_0=Q(y')(1+|y'|^2)^{-n/2}$. There is a unique harmonic $h_0\in X_{2,S}$ satisfying $\mathcal Rh_0=b_0$ weakly. Moreover,
\begin{equation}\label{I-eq:Robin-pairing-exact}
R_n\coloneqq \int_{\R^{n-1}}b_0h_0\,dy'
=\frac{n(n^2-2n-5)}
{(n-2)(n-1)^2(n+1)(n-3)(n-4)A_n^2}.
\end{equation}
\end{lemma}

\begin{proof}
The datum $b_0$ belongs to the dual trace space for $n>4$, since $b_0(y')=O(|y'|^{2-n})$. Coercivity of $\mathcal B_n$ on $X_{2,S}$ gives existence and uniqueness. The equation holds against all Dirichlet tests by angular orthogonality.

Under the Cayley map to the unit ball,
\[
\zeta=\frac{2(y+e_n)}s-e_n,\qquad
h_0(y)=s^{-(n-2)/2}H(\zeta),
\]
$e_n$ is the final Euclidean coordinate vector. The function $H$ is harmonic and its boundary condition is
\[
(\partial_{\nu_{\mathrm{out}}}-1)H=f,
\qquad f(\zeta)=\frac{Q(\zeta')}{2(1+\zeta_n)^2}.
\]
The change of variables gives $R_n=2^{-(n-2)}\int_{\Sph^{n-1}}fH$. Put $z=\zeta_n$ and $\zeta'=\sqrt{1-z^2}\theta$. Let $C_j^{(n+2)/2}$ be the Gegenbauer polynomial of degree $j$. The relevant spherical harmonics are
\[
(1-z^2)C_j^{(n+2)/2}(z)Q(\theta),\qquad j\ge0.
\]
Their degrees are $j+2$, so the harmonic Robin inverse has eigenvalue $(j+1)^{-1}$ on these modes. For general scalar-flat boundary pairings, see \cite[Proposition~3.3]{GongKimMussoWei2026}.

Set $\alpha=(n-1)/2$, and use $(b)_j=\Gamma(b+j)/\Gamma(b)$ for the rising factorial. The angular identity and the one-dimensional integrals are
\[
\int_{\Sph^{n-2}}Q(\theta)^2\,d\theta
=\frac{2|\Sph^{n-2}|}{(n-1)(n+1)},
\]
\[
I_j=\int_{-1}^1(1-z^2)^{(n+1)/2}(1+z)^{-2}
C_j^{(n+2)/2}(z)\,dz,\quad N_j=\int_{-1}^1(1-z^2)^{(n+1)/2}
\bigl(C_j^{(n+2)/2}(z)\bigr)^2\,dz.
\]
The Gegenbauer norm formula and its Rodrigues formula, followed by integration by parts, give
\[
\frac{I_j}{I_0}=(-1)^j\frac{(n+2)_j(2)_j}{j!(n+1)_j},
\qquad
\frac{N_j}{N_0}=\frac{(n+2)_j}{j!}\frac{n+2}{2j+n+2},
\]
where
\[
I_0=2^n\mathrm B(\alpha,\alpha+2),\qquad
N_0=2^{n+2}\mathrm B(\alpha+2,\alpha+2).
\]
Here $\mathrm B$ is the beta function. Thus
\begin{equation}\label{I-eq:Robin-pairing-series}
\begin{aligned}
R_n={}&\frac{|\Sph^{n-2}|}{2(n-1)(n+1)}
\frac{\mathrm B(\alpha,\alpha+2)^2}{\mathrm B(\alpha+2,\alpha+2)}\cdot\sum_{j=0}^\infty
\frac{(n+2)_j(2)_j^2}{j!(n+1)_j^2}
\frac{2j+n+2}{(n+2)(j+1)}.
\end{aligned}
\end{equation}
The summand simplifies to
\[
\frac{n!}{(n+1)(n+2)}
\frac{(j+n+1)(2j+n+2)}{(j+2)(j+3)\cdots(j+n)}.
\]
Using $\bigl(\prod_{k=2}^n(j+k)\bigr)^{-1}=\mathrm B(j+2,n-1)/(n-2)!$ and summing the geometric series and its first two derivatives gives
\[
\sum_{j=0}^\infty\frac{(j+n+1)(2j+n+2)}{\prod_{k=2}^n(j+k)}
=\frac{n(n^2-2n-5)}{(n-2)!(n-1)(n-3)(n-4)}.
\]
For $n>4$, convergence of the positive energy series gives the identity by approximation with finite spherical-harmonic sums. Substituting the last identity in \eqref{I-eq:Robin-pairing-series}, and using the normalization of $A_n$, proves \eqref{I-eq:Robin-pairing-exact}.
\end{proof}

We shall use the following consequence of Green's formula. Suppose $p$ has tangential spherical harmonic degree two, solves $-\Delta p=f$ in the half-space, and has Robin datum $b=\mathcal Rp$. Let $h_b$ be the harmonic solution with Robin datum $b$. Then $\mathcal B_n^{-1}f=p-h_b$, and
\begin{equation}\label{I-eq:Robin-Green-pairing}
f(\mathcal B_n^{-1}f)=f(p)-\int_{\R^{n-1}}bp
+\int_{\R^{n-1}}bh_b.
\end{equation}
Indeed, symmetry gives $\mathcal B_n(p,h_b)=\int bp$, whereas the equation for $p$ gives $\mathcal B_n(p,h_b)=f(h_b)+\int bh_b$. These identities are valid by density in the finite-energy space.

\begin{proposition}[Exact quadratic coefficients]
\label{I-prop:quadratic-exact}\label{I-prop:dkappa5-lower-bound}
For every $n\ge5$,
\begin{align}
\delta\kappa_3^{\mathrm{LS}}(n)
&=\frac{(2n-1)(n+2)}
{2(n-1)^2(n+1)(n-3)(n-4)},\label{I-eq:dkappa-exact}\\
\kappa_3^{\mathrm{red}}(n)
&=-\frac{(n-2)(n^2-2n-2)}
{2(n-1)^2(n+1)(n-3)(n-4)}<0.
\label{I-eq:kappa-red-exact}
\end{align}
In particular, $\delta\kappa_3^{\mathrm{LS}}(5)=21/128$ and $\kappa_3^{\mathrm{red}}(5)=-13/128$.
\end{proposition}

\begin{proof}
For $|S|=1$, let $K_1=(n-2)A_n/2$ and set
\[
p_1=K_1Q(a-2)s^{-n/2}.
\]
Direct differentiation gives
\[
-\Delta p_1=2n(n-2)A_ntQ s^{-(n+2)/2},\qquad
\mathcal Rp_1=-K_1b_0.
\]
Thus $w_S=p_1+K_1h_0$ solves $\mathcal B_n(w_S,\phi)=\mathscr L_S(\phi)$. All these functions have finite Dirichlet energy for $n\ge5$. The radial beta integral
\[
\int_0^\infty\frac{r^{m-1}}{(r^2+a^2)^k}\,dr
=\frac{a^{m-2k}}2\mathrm B(m/2,k-m/2),\qquad 0<m<2k,
\]
and integration over $a\in[1,\infty)$ give
\begin{align*}
&\mathscr L_S(p_1)=\frac{6-n}{2(n-1)(n-3)(n-4)},\\
\int_{\R^{n-1}}(-K_1b_0)p_1&=\frac{n-2}{4(n-1)(n-3)},\quad K_1^2R_n=\frac{n(n-2)(n^2-2n-5)}
{4(n-1)^2(n+1)(n-3)(n-4)}.
\end{align*}
Equation~\eqref{I-eq:Robin-Green-pairing} now yields \eqref{I-eq:dkappa-exact}. Subtracting this value from $\kappa_3^{\mathrm{bare}}(n)$ gives \eqref{I-eq:kappa-red-exact}. Its numerator is positive for $n\ge5$. At $n=5$, this quadratic form equals $\frac1{16}-\frac3{32}+\frac{25}{128}=\frac{21}{128}$.
\end{proof}

\begin{remark}[A dimension-five trial direction]
For $n=5$, the direction $\phi_M=Q_Sts^{-5/2}$ corresponds to Marques' perturbation in \cite[Eq.~(1.3)]{Marques2007}. The same beta integrals give
\[
\frac{\mathscr L_S(\phi_M)^2}{\mathcal B_5(\phi_M,\phi_M)}=\frac1{10},
\qquad
\frac1{16}-\frac1{10}=-\frac3{80}.
\]
The exact complementary correction gives the smaller reduced coefficient $-13/128$.
\end{remark}

The reduced quadratic coefficient is negative in every dimension $n\ge5$. At $n=6$ the quadratic coefficient of the standard bubble vanishes; at an umbilic point with zero mean curvature its cubic coefficient also vanishes. The next weighted moment is logarithmic in that dimension, so the fourth-order analysis below assumes $n\ge7$.

\subsection{Strict inequality for non-umbilic boundaries}

\begin{theorem}[Strict inequality for non-umbilic boundaries]
\label{I-thm:subcriticality-nonumbilic}
Assume $n \ge 5$, \textup{(Y$^+_{\partial}$)}, and \textup{(BG$^{3+}$)}. Suppose there exists a boundary point $p\in\partial M$ with $\mathring{\mathrm{II}}(p) \neq 0$. Then
\[
C^*_{\Esc}(M, g) < S_\ast.
\]
For $n\ge7$ this strict inequality is due to Escobar~\cite{EscobarAnnals92}, for $n=6$ to Escobar~\cite{Escobar1996a}, and for $n=5$ to Marques~\cite{Marques2007}. In particular, the Escobar infimum is attained by a smooth nonnegative minimizer, positive on each connected component on which it is nonzero.
\end{theorem}

\begin{proof}
Choose $\widehat g=\varphi^{4/(n-2)}g$ with $\varphi|_{\partial M}=1$ and $H_{\widehat g}=0$, using \cite[Lemma~3.13]{EscobarGeometry}. Then $\mathring{\mathrm{II}}_{\widehat g}(p)=\mathring{\mathrm{II}}_g(p)\ne0$. Proposition~\ref{I-rem:reduced-vs-bare} gives
\[
\mathcal J_1(p,\varepsilon)/S_\ast
=1+\kappa_3^{\mathrm{red}}(n)|\mathring{\mathrm{II}}(p)|^2\varepsilon^2+o(\varepsilon^2)<1
\]
for small $\varepsilon$, since the reduced coefficient is negative in every dimension $n\ge5$. Conformal invariance transfers this test function to $g$ and proves the strict inequality. Corollary~\ref{I-cor:escobar-precompact-cov} gives attainment. Replace a minimizer by its absolute value. Cherrier's boundedness argument \cite[proof of Theorem~1, pp.~160--165]{Cherrier1984} gives a bounded minimizer. Boundary regularity, the strong maximum principle, and the Hopf boundary lemma give positivity on each nonzero component. The boundary nonlinearity is smooth on the positive range, so elliptic bootstrap gives a smooth minimizer; see also \cite[p.~2]{AlmarazWang2022}.
\end{proof}

\section{Differentiated one-bubble reduction}\label{B-sec:differentiated}
Throughout this section assume $n\ge5$, \textup{(Y$^+_{\partial}$)}, and \textup{(BG$^{3+}$)}.

The scale expansion relates the norm of the constrained gradient to the distances from the standard and corrected bubble families. The complementary gap controls the error transverse to the corrected family; the scale equation supplies the geometric obstruction away from the umbilic locus.

Retain the normalized fields and complement of Lemma~\ref{I-lem:residual-kernel-pure}, and put $P_{x,\varepsilon}=\Pi_{x,\varepsilon}$ and $\omega_\nabla(\varepsilon)=\omega(R(\varepsilon))$.

\subsection{Normalized modulation}
\begin{lemma}[Normalized one-bubble modulation]\label{B-lem:normalized-modulation}
Under the hypotheses of Lemma~\ref{I-lem:modulation}, let $u\in\mathcal S$ be sufficiently close in $H^1(M)$ to a normalized concentrating boundary bubble. Then there are locally unique parameters $(x,\varepsilon)$ and $w\in X_{x,\varepsilon}$ such that
\[
u=\mathsf{ret}_{\mathcal U_{x,\varepsilon}}(w),\qquad
\|w\|_{H^1}\le C\operatorname{dist}_{H^1}(u,\{\mathcal U_{x,\varepsilon}\}).
\]
The parameters depend $C^1$ on $u$ in the local chart.
\end{lemma}
\begin{proof}
Choose a nearby base parameter $(x^0,\varepsilon^0)$ and use
\[
\varepsilon(s)=e^s\varepsilon^0,\qquad
x(\xi)=\exp_{x^0}^{\partial}(\varepsilon^0\xi),\qquad
\Xi(s,\xi)=\mathcal U_{x(\xi),\varepsilon(s)}.
\]
On a fixed parameter ball, Fermi rescaling, the cutoff estimates, and the uniform Gram estimate in the proof of Lemma~\ref{I-lem:modulation} give
\begin{equation}\label{II-eq:uniform-modulation-chart-estimates-onebubble}
\|D\Xi\|+\|D^2\Xi\|\le C,\qquad
\|D\Xi(0)h\|_{H^1}\ge c|h|.
\end{equation}
For $V,u\in\mathcal S$ sufficiently close, put
\[
\nu_V(u)=\frac{u}{\int_{\partial M}|V|^{q-2}Vu\,d\sigma_g}-V.
\]
Then $\nu_V(u)\in T_V\mathcal S$ and $u=\mathsf{ret}_V(\nu_V(u))$. For the weighted modulation frame $Z_\mu$ set
\[
\mathcal M_\mu(u,z)=\langle\nu_{\Xi(z)}(u),Z_\mu(z)\rangle_{E,g}.
\]
At $u=\Xi(0)$, the derivative with respect to $z$ is, up to sign, the Gram matrix of the frame, uniformly invertible by \eqref{II-eq:uniform-modulation-chart-estimates-onebubble}. The quantitative implicit-function theorem gives the decomposition, local uniqueness, and the displayed bound, exactly as in the proof of Lemma~\ref{I-lem:modulation}.
\end{proof}

\subsection{The one-bubble scale derivative}

\begin{remark}[Model rescaling and the complementary inverse]
\label{II-rem:blowup-trivialization-onebubble}
Write $\mathbf p=(x,\varepsilon)$ and $X_+=X^\flat$, let $P_+^\perp$ be the Dirichlet-orthogonal projection onto $X_+$, and put $\mathbf D=\varepsilon\partial_\varepsilon$ at a fixed model coordinate $Y$. The notation $\mathbf D^{\mathrm{ph}}$ denotes the logarithmic scale derivative at a fixed point of $M$. If $\widetilde V_\varepsilon(Y)=\varepsilon^{(n-2)/2}V_\varepsilon(\Phi_x(\varepsilon Y))$, then
\[
\varepsilon^{(n-2)/2}(\mathbf D^{\mathrm{ph}}V_\varepsilon)(\Phi_x(\varepsilon Y))
=\mathbf D\widetilde V_\varepsilon-\mathcal A\widetilde V_\varepsilon,
\qquad \mathcal A=\frac{n-2}{2}+Y\cdot\nabla.
\]
Let $\mathcal T_{x,\varepsilon}$ denote critical Fermi rescaling, and define the bounded model map
\[
\mathfrak B_{\mathbf p}z
=P_{\mathbf p}\mathcal T_{x,\varepsilon}(\chi_{R(\varepsilon)}z),
\qquad z\in X_+.
\]
For $R=\varepsilon^{-\beta}$ and $z,\mathcal Az\in\mathcal D^{1,2}$, direct differentiation gives
\begin{equation}\label{II-eq:fixed-model-rescaling-derivative}
\mathbf D^{\mathrm{ph}}\mathcal T_{x,\varepsilon}(\chi_Rz)
=\mathcal T_{x,\varepsilon}\bigl(-\chi_R\mathcal Az
+(\beta-1)(Y\cdot\nabla\chi_R)z\bigr).
\end{equation}
Hardy's inequality makes the annular term tend to zero in the Dirichlet norm. For the $\beta$-cutoff family use $R_\beta(\varepsilon)$ in the definition of $\mathfrak B$. For each fixed $z\in X_+$, the energy norm of $\mathfrak B_{\mathbf p}z$ tends to $\|\nabla z\|_2$ by coordinate comparison and the inverse Gram formula. The limiting complementary Hessian, with model directions rescaled without an amplitude factor, is
\[
\langle\mathscr L_0z,\phi\rangle
=2\mathcal B_n(z,\phi),\qquad
\lambda_{\mathrm{flat}}=S_\ast^{q/2}.
\]
The hemisphere gap gives
\begin{equation}\label{II-eq:flat-positive-gap-onebubble}
\mathcal B_n(z,z)\ge c_{\mathrm{flat}}\|\nabla z\|_2^2,
\qquad z\in X_+.
\end{equation}

We also use a complementary inverse on the fixed space $\mathcal H=H^1(M)$, identified with its dual by the energy form. If $\mathcal L_{\mathbf p}^{\mathrm{ph}}$ is the complementary Hessian of $\mathcal J_g\circ\mathsf{ret}_{\mathcal U_{\mathbf p}}$ at zero, put
\begin{equation}\label{II-eq:ambient-stabilized-operator-onebubble}
\mathbb L_{\mathbf p}
=P_{\mathbf p}\mathcal L_{\mathbf p}^{\mathrm{ph}}P_{\mathbf p}
+(I-P_{\mathbf p}),\qquad
\mathbb G_{\mathbf p}=\mathbb L_{\mathbf p}^{-1}.
\end{equation}
Uniform complementary coercivity and differentiation of the inverse give
\begin{equation}\label{II-eq:ambient-resolvent-bound-onebubble}
\|\mathbb G_{\mathbf p}\|\le C,
\end{equation}
\begin{equation}\label{II-eq:ambient-resolvent-derivative-onebubble}
\partial_\theta\mathbb G_{\mathbf p}
=-\mathbb G_{\mathbf p}(\partial_\theta\mathbb L_{\mathbf p})\mathbb G_{\mathbf p}.
\end{equation}
Let $\mathcal S_{\mathbf p}^{\mathrm{ph}}\in X_{\mathbf p}$ be the energy Riesz representative of the complementary differential. For a local identification $\mathcal T_{\mathbf p}:X_0\to X_{\mathbf p}$, put
$\mathcal L_{\mathbf p}=\mathcal T_{\mathbf p}^*\mathcal L_{\mathbf p}^{\mathrm{ph}}\mathcal T_{\mathbf p}$ and
$\mathcal S_{\mathbf p}=\mathcal T_{\mathbf p}^*\mathcal S_{\mathbf p}^{\mathrm{ph}}$. Then
\begin{equation}\label{II-eq:fixed-physical-quadratic-invariance-onebubble}
\langle\mathcal L_{\mathbf p}^{-1}\mathcal S_{\mathbf p},\mathcal S_{\mathbf p}\rangle
=\langle\mathbb G_{\mathbf p}\mathcal S_{\mathbf p}^{\mathrm{ph}},\mathcal S_{\mathbf p}^{\mathrm{ph}}\rangle.
\end{equation}
\end{remark}

\begin{lemma}[Differentiation after completing the square]
\label{II-lem:differentiated-quadratic-completion}
Let $a_\varepsilon$ be a $C^1$ family of symmetric bilinear forms on a fixed Hilbert space, with
$\|a_\varepsilon\|+\|\mathbf D^{\mathrm{ph}}a_\varepsilon\|\le C$, where
$\mathbf D^{\mathrm{ph}}=\varepsilon\partial_\varepsilon$. Let $X_\varepsilon$ be closed subspaces, and suppose $C^1$ families $v_\varepsilon,V_\varepsilon\in X_\varepsilon$ satisfy
\[
a_\varepsilon(v_\varepsilon,h)=s_\varepsilon[h]
\quad(h\in X_\varepsilon),\qquad
e_\varepsilon=v_\varepsilon-V_\varepsilon.
\]
Then
\begin{equation}\label{II-eq:quadratic-completion-identity}
s_\varepsilon[v_\varepsilon]
=2s_\varepsilon[V_\varepsilon]
-a_\varepsilon(V_\varepsilon,V_\varepsilon)
+a_\varepsilon(e_\varepsilon,e_\varepsilon).
\end{equation}
If $\|e_\varepsilon\|=o(\varepsilon^k)$ and
$\|\mathbf D^{\mathrm{ph}}e_\varepsilon\|=O(\varepsilon^k)$, then the last term and its weighted derivative are $o(\varepsilon^{2k})$.
\end{lemma}
\begin{proof}
Expand the last term and use the equation with $h=v_\varepsilon,V_\varepsilon$. Moreover,
\[
\bigl|\mathbf D^{\mathrm{ph}}a_\varepsilon(e_\varepsilon,e_\varepsilon)\bigr|
\le C\|e_\varepsilon\|^2
+2C\|e_\varepsilon\|\,
\|\mathbf D^{\mathrm{ph}}e_\varepsilon\|
=o(\varepsilon^{2k}).
\]
These derivatives are taken in the fixed ambient Hilbert space on which the forms are defined.
\end{proof}

\begin{proposition}[Scale derivative of the one-bubble reduced functional]
\label{II-prop:onebubble-scale-derivative-II}
Assume $n\ge5$, \textup{(Y$^+_{\partial}$)}, and \textup{(BG$^{3+}$)}, and let $\mathcal C\subset\partial M$ be an open boundary region admitting a fixed Fermi collar, with $H_g=0$ on $\mathcal C$. For
\[
\mathcal J_1(x,\varepsilon)
\coloneqq \mathcal J_g\bigl(
\mathsf{ret}_{\mathcal U_{x,\varepsilon}}(w_{x,\varepsilon})
\bigr),
\qquad x\in\mathcal C,
\]
one has, uniformly on compact subsets of $\mathcal C$,
\begin{equation}\label{II-eq:onebubble-scale-derivative-II}
\frac1{S_\ast}\partial_\varepsilon\mathcal J_1(x,\varepsilon)
=2\mathfrak R_g^{\mathrm{red}}(x)\varepsilon+o(\varepsilon).
\end{equation}
\end{proposition}

\begin{proof}
Put $\alpha=(n-2)/2$, $\gamma=q-2$, and $\mathcal U=\mathcal U_{x,\varepsilon}$, and use the first variation and Hessian restricted to the complement on $M$
\begin{equation}\label{II-eq:onebubble-physical-data}
\mathscr S=D(\mathcal J_g\circ\mathsf{ret}_{\mathcal U})(0)|_{X_{x,\varepsilon}},
\qquad
\mathscr L=D^2(\mathcal J_g\circ\mathsf{ret}_{\mathcal U})(0)|_{X_{x,\varepsilon}}.
\end{equation}
The Hessian on the complement is the restriction of the symmetric form
\begin{equation}\label{B-eq:physical-PDE-form}
a_{x,\varepsilon}(h,\phi)
=2\mathcal N_g^\circ(h,\phi)
-2(q-1)\mathcal J_g[\mathcal U]
\int_{\partial M}\mathcal U^{q-2}h\phi\,d\sigma_g.
\end{equation}
The first model source is
$\mathscr S_{1,x}=t_\ast F_1(x)=-2t_\ast\mathscr L_{\mathring{\mathrm{II}}(x)}$.
Its degree-$2$ angular dependence gives the full moment cancellations
\begin{equation}\label{II-eq:first-source-moment-cancellation}
\mathscr S_{1,x}[U_+]=\mathscr S_{1,x}[Z_j^\flat]=0
\qquad(0\le j\le n-1).
\end{equation}
The differentiated standard-bubble expansion gives
\begin{equation}\label{II-eq:uncorrected-scale-derivative-II}
S_\ast^{-1}\partial_\varepsilon\mathcal J_g[\mathcal U]
=2\mathfrak R_g^{\mathrm{bare}}(x)\varepsilon+o(\varepsilon).
\end{equation}
This is \eqref{I-eq:standard-bubble-scale-input}; its proof estimates the moving cutoff and the geometric remainder before differentiating.

Let $v=\mathscr L^{-1}\mathscr S$ and
$z_{1,x}=\mathscr L_0^{-1}\mathscr S_{1,x}$. Write
$\mathscr S^{\mathrm{full}}=D\mathcal J_g[\mathcal U]$ on $H^1(M)$, whose restriction is $\mathscr S$. We first estimate derivatives at a fixed point of $M$. If $\widetilde h_\varepsilon(Y)=\varepsilon^\alpha h(\Phi_x(\varepsilon Y))$ and
$s_\varepsilon[h]=\varepsilon^k\int f(Y)\widetilde h_\varepsilon(Y)\,dY$, then integration by parts gives
\begin{equation}\label{II-eq:physical-source-differentiation}
\mathbf D^{\mathrm{ph}}s_\varepsilon[h]
=\varepsilon^k\int\bigl((k+\alpha-n)f-Y\cdot\nabla f\bigr)
\widetilde h_\varepsilon\,dY.
\end{equation}
The boundary contribution is zero because $Y_n=0$ there. For the first source, tangential integration by parts gives a scalar forcing $f$ with
$|f|+|Y\cdot\nabla f|\le C(1+|Y|)^{1-n}$. Hardy's inequality applies since
\[
\int_{|Y|>1}|Y|^2\bigl(|f|^2+|Y\cdot\nabla f|^2\bigr)\,dY
\le C\int_1^\infty r^{3-n}\,dr<\infty.
\]
For the explicit metric and cutoff errors, the same differentiation and integration by parts preserve the bounds $C(\eta_n(\varepsilon)+R^{-\alpha})$ established in Lemma~\ref{I-lem:spin2-profile}. Derivatives of normalization and quotient factors are estimated by their defining integrals. Hence
\begin{equation}\label{II-eq:first-physical-derivative-bounds}
\|\mathscr S^{\mathrm{full}}\|+\|\mathbf D^{\mathrm{ph}}\mathscr S^{\mathrm{full}}\|
\le C\varepsilon,
\qquad
\|a_{x,\varepsilon}\|+\|\mathbf D^{\mathrm{ph}}a_{x,\varepsilon}\|\le C.
\end{equation}
For the form estimate, only the boundary potential
$V_\varepsilon=\mathcal J_g[\mathcal U]\mathcal U^\gamma$ varies. Off the cutoff annulus its logarithmic derivative is bounded by $C\mathcal U^\gamma$; on the annulus, $\chi_R^\gamma=\eta_R^2$ gives the $O(R^{-1})$ coefficient bound in Convention~\ref{II-conv:cutoff-choice}. Apply \eqref{I-eq:boundary-potential-form-bound} and
$\mathbf D\mathcal J_g[\mathcal U]=O(\varepsilon^2)$ from Proposition~\ref{I-prop:standard-bubble-inputs}.

Let $Z$ denote the map whose columns are the normal frame in Remark~\ref{B-rem:inverse-Gram}, and $G=Z^*Z$. Then
\[
P=I-ZG^{-1}Z^*,\qquad
\mathbf D^{\mathrm{ph}}G^{-1}=-G^{-1}(\mathbf D^{\mathrm{ph}}G)G^{-1},
\qquad \|\mathbf D^{\mathrm{ph}}P\|\le C.
\]
The frame and its derivative are bounded by the explicit Jacobi fields and the exact constraint-normal identity from that remark, using \eqref{II-eq:first-physical-derivative-bounds}. The stabilized operator \eqref{II-eq:ambient-stabilized-operator-onebubble} consequently satisfies $\|\mathbb L\|+\|\mathbf D^{\mathrm{ph}}\mathbb L\|\le C$. The derivative of its inverse is given by \eqref{II-eq:ambient-resolvent-derivative-onebubble}, so
\begin{equation}\label{II-eq:first-response-bounded-derivative}
\|v\|_{H^1}+\|\mathbf D^{\mathrm{ph}}v\|_{H^1}\le C\varepsilon.
\end{equation}
Complementary coordinate derivatives obey the same bounds after the local identification of Remark~\ref{I-rem:banach-trivialization}.

For the solution $z_{1,x}$ of the model linearized equation, the interior forcing $f_{1,x}$ satisfies
\[
|f_{1,x}(Y)|\le C(1+|Y|)^{1-n},\qquad
\int_{|Y|>1}|Y|^2|f_{1,x}|^2\,dY<\infty,
\]
and it satisfies a homogeneous Robin boundary condition by angular orthogonality. Write
$A_r=B^+_{2r}\setminus\overline{B^+_r}$ and
$\partial' A=A\cap\partial\mathbb R^n_+$. The unit half-annulus $H^2$ estimate, after dilation, gives for interior forcing $f$ and boundary forcing $g$
\begin{equation}\label{II-eq:weighted-annular-H2}
\begin{aligned}
r\|\nabla^2h\|_{L^2(A_r)}
\le C\bigl(&\|\nabla h\|_{L^2(A_{r/2,4r})}
+r^{-1}\|h\|_{L^2(A_{r/2,4r})}\\
+&r\|f\|_{L^2(A_{r/2,4r})}
+r^{n/2}\|g(r\,\cdot)\|_{H^{1/2}(\partial' A_{1/2,4})}\bigr),
\end{aligned}
\end{equation}
where $A_{a,b}=B_b^+\setminus\overline{B_a^+}$. To obtain this estimate, scale to the unit half-annulus and use tangential difference quotients in the weak equation with a compactly supported cutoff. Trace interpolation absorbs the terms containing the Robin coefficient; the remaining terms are controlled by the displayed $H^1$, forcing, and boundary norms. The equation then recovers the second normal derivative. The rescaled Robin coefficient $nr/(1+r^2|Y'|^2)$ and its first derivatives are uniformly bounded on the boundary half-annulus, so the constant is independent of $r\ge1$. Squaring and summing over dyadic annuli, whose enlarged versions have bounded overlap, and using Hardy's inequality proves
$Y\nabla^2z_{1,x}\in L^2$, and hence $\mathcal A z_{1,x}\in\mathcal D^{1,2}$. This estimate is uniform in $x$ on the stated compact sets.

Set
\[
V_\varepsilon=\varepsilon\mathfrak B_{x,\varepsilon}z_{1,x},
\qquad
C_x=\langle\mathscr S_{1,x},z_{1,x}\rangle.
\]
The pairings satisfy
\begin{align}
\mathscr S[V_\varepsilon]=\varepsilon^2C_x+o(\varepsilon^2),\qquad
a_{x,\varepsilon}(V_\varepsilon,V_\varepsilon)=\varepsilon^2C_x+o(\varepsilon^2),
\label{II-eq:first-fixed-trial-pairings}
\end{align}
with weighted derivatives of both remainders $o(\varepsilon^2)$, and
$\|\mathbf D^{\mathrm{ph}}V_\varepsilon\|_{H^1}=O(\varepsilon)$.
For the source pairing, rescale to the fixed model and differentiate its explicit coefficients. The normalized errors are bounded by
\[
C\left(\frac{\eta_n(\varepsilon)+\omega(R)}{\varepsilon}
+R^{-(n-4)/2}\right)=o(1).
\]
For the form pairing, the principal and scalar errors are $O(\varepsilon R)$ and $O((\varepsilon R)^2)$, including their logarithmic derivatives. In model coordinates the boundary potential tends to $V_0=\lambda_{\mathrm{flat}}U_+^\gamma$. Its logarithmic scale derivative at fixed points of $M$, expressed in rescaled coordinates, tends to $-V_0-Y'\cdot\nabla V_0$. Including the derivatives of both transported test functions gives the leading sum
\[
(-V_0-Y'\cdot\nabla V_0)z^2-2V_0z(\alpha z+Y'\cdot\nabla z)
=-\operatorname{div}_{Y'}(Y'V_0z^2),
\]
whose integral is zero, first for compactly supported tests and then by approximation for $z,\mathcal Az\in\mathcal D^{1,2}$. The analogous Dirichlet terms cancel by dilation invariance. The coefficient errors and cutoff tails therefore tend to zero also in the differentiated fixed-profile pairing.

These assertions are unchanged by the finite-dimensional projection in $\mathfrak B$. To see the differentiated statement, let $T_\varepsilon z=\mathcal T_{x,\varepsilon}(\chi_Rz)$ and use a normal frame $\zeta_j$ in the inverse Gram formula. Its rescaled limit is the flat normal frame, and the logarithmic scale derivative at fixed points of $M$, expressed in rescaled coordinates, tends to $-\mathcal A\zeta_j^\flat$. The energy Riesz representative $n_\varepsilon$ of $h\mapsto\int\mathcal U^{q-1}h$ satisfies the exact identity
\[
n_\varepsilon=\frac{\mathcal U}{\mathcal J_g[\mathcal U]}
-\frac{R_E\mathscr S^{\mathrm{full}}}{2\mathcal J_g[\mathcal U]}.
\]
By \eqref{II-eq:first-physical-derivative-bounds} and \eqref{II-eq:uncorrected-scale-derivative-II},
$n_\varepsilon-\mathcal U/\mathcal J_g[\mathcal U]=O_{H^1}(\varepsilon)$ and
$\mathbf D^{\mathrm{ph}}n_\varepsilon-(\mathbf D^{\mathrm{ph}}\mathcal U)/\mathcal J_g[\mathcal U]=O_{H^1}(\varepsilon)$.
This proves the differentiated limit of the constraint normal; the remaining frame elements are explicit Jacobi fields. The projection coefficients of $T_\varepsilon z_{1,x}$ tend to zero. Their differentiated limits also vanish, since critical dilation is a Dirichlet isometry:
\[
\int\nabla(\mathcal A z_{1,x})\cdot\nabla\zeta_j^\flat
+\int\nabla z_{1,x}\cdot\nabla(\mathcal A\zeta_j^\flat)=0.
\]
Thus the removed projection and its logarithmic scale derivative at fixed points of $M$ are $o_{H^1}(1)$. This proves all the assertions in \eqref{II-eq:first-fixed-trial-pairings}.

Lemma~\ref{I-lem:resolvent-limit} gives
$\mathscr S[v]=\varepsilon^2C_x+o(\varepsilon^2)$. Applying
\eqref{II-eq:quadratic-completion-identity} with the trial function $V_\varepsilon$, and using complementary coercivity, therefore gives
\[
\|v-V_\varepsilon\|_{H^1}=o(\varepsilon).
\]
The logarithmic scale derivatives of both terms at fixed points of $M$ are $O(\varepsilon)$ by the preceding estimates. Lemma~\ref{II-lem:differentiated-quadratic-completion} now yields
\begin{equation}\label{II-eq:first-inverse-pairing-C1}
\langle\mathscr L^{-1}\mathscr S,\mathscr S\rangle
=\varepsilon^2C_x+o(\varepsilon^2),
\qquad
\mathbf D\langle\mathscr L^{-1}\mathscr S,\mathscr S\rangle
=2\varepsilon^2C_x+o(\varepsilon^2).
\end{equation}

It remains to differentiate the nonlinear remainder. Put $\gamma=q-2$ and, for $a\ge0$, define
\begin{align*}
N(a,w)&=|a+w|^\gamma(a+w)-a^{1+\gamma}-(1+\gamma)a^\gamma w,\\
Q(a,w)&=|a+w|^{2+\gamma}-a^{2+\gamma}-(2+\gamma)a^{1+\gamma}w
-\frac{(2+\gamma)(1+\gamma)}2a^\gamma w^2.
\end{align*}
Splitting into the cases $|w|\le a/2$ and $|w|>a/2$ gives
\[
|a\partial_aN(a,w)|\le C|w|^{1+\gamma},\qquad
|a\partial_aQ(a,w)|\le C|w|^{2+\gamma},
\]
with continuous extension at $a=0$. The explicit bubble derivative has the form
$\mathbf D^{\mathrm{ph}}\mathcal U=b(Y)\mathcal U+\mathcal A^{\mathrm{ann}}$,
where $b$ is bounded and the annular term has critical norm $O(R^{-\alpha})$, $\alpha=(n-2)/2$. On the part of the annulus where $\eta>0$, the coefficient satisfies
\begin{equation}\label{II-eq:annular-nonlinear-cancellation}
(\eta^{n-2}U_+)^{\gamma-1}(n-2)\eta^{n-3}(\mathbf D\eta)U_+
=(n-2)\eta(\mathbf D\eta)U_+^\gamma.
\end{equation}
The right-hand side extends continuously across $\eta=0$, and its $L^{n-1}$ norm is $O(R^{-1})$. The remaining annular term is estimated by H\"older's inequality using $\|\mathcal A^{\mathrm{ann}}\|_q=O(R^{-\alpha})$. Thus the parameter derivative of the nonlinear gradient at fixed $w$ is bounded by
\begin{equation}\label{II-eq:nonlinear-gradient-parameter-bound}
C\bigl(\|w\|^{1+\gamma}+R^{-\alpha}\|w\|^\gamma+R^{-1}\|w\|\bigr).
\end{equation} Differentiating the complementary equation and using uniform inversion gives $\|\mathbf D^{\mathrm{ph}}w\|_{H^1}=O(\varepsilon)$. Stationarity and the preceding pointwise estimates then give
\begin{equation}\label{II-eq:onebubble-differentiated-quadratic}
\mathcal J_1=\mathcal J_g[\mathcal U]
-\frac12\langle\mathscr L^{-1}\mathscr S,\mathscr S\rangle+\mathscr E,
\qquad
|\mathscr E|+|\mathbf D\mathscr E|
\le C\left(\varepsilon^q+R^{-\alpha}\varepsilon^{q-1}+R^{-1}\varepsilon^2\right)
=o(\varepsilon^2).
\end{equation}
Here $3(n-2)/8+q-1>2$ for every $n\ge5$.
Finally, the source normalization and \eqref{I-eq:spin2-schur-general} give
$C_x=2S_\ast\delta\kappa_3^{\mathrm{LS}}(n)|\mathring{\mathrm{II}}(x)|^2$.
Differentiate the last expansion, apply \eqref{II-eq:first-inverse-pairing-C1} and \eqref{II-eq:uncorrected-scale-derivative-II}, and subtract this correction from the standard-bubble coefficient. This proves \eqref{II-eq:onebubble-scale-derivative-II}.
\end{proof}

\subsection{Quantitative stability near a boundary bubble}
\label{II-subsec:residual-stability}

Write
\[
r_g(u)=\|\mathsf{grad}_{\mathcal S}\mathcal J_g(u)\|_{E,g},
\qquad
\mathcal B_0=\{\sigma\mathcal U_{x,\varepsilon}:\sigma\in\{\pm1\}\},
\qquad
\mathcal B_{\mathrm{LS}}=\{\sigma\Phi(x,\varepsilon):\sigma\in\{\pm1\}\},
\]
where the scales are sufficiently small and the centers lie in the boundary region under consideration. Distances in this subsection are measured in the energy norm. For $p=(x,\varepsilon)$ put
\[
F_p(v)=\mathcal J_g(\mathsf{ret}_{\mathcal U_p}(v)),\qquad v\in X_p.
\]
The correction $w_p$ of Lemma~\ref{I-lem:LS} satisfies $D_XF_p(w_p)=0$.

\begin{lemma}[Complementary estimates and parameter derivatives]
\label{II-lem:residual-coordinates}
There are uniform constants $c,C,\rho>0$ such that, for sufficiently small scales and $v,w_p\in X_p$ with $\|v\|_E,\|w_p\|_E<\rho$, setting $u=\mathsf{ret}_{\mathcal U_p}(v)$ and $\eta=v-w_p$ gives
\begin{align}
\|\eta\|_E+\|u-\Phi(p)\|_E&\le C r_g(u),\qquad
r_g(\Phi(p))\le C r_g(u),\label{II-eq:residual-complement}\\
c\|\eta\|_E^2&\le F_p(v)-\mathcal J_1(p)\le C\|\eta\|_E^2.
\label{II-eq:residual-normal-energy}
\end{align}
Moreover, in local boundary coordinates,
\begin{equation}\label{II-eq:residual-parameter-equivalence}
r_g(\Phi(p))\asymp
\left(|\varepsilon\partial_\varepsilon\mathcal J_1(p)|^2
+\varepsilon^2|\nabla_\partial\mathcal J_1(p)|^2\right)^{1/2}.
\end{equation}
The constants are uniform on compact subsets of the boundary region. The same assertions hold after multiplication by $-1$.
\end{lemma}

\begin{proof}
The complementary Hessian is uniformly positive by Lemma~\ref{I-lem:LS}. Its uniform continuity on a fixed small ball gives
\[
c\|h\|_E^2\le D^2F_p(z)[h,h]\le C\|h\|_E^2
\qquad(z,h\in X_p,\ \|z\|_E<2\rho).
\]
Integrating along the segment from $w_p$ to $v$ and using $D_XF_p(w_p)=0$, we obtain
\[
c\|\eta\|_E^2\le D_XF_p(v)[\eta]
\le C r_g(u)\|\eta\|_E.
\]
The normalization map has uniformly bounded first and second derivatives here. This proves the first estimate in \eqref{II-eq:residual-complement}. The identity
\[
F_p(v)-F_p(w_p)=\int_0^1(1-s)D^2F_p(w_p+s\eta)[\eta,\eta]\,ds
\]
proves both bounds in \eqref{II-eq:residual-normal-energy}. For the gradient estimate let $n(u)=R_E\ell_u$, with $\ell_u$ as in \eqref{I-eq:exact-PS-multiplier}. Since $\langle n(u),u\rangle_E=1$, one has $\|n(u)\|_E\ge\|u\|_E^{-1}$, and the constrained gradient is exactly
\[
\mathsf{grad}_{\mathcal S}\mathcal J_g(u)
=2\left(u-\frac{n(u)}{\|n(u)\|_E^2}\right).
\]
The trace inequality bounds $Dn$ on the energy-bounded neighborhood. This formula therefore makes the constrained gradient uniformly Lipschitz. Applying this to $u$ and $\Phi(p)$ proves the remaining estimate in \eqref{II-eq:residual-complement}.

For \eqref{II-eq:residual-parameter-equivalence}, let $N(z)=z/\|z\|_q$ and
\[
E_p=DN(\mathcal U_p+w_p)X_p\subset T_{\Phi(p)}\mathcal S.
\]
Complementary stationarity says that the constrained gradient at $\Phi(p)$ is orthogonal to $E_p$. This subspace has codimension $n$ in $T_{\Phi(p)}\mathcal S$. Let $P_p$ be the energy projection onto $X_p$, and let $D_0=\varepsilon\partial_\varepsilon$ and $D_a=\varepsilon\partial_{x^a}$. The normalized modulation estimates and the inverse Gram formula give $\|D_\alpha P_p\|\le C$. Differentiating $P_pw_p=w_p$ yields
\[
(I-P_p)D_\alpha w_p=(D_\alpha P_p)w_p=O_E(\varepsilon).
\]
The term $DN(\mathcal U_p+w_p)P_pD_\alpha w_p$ belongs to $E_p$. Consequently the projection of $D_\alpha\Phi(p)$ onto $E_p^\perp\cap T_{\Phi(p)}\mathcal S$ equals the projection of
\[
DN(\mathcal U_p+w_p)Z_\alpha+O_E(\varepsilon).
\]
Since $\|w_p\|_E=O(\varepsilon)$, these projected vectors differ by $O_E(\varepsilon)$ from the normalized bubble frame. Its Gram matrix has a uniformly positive flat limit by Lemma~\ref{B-lem:normalized-modulation}. They therefore form a uniformly nondegenerate basis of $E_p^\perp\cap T_{\Phi(p)}\mathcal S$. Let $K$ be the energy Gram matrix of these projected vectors and let $b_\alpha=D_\alpha\mathcal J_1$. Pairing with the constrained gradient gives the exact identity
\[
r_g(\Phi(p))^2=b^{\mathsf T}K^{-1}b.
\]
The uniform bounds for $K$ and $K^{-1}$ prove \eqref{II-eq:residual-parameter-equivalence}. Evenness of the functional proves the assertions for both signs.
\end{proof}

\begin{theorem}[Quantitative stability away from the umbilic locus]
\label{II-thm:residual-stability}
Let $K\Subset\mathcal C$, where $H_g=0$ on the open boundary region $\mathcal C$, and assume $\inf_K|\mathring{\mathrm{II}}_g|>0$. There are $\rho,\varepsilon_0,C>0$ such that, whenever normalized modulation gives
\[
u=\sigma\mathsf{ret}_{\mathcal U_{x,\varepsilon}}(v),\qquad
v\in X_{x,\varepsilon},\quad \|v\|_E<\rho,\quad
x\in K,\quad 0<\varepsilon<\varepsilon_0,\quad \sigma\in\{\pm1\},
\]
one has
\begin{align}
\|u-\sigma\Phi(x,\varepsilon)\|_E&\le Cr_g(u),\label{II-eq:residual-corrected-linear}\\
\varepsilon^2&\le Cr_g(u),\label{II-eq:residual-scale}\\
\dist_E(u,\mathcal B_0)&\le Cr_g(u)^{1/2}.
\label{II-eq:residual-standard-sqrt}
\end{align}
In particular, $\dist_E(u,\mathcal B_{\mathrm{LS}})\le Cr_g(u)$. If $\mathcal J_g[u]\ge S_\ast$, then
\begin{equation}\label{II-eq:residual-above-linear}
\dist_E(u,\mathcal B_0)\le Cr_g(u).
\end{equation}
\end{theorem}

\begin{proof}
By evenness we may take $\sigma=1$. The first estimate is Lemma~\ref{II-lem:residual-coordinates}. The exact coefficient in Proposition~\ref{I-prop:quadratic-exact} is negative, so the differentiated expansion gives
\[
c\varepsilon^2\le|\varepsilon\partial_\varepsilon\mathcal J_1(x,\varepsilon)|
\le Cr_g(\Phi(x,\varepsilon))\le Cr_g(u).
\]
This proves \eqref{II-eq:residual-scale}. Since $\|w_{x,\varepsilon}\|_E\le C\varepsilon$,
\[
\dist_E(u,\mathcal B_0)\le C\bigl(r_g(u)+\varepsilon\bigr)
\le Cr_g(u)^{1/2}
\]
for small $r_g(u)$. The remaining part of the neighborhood is covered by enlarging $C$.

On $K$, the expansion of $\mathcal J_1$ gives $\mathcal J_1(x,\varepsilon)\le S_\ast-c\varepsilon^2$. If $\mathcal J_g[u]\ge S_\ast$, the upper estimate in \eqref{II-eq:residual-normal-energy} implies
\[
c\varepsilon^2\le\mathcal J_g[u]-\mathcal J_1(x,\varepsilon)
\le C\|v-w_{x,\varepsilon}\|_E^2\le Cr_g(u)^2.
\]
Thus $\varepsilon\le Cr_g(u)$, which proves \eqref{II-eq:residual-above-linear}.
\end{proof}

\begin{proposition}[Sharpness of the stability exponents]
\label{II-prop:residual-sharpness}
Suppose $H_g=0$ on a connected boundary component $\Sigma$ and $\mathring{\mathrm{II}}_g$ is not identically zero there. There are corrected bubbles $u_\varepsilon=\Phi(x_\varepsilon,\varepsilon)$ with centers in a compact non-umbilic subset of $\Sigma$ such that
\[
\mathcal J_g[u_\varepsilon]<S_\ast,\qquad
r_g(u_\varepsilon)\asymp\varepsilon^2,\qquad
\dist_E(u_\varepsilon,\mathcal B_0)\asymp\varepsilon.
\]
Consequently the exponent $1/2$ in \eqref{II-eq:residual-standard-sqrt} cannot be increased. The exponent $1$ in \eqref{II-eq:residual-above-linear} is also optimal.
\end{proposition}

\begin{proof}
For each small $\varepsilon$, choose $x_\varepsilon$ minimizing $\mathcal J_1(x,\varepsilon)$ on the compact manifold $\Sigma$. The uniform expansion of $\mathcal J_1$ and $\kappa_3^{\mathrm{red}}<0$ imply that every limit point of $x_\varepsilon$ maximizes $|\mathring{\mathrm{II}}|^2$ on $\Sigma$. Hence $|\mathring{\mathrm{II}}(x_\varepsilon)|$ is bounded below, and $\mathcal J_1(x_\varepsilon,\varepsilon)<S_\ast$. At these centers $\nabla_\partial\mathcal J_1=0$. The scale expansion and \eqref{II-eq:residual-parameter-equivalence} therefore give $r_g(u_\varepsilon)\asymp\varepsilon^2$.

The source estimate and injectivity of the source map in Lemma~\ref{I-lem:spin2-profile}, together with compactness of the unit sphere in $\mathrm{Sym}^2_0$, give $\|D_XF_p(0)\|_{X_p^*}\ge c\varepsilon$ at these centers. Since $D_XF_p(w_p)=0$ and the Hessian is uniformly bounded, $\|w_p\|_E\ge c\varepsilon$. The reverse inequality is Lemma~\ref{I-lem:LS}. Uniform normalized modulation gives
\[
\dist_E(\Phi(p),\mathcal B_0)\asymp\|w_p\|_E\asymp\varepsilon.
\]
Indeed, sufficiently close bubbles belong to the same normalized parameter chart, where Lemma~\ref{B-lem:normalized-modulation} gives the lower estimate; bubbles outside that chart stay a fixed distance away. This proves sharpness of $1/2$.

To prove sharpness above $S_\ast$, choose $h_\varepsilon\in X_{x_\varepsilon,\varepsilon}$ with $\|h_\varepsilon\|_E=1$, fix a sufficiently large constant $L$, and put
\[
\widetilde u_\varepsilon
=\mathsf{ret}_{\mathcal U_{x_\varepsilon,\varepsilon}}
(w_{x_\varepsilon,\varepsilon}+L\varepsilon h_\varepsilon).
\]
The lower bound in \eqref{II-eq:residual-normal-energy} makes $\mathcal J_g[\widetilde u_\varepsilon]\ge S_\ast$ for large fixed $L$. Taking $L$ also larger than twice a uniform upper bound for $\varepsilon^{-1}\|w_{x_\varepsilon,\varepsilon}\|_E$, modulation gives $\dist_E(\widetilde u_\varepsilon,\mathcal B_0)\asymp\varepsilon$. Complementary coercivity gives $r_g(\widetilde u_\varepsilon)\ge c\varepsilon$, and the Lipschitz estimate for the gradient gives the reverse bound. Thus no exponent larger than $1$ is possible on this side. Moreover, no sign conditions are required on these sharpness sequences.\end{proof}

\begin{remark}[A quantitative obstruction to concentration]
\label{II-rem:residual-obstruction}
Without a lower bound on $|\mathring{\mathrm{II}}|$, the same argument gives, uniformly on compact subsets of a minimal boundary region, $\varepsilon^2|\mathring{\mathrm{II}}(x)|^2\le Cr_g(u)+o(\varepsilon^2)$. In particular, if a modulated concentrating sequence has $r_g(u_j)=o(\varepsilon_j^2)$ and $x_j\to p$, then $p$ is umbilic.
\end{remark}

\section{Critical points and compactness}
\label{II-sec:threshold}
Throughout this section assume $n\ge5$, \textup{(Y$^+_{\partial}$)}, and \textup{(BG$^{3+}$)}, unless stronger hypotheses are stated.

\subsection{Critical points and the reduced functional}
\label{II-sec:balanced}

The quadratic scale coefficient excludes sufficiently small stationary scales away from the umbilic locus and vanishes on it. We use modulation and uniqueness to pass from criticality to stationarity of the reduced functional.

\begin{proposition}[One-bubble critical points and the reduced functional]
\label{II-prop:isolated-simple-relation}\label{I-lem:onebubble-to-reduced-criticality}
Let $(u_\ell)\subset\mathcal S$ be constrained critical points of $\mathcal J_g$. Suppose that, after passing to a subsequence and making one fixed global sign change if necessary, there are $x_\ell^0\in\partial M$ and $\varepsilon_\ell^0\downarrow0$ such that
\[
\|u_\ell-\mathcal U_{x_\ell^0,\varepsilon_\ell^0}\|_{H^1(M)}
\longrightarrow0.
\]
Then there are $x_\ell\in\partial M$ and $\varepsilon_\ell\downarrow0$ such that
\[
u_\ell=\Phi(x_\ell,\varepsilon_\ell),
\qquad
D\mathcal J_1(x_\ell,\varepsilon_\ell)=0.
\]

More generally, whenever a constrained critical point admits a representation $u=\mathsf{ret}_{\mathcal U_{x,\varepsilon}}(w)$ with a sufficiently small correction $w\in X_{x,\varepsilon}$, one has $w=w_{x,\varepsilon}$, $u=\Phi(x,\varepsilon)$, and $D\mathcal J_1(x,\varepsilon)=0$.
\end{proposition}

\begin{proof}
For the local statement, criticality makes $w$ solve \eqref{I-eq:LS-projected-grad}; uniqueness identifies it with $w_{x,\varepsilon}$, and the chain rule on $\mathcal S$ gives $D\mathcal J_1=0$. For the sequence, normalized modulation (Lemma~\ref{B-lem:normalized-modulation}) gives
\[
u_\ell=\mathsf{ret}_{\mathcal U_{x_\ell,\varepsilon_\ell}}(\widetilde w_\ell),
\qquad \widetilde w_\ell\in X_{x_\ell,\varepsilon_\ell},
\qquad \|\widetilde w_\ell\|_{H^1}\to0,
\qquad \varepsilon_\ell\downarrow0.
\]
For large $\ell$, criticality makes $\widetilde w_\ell$ solve \eqref{I-eq:LS-projected-grad}. Uniqueness gives $\widetilde w_\ell=w_{x_\ell,\varepsilon_\ell}$, hence $u_\ell=\Phi(x_\ell,\varepsilon_\ell)$. Since $\Phi$ takes values in $\mathcal S$, the chain rule gives $D\mathcal J_1(x_\ell,\varepsilon_\ell)=0$.
\end{proof}

\subsection{Blow-up points and compactness}

\begin{theorem}[Blow-up points at the level $S_\ast$]
\label{II-thm:threshold-Sstar-onebubble}
Assume $n\ge5$, \textup{(Y$^+_{\partial}$)}, and
\textup{(BG$^{3+}$)}. Let $(u_\ell)\subset\mathcal S$ be
constrained critical points such that
\[
\mathcal J_g[u_\ell]\longrightarrow S_\ast.
\]
If $(u_\ell)$ is not precompact in $H^1(M)$, then, after passing to a subsequence and making one fixed global sign change,
\[
u_\ell=\Phi(x_\ell,\varepsilon_\ell),
\qquad
\varepsilon_\ell\downarrow0,
\qquad
x_\ell\longrightarrow p\in\partial M.
\]
In a conformal representative with $H_g=0$,
\[
\mathfrak R_g^{\mathrm{red}}(p)=0,
\qquad
\mathring{\mathrm{II}}_g(p)=0,
\qquad
\nabla_\partial\mathfrak R_g^{\mathrm{red}}(p)=0.
\]
In particular, the limiting concentration point is umbilic in every conformal representative.
\end{theorem}

\begin{proof}
Lemma~\ref{I-lem:extremal-profile-onesigned} and normalized modulation give the concentrating one-bubble representation, after one global sign change. Proposition~\ref{II-prop:isolated-simple-relation} gives $u_\ell=\Phi(x_\ell,\varepsilon_\ell)$ and $D\mathcal J_1(x_\ell,\varepsilon_\ell)=0$. Extract $x_\ell\to p$.

The conformal change with boundary factor equal to one preserves the limiting center. Apply modulation again to the sequence of constrained critical points after this conformal change in the representative with $H_g=0$. The scale derivative gives
\[
0=S_\ast^{-1}\partial_\varepsilon\mathcal J_1(x_\ell,\varepsilon_\ell)
=2\mathfrak R_g^{\mathrm{red}}(x_\ell)\varepsilon_\ell+o(\varepsilon_\ell).
\]
Division by $\varepsilon_\ell$ gives $\mathfrak R_g^{\mathrm{red}}(p)=0$. Since $\mathfrak R_g^{\mathrm{red}}=\kappa_3^{\mathrm{red}}|\mathring{\mathrm{II}}|^2$ with $\kappa_3^{\mathrm{red}}<0$, both $\mathring{\mathrm{II}}(p)$ and $\nabla_\partial\mathfrak R_g^{\mathrm{red}}(p)$ vanish. The umbilic locus is conformally invariant.
\end{proof}

\begin{lemma}[Palais--Smale sequences at level $S_\ast$ with prescribed centers]
\label{II-lem:PS-obstruction}
Under the standing hypotheses, use the normalized bubble $\mathcal U_{x,\varepsilon}$ with the cutoff $R(\varepsilon)=\varepsilon^{-3/4}$. Uniformly for $x\in\partial M$,
\begin{equation}\label{II-eq:PS-obstruction-grad}
\|\mathsf{grad}_{\mathcal S}
\mathcal J_g(\mathcal U_{x,\varepsilon})\|_{H^1(M)}
\le
C\bigl(\varepsilon+\omega_\nabla(\varepsilon)\bigr)
=
O(\varepsilon),
\end{equation}
and
\[
\mathcal J_g(\mathcal U_{x,\varepsilon})
=
S_\ast+o(1).
\]
Consequently, every $p\in\partial M$ is the limiting center of a nonnegative Palais--Smale sequence at level $S_\ast$.
\end{lemma}

\begin{proof}
Lemma~\ref{I-lem:first-order-drift} gives the gradient estimate, and the standard bubble expansion gives $\mathcal J_g[\mathcal U_{x,\varepsilon}]=S_\ast+o(1)$. Taking $x=p$ and $\varepsilon\downarrow0$ proves the assertion at any prescribed point.
\end{proof}

\begin{remark}[Minimizing sequences of critical points and small scales]\label{II-thm:threshold-leading}\label{II-cor:small-scale-exclusion-onebubble}
A noncompact minimizing sequence of constrained critical points has level $S_\ast$ by Theorem~\ref{I-thm:compactness-escobar}, so Theorem~\ref{II-thm:threshold-Sstar-onebubble} applies. On a compact minimal boundary region where $|\mathfrak R_g^{\mathrm{red}}|$ is bounded below, the uniform scale expansion gives $\partial_\varepsilon\mathcal J_1<0$ at sufficiently small scales, excluding reduced critical points there.
\end{remark}

\begin{corollary}[Compactness below the first mixed level]
\label{II-cor:threshold-compactness}
Under the standing hypotheses, suppose that $\partial M$ is nowhere umbilic. For every $c_0<c_{\mathrm{mix}}$, with $c_{\mathrm{mix}}$ as in \eqref{I-eq:first-mixed-level}, the set of normalized constrained critical points satisfying $\mathcal J_g[u]\le c_0$ is compact in $H^1(M)$. In particular, every sequence of normalized constrained critical points with $\mathcal J_g[u_\ell]\to S_\ast$ is precompact.
\end{corollary}

\begin{proof}
Pass to a subsequence whose values converge. Below $S_\ast$ apply Theorem~\ref{I-thm:compactness-escobar}; at $S_\ast$, Theorem~\ref{II-thm:threshold-Sstar-onebubble} excludes blow-up because there are no umbilic points; above $S_\ast$ apply Theorem~\ref{I-thm:signfree-first-gap-compactness}. Strong convergence preserves the constraint, criticality, and the energy bound.
\end{proof}

Write $\operatorname{Met}^m(M)$ for the Banach manifold of $C^m$ Riemannian metrics. Nowhere-umbilic genericity is also recorded in \cite[p.~1765]{KimMussoWei2021}; the finite-regularity proof is given in Appendix~\ref{app:genericity-proof}.

\begin{theorem}[Generic strict subcriticality and threshold compactness]
\label{II-thm:generic-compactness}
Assume $n\ge5$ and fix an integer $m\ge3$. Let
\[
\operatorname{Met}_{Y,BG}^m(M)
\coloneqq 
\left\{
g\in\operatorname{Met}^m(M):
g\text{ satisfies \textup{(Y$^+_{\partial}$)}
and \textup{(BG$^{3+}$)}}
\right\}.
\]
The bounds and chart radii in \textup{(BG$^{3+}$)} may depend on the metric. Every $C^m$ metric on the fixed compact manifold satisfies these bounds, and strict coercivity is open under $C^2$ perturbations. Thus this is a $C^m$-open subset of $\operatorname{Met}^m(M)$. There is a $C^m$-open dense set
\[
\mathcal G\subset\operatorname{Met}_{Y,BG}^m(M)
\]
such that every $g\in\mathcal G$ satisfies
\begin{equation}\label{II-eq:generic-geometric-properties}
\mathring{\mathrm{II}}_g(x)\ne0
\quad(x\in\partial M),
\qquad
W_g\not\equiv0
\quad\text{in }M^\circ.
\end{equation}
In particular, $(M,g)$ is not conformally diffeomorphic to the round hemisphere. Moreover,
\begin{equation}\label{II-eq:generic-subcritical-attainment}
C^*_{\Esc}(M,g)<S_\ast,
\qquad
C^*_{\Esc}(M,g)=\mathcal J_g[u_g]
\quad\text{for some }u_g\in\mathcal S.
\end{equation}
Every normalized Palais--Smale sequence at $C^*_{\Esc}(M,g)$ is precompact in $H^1(M)$. For every $c_0<c_{\mathrm{mix}}$, the normalized constrained critical points with values at most $c_0$ form a compact subset of $H^1(M)$.

\end{theorem}
\begin{proof}
Appendix~\ref{app:genericity-proof} proves the boundary transversality, the independent interior Weyl perturbation, and the analytic conclusions at finite $C^m$ regularity.
\end{proof}

\section{The fourth-order reduced coefficient}\label{B-sec:quartic}

The quadratic scale equation forces blow-up points of sequences of constrained critical points at $S_\ast$ to be umbilic. We now continue the same reduction to fourth order at umbilic points. The expansion is constructed on the umbilic locus; its explicit curvature formula requires total umbilicity in a neighborhood. Total umbilicity also permits the expansion to be applied uniformly along varying centers of sequences of constrained critical points.

\begin{lemma}[Vanishing at umbilic centers]\label{I-lem:umbilic-lower-coefficients}
Assume $n\ge6$, \textup{(Y$^+_{\partial}$)}, and \textup{(BG$^{4+}$)}, and work in a boundary region with $H_g=0$. Write $\Theta_g$ for the cubic coefficient of the normalized standard-bubble energy $\mathcal J_g/S_\ast$, as in \cite[Lemma~4.11]{EscobarGeometry}. At every umbilic point $p$ in this region,
\[
\mathfrak R_g^{\mathrm{bare}}(p)=\mathfrak R_g^{\mathrm{red}}(p)
=\nabla_\partial\mathfrak R_g^{\mathrm{red}}(p)=\Theta_g(p)=0.
\]
\end{lemma}

\begin{proof}
The standard-bubble coefficients vanish by Proposition~\ref{I-prop:standard-bubble-inputs}. The formula $\mathfrak R_g^{\mathrm{red}}=\kappa_3^{\mathrm{red}}|\mathring{\mathrm{II}}|^2$ from Proposition~\ref{I-rem:reduced-vs-bare} gives the remaining equalities.
\end{proof}

\subsection{The second-order term and the linearized equation}

For $n\ge6$ and \textup{(BG$^{4+}$)}, the quadratic and geometric cubic coefficients vanish at an umbilic point with zero mean curvature by Lemma~\ref{I-lem:umbilic-lower-coefficients}. The order-$\varepsilon$ geometric source vanishes there as well. For $n\ge7$, with the cutoff below, the complementary source is $O(\varepsilon^2)$, and the leading correction obtained by completing the square has order $\varepsilon^4$.

For the rest of this section assume $n\ge7$ and fix $4/(n-2)<\beta<1$. We refine the cutoff to separate this contribution from the tails: $R_\beta(\varepsilon)=\varepsilon^{-\beta}$ makes the cutoff radius on $M$ tend to zero, the flat energy tail $o(\varepsilon^4)$, and the tail in the first variation $o(\varepsilon^2)$.

\begin{remark}[Notation for the $\beta$-cutoff family]
\label{II-rem:beta-fixedspace-data}
Use $\mathcal U^{(\beta)}_{p,\varepsilon}$ from Definition~\ref{B-recall-I-def:bubbles}, and let $X^{(\beta)}_{p,\varepsilon}$ and $P^{(\beta)}_{p,\varepsilon}$ be its complement and energy projection as in Lemma~\ref{II-lem:residual-kernel-pure}. Define $\mathscr S^{(\beta)}$ and $\mathscr L^{(\beta)}$ by \eqref{II-eq:onebubble-physical-data}, using this normalized bubble. The model comparison map is
\[
\mathfrak B^{(\beta)}_{p,\varepsilon}z
\coloneqq P^{(\beta)}_{p,\varepsilon}
\mathcal T_{p,\varepsilon}(\chi_{R_\beta(\varepsilon)}z).
\]
As in Remark~\ref{II-rem:blowup-trivialization-onebubble}, weighted derivatives of objects on $M$ are taken after identifying nearby complements by the inverse of the restricted energy projection; the superscript $\mathrm{ph}$ denotes differentiation at a fixed point of $M$.

Put $T_w^{(\beta)}=D\mathsf{ret}_{\mathcal U^{(\beta)}_{p,\varepsilon}}(w)|_{X^{(\beta)}_{p,\varepsilon}}$. The unique small correction satisfies
\begin{equation}\label{II-eq:beta-projected-equation}
\bigl(T_{w^{(\beta)}_{p,\varepsilon}}^{(\beta)}\bigr)^*
\mathsf{grad}_{\mathcal S}\mathcal J_g
\left(\mathsf{ret}_{\mathcal U^{(\beta)}_{p,\varepsilon}}
(w^{(\beta)}_{p,\varepsilon})\right)=0
\qquad\text{in }X^{(\beta)}_{p,\varepsilon}.
\end{equation}
The two cutoff conventions satisfy, uniformly in $p$,
\begin{equation}\label{II-eq:beta-standard-cutoff-closeness}
\|\mathcal U^{(\beta)}_{p,\varepsilon}-\mathcal U_{p,\varepsilon}\|_{H^1(M)}=o(1).
\end{equation}
\end{remark}

\begin{lemma}[Second-order source and correction]
\label{II-lem:second-source-response-expansion}
Assume $n\ge7$, \textup{(Y$^+_{\partial}$)}, and \textup{(BG$^{5+}$)}, and let $\mathcal C\subset\partial M$ be an open boundary region admitting a fixed Fermi collar, with $H_g=0$ on $\mathcal C$. Use the $\beta$-cutoff family and objects of Remark~\ref{II-rem:beta-fixedspace-data}. In each smoothly varying local Fermi frame there is a continuous family $\mathscr S_{2,p}\in X_+^*$ on the umbilic locus $\{p\in\mathcal C:\mathring{\mathrm{II}}_g(p)=0\}$ such that
\begin{equation}\label{II-eq:second-source-expansion}
\left\|\varepsilon^{-2}
(\mathfrak B^{(\beta)}_{p,\varepsilon})^*
\mathscr S^{(\beta)}_{p,\varepsilon}-\mathscr S_{2,p}\right\|_{X_+^*}
=o(1).
\end{equation}
Let $\mathcal A=(n-2)/2+Y\cdot\nabla$ and
$\mathcal D(\mathcal A)=\{z\in\mathcal D^{1,2}:\mathcal Az\in\mathcal D^{1,2}\}$.
For each fixed $z\in X_+\cap\mathcal D(\mathcal A)$,
\begin{equation}\label{II-eq:second-source-fixed-derivative}
\mathbf D\left[\varepsilon^{-2}
\langle\mathscr S^{(\beta)}_{p,\varepsilon},
\mathfrak B^{(\beta)}_{p,\varepsilon}z\rangle\right]=o(1).
\end{equation}
This last convergence is uniform for $z$ in sets compact with respect to the graph norm
$\|z\|_{\mathcal D^{1,2}}+\|\mathcal Az\|_{\mathcal D^{1,2}}$.
Both statements are uniform for $p$ in compact subsets of the umbilic locus. A change of tangential orthonormal frame acts by the corresponding orthogonal change of variables on $X_+$ and its dual.

Put $z_{2,p}=\mathscr L_0^{-1}\mathscr S_{2,p}$ and $v^{(\beta)}_{p,\varepsilon}=(\mathscr L^{(\beta)}_{p,\varepsilon})^{-1}\mathscr S^{(\beta)}_{p,\varepsilon}$. Then
\begin{equation}\label{II-eq:second-response-comparison}
\|v^{(\beta)}_{p,\varepsilon}
-\varepsilon^2\mathfrak B^{(\beta)}_{p,\varepsilon}z_{2,p}\|_{H^1(M)}
=o(\varepsilon^2),
\qquad
\|\mathbf D^{\mathrm{ph}}v^{(\beta)}_{p,\varepsilon}\|_{H^1(M)}
\le C\varepsilon^2.
\end{equation}
Moreover,
\begin{equation}\label{II-eq:second-pairing-C1}
\begin{gathered}
\langle\mathscr S^{(\beta)},(\mathscr L^{(\beta)})^{-1}\mathscr S^{(\beta)}\rangle
=\varepsilon^4\langle\mathscr S_{2,p},\mathscr L_0^{-1}\mathscr S_{2,p}\rangle+e_2,\\
|e_2|+|\mathbf De_2|=o(\varepsilon^4).
\end{gathered}
\end{equation}
The estimates are uniform when $\beta$ ranges over a compact subinterval of $(4/(n-2),1)$, with the same fixed cutoff profile. In particular, $\mathscr S_{2,p}$ and $z_{2,p}$ are independent of $\beta$.
\end{lemma}

\begin{proof}
Put $\alpha=(n-2)/2$ and $c_n=(n-2)/(4(n-1))$. All curvature tensors in the following formulas are evaluated at $p$; barred tensors belong to the induced boundary metric, and a semicolon denotes covariant differentiation.

\smallskip\noindent \emph{The source.}
At an umbilic point with zero mean curvature, $\mathrm{II}(p)=0$. The second-order Fermi coefficients are
\begin{equation}\label{II-eq:second-Fermi-coefficients}
\begin{aligned}
b_2(Y')&=-\tfrac16\overline{\Ric}_{cd}Y_cY_d,&
d_2(Y)&=b_2(Y')-\tfrac12\Ric_{nn}t^2,\\
A_2^{ab}(Y)&=\tfrac13\overline R_{acbd}Y_cY_d
+2\mathrm{II}_{ab;c}Y_ct+R_{anbn}t^2+d_2(Y)\delta_{ab},\\
A_2^{nn}(Y)&=d_2(Y),& A_2^{an}(Y)&=0.
\end{aligned}
\end{equation}
Thus $\sqrt{\det g(\varepsilon Y)}g^{ij}(\varepsilon Y)=\delta^{ij}+\varepsilon^2A_2^{ij}(Y)+O(\varepsilon^3|Y|^3)$ and the boundary density is $1+\varepsilon^2b_2(Y')+O(\varepsilon^3|Y'|^3)$. On $X_+$ the source is
\begin{equation}\label{II-eq:second-source-explicit}
\mathscr S_{2,p}[\phi]=2t_\ast\bigg(\int_{\R^n_+}A_2^{ij}\partial_iU_+\partial_j\phi
+c_n\Scal(p)\int_{\R^n_+}U_+\phi-\lambda_{\mathrm{flat}}\int_{\partial\R^n_+}b_2U_+^{q-1}\phi\bigg).
\end{equation}
Trace normalization adds only a multiple of the constraint covector, which vanishes on $X_+$. The terms $\mathrm{II}_{ab;c}$ remain at a general umbilic point. They vanish when the boundary is totally umbilic near $p$ in this representative with minimal boundary.

Since $A_2^{nn}(Y',0)=b_2(Y')$, $A_2^{an}=0$, and $-\partial_tU_+=\lambda_{\mathrm{flat}}U_+^{q-1}$ on the boundary, integration by parts cancels the last term in \eqref{II-eq:second-source-explicit}. Its extension to the full Dirichlet space therefore has interior forcing
\[
f_{2,p}=2t_\ast\bigl[-\operatorname{div}(A_2\nabla U_+)+c_n\Scal(p)U_+\bigr],
\qquad |f_{2,p}|+|Y\cdot\nabla f_{2,p}|\le C(1+|Y|)^{2-n}.
\]
The model moments converge because
\begin{equation}\label{II-eq:second-source-coefficient-tail}
\bigl\||Y|^2\nabla U_+\bigr\|_{L^2(|Y|>R)}
+\bigl\||Y|U_+\bigr\|_{L^2(|Y|>R)}
\le CR^{-(n-6)/2}.
\end{equation}
The tail of the second-order boundary term has the stronger bound
$\||Y'|^2U_+^{q-1}\|_{L^{q'}(|Y'|>R)}\le CR^{-(n-4)/2}$, where $q'=q/(q-1)$. The next term in the interior Taylor expansion contributes $C\varepsilon^3\ell_n(R)$, where
\[
\ell_n(R)=
\begin{cases}R^{1/2},&n=7,\\\sqrt{\log(2R)},&n=8,\\1,&n\ge9.
\end{cases}
\]
Indeed, the required norms are $\||Y|^3\nabla U_+\|_{L^2(B_{2R}^+)}$ and $\||Y|^2U_+\|_{L^2(B_{2R}^+)}$. The third-order boundary term is $O(\varepsilon^3)$ since $|Y'|^3U_+^{q-1}\in L^{q'}$. If $\mathscr C_R$ denotes the flat cutoff covector and $\mathscr E_{p,\varepsilon,R}$ the Taylor error, then
\begin{equation}\label{II-eq:second-source-Taylor-tail}
\|\mathscr C_R\|+\|R\partial_R\mathscr C_R\|\le CR^{-\alpha},
\qquad
\|\mathscr E_{p,\varepsilon,R}\|+
\|(\varepsilon\partial_\varepsilon-\beta R\partial_R)\mathscr E_{p,\varepsilon,R}\|
\le C\varepsilon^3\ell_n(R).
\end{equation}
Consequently the normalized error, also after one logarithmic scale derivative in model coordinates, is bounded by
\[
C\bigl(\varepsilon^{-2}R^{-\alpha}+R^{-(n-6)/2}
+\varepsilon\ell_n(R)\bigr)=o(1),
\qquad R=\varepsilon^{-\beta}.
\]
Here $\beta\alpha>2$ and $\beta<1$. The inverse Gram matrices are uniformly bounded; the frames converge strongly, and the constraint normal is controlled by Remark~\ref{B-rem:inverse-Gram}. The $O(\varepsilon)$ bound for the full first variation and the Riesz identity give the strong frame convergence. Projecting the source therefore preserves this limit and proves \eqref{II-eq:second-source-expansion}. Continuity in $p$ follows from the displayed finite collection of curvature jets.

\smallskip\noindent \emph{Logarithmic scale derivatives at fixed points of $M$.}
Apply \eqref{II-eq:physical-source-differentiation} with $k=2$. Hardy's inequality gives
\[
\int_{|Y|>1}|Y|^2\bigl(|f_{2,p}|^2+|Y\cdot\nabla f_{2,p}|^2\bigr)\,dY
\le C\int_1^\infty r^{5-n}\,dr<\infty.
\]
The explicit cutoff and Taylor expressions obey the differentiated bounds above. The constraint contribution from normalization is estimated by its scalar coefficient and the bounded constraint covector. Thus the full first variation and bilinear form on $M$ satisfy
\begin{equation}\label{II-eq:second-global-source-operator-bound}
\|\mathscr S^{\mathrm{full}}\|+\|\mathbf D^{\mathrm{ph}}\mathscr S^{\mathrm{full}}\|
\le C\varepsilon^2,
\qquad \|a_\varepsilon\|+\|\mathbf D^{\mathrm{ph}}a_\varepsilon\|\le C.
\end{equation}
As in the proof of Proposition~\ref{II-prop:onebubble-scale-derivative-II}, the Riesz identity and inverse Gram formula give bounded logarithmic scale derivatives of the projection at fixed points of $M$. Stabilized inversion yields
\begin{equation}\label{II-eq:second-global-response-bound}
\|v\|_{H^1}+\|\mathbf D^{\mathrm{ph}}v\|_{H^1}\le C\varepsilon^2.
\end{equation}

\smallskip\noindent \emph{The solution of the linearized equation.}
The fixed-ball and tail argument of Lemma~\ref{I-lem:resolvent-limit}, using \eqref{II-eq:second-source-coefficient-tail}, identifies every rescaled weak limit of $\varepsilon^{-2}v$ with $z_{2,p}=\mathscr L_0^{-1}\mathscr S_{2,p}$. In particular,
\[
\varepsilon^{-4}\langle\mathscr S,v\rangle
\longrightarrow C_p\coloneqq \langle\mathscr S_{2,p},z_{2,p}\rangle.
\]
For the dilation regularity of $z_{2,p}$, extend the complementary equation to the full Dirichlet space by adding finitely many constraint and Jacobi covectors. Their Dirichlet representatives are harmonic and add only Robin boundary forcing of size $O((1+|Y'|)^{-n})$, with the corresponding differentiated decay. These terms have finite boundary norms in \eqref{II-eq:weighted-annular-H2}. The displayed Hardy integral and that estimate give $Y\nabla^2z_{2,p}\in L^2$; local regularity, Hardy's inequality, and cutoff approximation then give $\mathcal Az_{2,p}\in\mathcal D^{1,2}$.

For a fixed $z\in X_+\cap\mathcal D(\mathcal A)$, the transported test and its logarithmic scale derivative at a fixed point of $M$ are bounded by its graph norm, by \eqref{II-eq:fixed-model-rescaling-derivative}. The coefficient estimates give \eqref{II-eq:second-source-fixed-derivative} before projection. The removed moments and their derivatives tend to zero: for each limiting normal field $e$,
\[
\langle\mathcal Az,e\rangle_{\mathcal D^{1,2}}
+\langle z,\mathcal Ae\rangle_{\mathcal D^{1,2}}=0
\]
by the energy isometry of critical dilations. The exact Riesz identity handles the constraint normal on $M$ with an $O(\varepsilon^2)$ error, also after differentiation. The inverse Gram formula proves the same statement after projection. Uniform graph-norm bounds and a finite-net argument give uniformity on compact subsets of the graph-norm space.

\smallskip\noindent \emph{Completing the square.}
Set $V_\varepsilon=\varepsilon^2\mathfrak B^{(\beta)}_{p,\varepsilon}z_{2,p}$. The preceding source comparison and the calculation of the bilinear-form pairing in the proof of Proposition~\ref{II-prop:onebubble-scale-derivative-II} give
\[
\mathscr S[V_\varepsilon]=\varepsilon^4C_p+o(\varepsilon^4),
\qquad a_\varepsilon(V_\varepsilon,V_\varepsilon)=\varepsilon^4C_p+o(\varepsilon^4),
\]
including $o(\varepsilon^4)$ weighted derivatives of both errors, and
$\|\mathbf D^{\mathrm{ph}}V_\varepsilon\|_{H^1}\le C\varepsilon^2$. Equation~\eqref{II-eq:quadratic-completion-identity} and coercivity imply $\|v-V_\varepsilon\|_{H^1}=o(\varepsilon^2)$. Both logarithmic scale derivatives at fixed points of $M$ are $O(\varepsilon^2)$, so Lemma~\ref{II-lem:differentiated-quadratic-completion} proves \eqref{II-eq:second-pairing-C1}. The finite-dimensional dependence on the curvature jets gives uniformity in $p$. All exponents remain positive on compact subintervals of $4/(n-2)<\beta<1$, and the limiting equation has no $\beta$ dependence.
\end{proof}

\subsection{The fourth-order expansion}

\begin{lemma}[Remainder estimate and critical points for the cutoff family]
\label{II-lem:beta-cutoff-relation}
Under the hypotheses of Lemma~\ref{II-lem:second-source-response-expansion}, let $w^{(\beta)}_{x,\varepsilon}$ be the unique small one-bubble Lyapunov--Schmidt correction and define
\[
\Phi_\beta(x,\varepsilon)
\coloneqq 
\mathsf{ret}_{\mathcal U^{(\beta)}_{x,\varepsilon}}
\bigl(w^{(\beta)}_{x,\varepsilon}\bigr),
\qquad
\mathcal J_{1,\beta}
\coloneqq 
\mathcal J_g\circ\Phi_\beta.
\]
Put
\[
\mathscr E_\beta(x,\varepsilon)
\coloneqq 
\mathcal J_{1,\beta}(x,\varepsilon)
-\mathcal J_g[\mathcal U^{(\beta)}_{x,\varepsilon}]
+\frac12
\left\langle
(\mathscr L^{(\beta)}_{x,\varepsilon})^{-1}
\mathscr S^{(\beta)}_{x,\varepsilon},
\mathscr S^{(\beta)}_{x,\varepsilon}
\right\rangle .
\]
If $\mathring{\mathrm{II}}_g(x)=0$, then
\[
\|\mathscr S^{(\beta)}_{x,\varepsilon}\|_
{(X^{(\beta)}_{x,\varepsilon})^*}
+
\|w^{(\beta)}_{x,\varepsilon}\|_{H^1(M)}
\le C\varepsilon^2
\]
and, with $\mathbf D\coloneqq \varepsilon\partial_\varepsilon$,
\begin{equation}\label{II-eq:beta-holder-tail}
|\mathscr E_\beta|
+
|\mathbf D\mathscr E_\beta|
=
o(\varepsilon^4).
\end{equation}
These estimates are uniform on compact subsets of the umbilic locus in $\mathcal C$.

Suppose that the boundary is totally umbilic in the representative satisfying $H_g=0$. If $(u_\ell)\subset\mathcal S$ is a non-precompact sequence of constrained critical points with
\[
\mathcal J_g[u_\ell]\longrightarrow S_\ast,
\]
then, after passing to a subsequence and making one fixed global sign change, there are
\[
x_\ell\longrightarrow p\in\partial M,
\qquad
\varepsilon_\ell\downarrow0,
\]
such that, for all sufficiently large $\ell$,
\begin{equation}\label{II-eq:beta-criticality-relation}
u_\ell=\Phi_\beta(x_\ell,\varepsilon_\ell),
\qquad
D\mathcal J_{1,\beta}(x_\ell,\varepsilon_\ell)=0.
\end{equation}
\end{lemma}

\begin{proof}
The bubble truncation error in rescaled coordinates has Dirichlet norm $O(R_\beta^{-\alpha})=o(\varepsilon^2)$, $\alpha=(n-2)/2$. The same bound holds after logarithmic scale differentiation and for the normalized Jacobi fields. The frame and gap estimates of Theorem~\ref{I-thm:intrinsic-single-gap}, Lemma~\ref{B-lem:normalized-modulation}, and Remark~\ref{B-rem:inverse-Gram} therefore apply uniformly. Lemma~\ref{II-lem:second-source-response-expansion}, followed by the inversion in Lemmas~\ref{I-lem:LS-C1-inversion} and~\ref{I-lem:LS}, gives $\|\mathscr S^{(\beta)}\|+\|w^{(\beta)}\|_{H^1}\le C\varepsilon^2$.

Suppress $x,\varepsilon,\beta$ and put $v=\mathscr L^{-1}\mathscr S$. At a fixed point of $M$,
\[
\mathbf D^{\mathrm{ph}}\mathcal U=b(Y)\mathcal U+\mathcal A^{\mathrm{ann}},
\qquad \|b\|_\infty\le C,\qquad
\|\mathcal A^{\mathrm{ann}}\|_E\le CR_\beta^{-\alpha}.
\]
Where the bubble cutoff equals one, $b=\mathbf D\log t-\alpha-Y\cdot\nabla\log U_+$. The estimates for $N$ and $Q$ in Proposition~\ref{II-prop:onebubble-scale-derivative-II} control this bounded term. On the annulus, \eqref{II-eq:annular-nonlinear-cancellation} gives a boundary-potential coefficient with $L^{n-1}$ norm $O(R_\beta^{-1})$. The weighted finite-rank coordinate derivatives are bounded. Thus the nonlinear gradient remainder satisfies
\begin{align*}
\|\mathcal N_\beta(w)\|&\le C\|w\|^{1+\gamma},&
\|D_w\mathcal N_\beta(w)\|&\le C\|w\|^\gamma,\\
\|\mathbf D^{\mathrm{par}}\mathcal N_\beta(w)\|
&\le C\bigl(\|w\|^{1+\gamma}
+R_\beta^{-\alpha}\|w\|^\gamma+R_\beta^{-1}\|w\|\bigr).
\end{align*}
Here $\mathbf D^{\mathrm{par}}$ differentiates the parameter at fixed $w$ in the complement coordinates. For $\|w\|=O(\varepsilon^2)$ the three orders are
\[
\varepsilon^{2+2\gamma},\qquad
\varepsilon^{\beta\alpha+2\gamma},\qquad
\varepsilon^{\beta+2}.
\]
Since $\beta\alpha>2$ and $\beta>2\gamma$, all are $O(\varepsilon^{2+2\gamma})$. Using \eqref{II-eq:second-global-source-operator-bound}--\eqref{II-eq:second-global-response-bound} and differentiating the projected equation gives, with $\mathbf D$ denoting the total logarithmic scale derivative in the complement coordinates,
\[
\|w\|+\|\mathbf Dw\|+\|v\|+\|\mathbf Dv\|\le C\varepsilon^2.
\]
The equation $\mathscr L(w+v)=-\mathcal N_\beta(w)$ and its derivative then give
\[
\|w+v\|+\|\mathbf D(w+v)\|\le C\varepsilon^{2+2\gamma}.
\]

For the energy remainder after the linear and quadratic terms,
\begin{align*}
|\mathcal Q_\beta(w)|&\le C\|w\|^{2+\gamma},\\
|\mathbf D\mathcal Q_\beta(w)|
&\le C\bigl(\|w\|^{2+\gamma}+R_\beta^{-\alpha}\|w\|^{1+\gamma}
+R_\beta^{-1}\|w\|^2+\|\mathcal N_\beta(w)\|\|\mathbf Dw\|\bigr)
=o(\varepsilon^4).
\end{align*}
The energy terms have orders $\varepsilon^{4+2\gamma}$, $\varepsilon^{\beta\alpha+2+2\gamma}$, and $\varepsilon^{\beta+4}$, all $o(\varepsilon^4)$. Lemma~\ref{II-lem:holder-quadratic-expansion} and completing the square give
\[
\mathscr E_\beta=\frac12\langle\mathscr L(w+v),w+v\rangle+\mathcal Q_\beta(w).
\]
The preceding bounds and the bounded weighted scale derivative of $\mathscr L$ prove \eqref{II-eq:beta-holder-tail}.

For the sequence of constrained critical points, Lemma~\ref{I-lem:extremal-profile-onesigned} gives a single signed bubble. The $o(1)$ cutoff comparison \eqref{II-eq:beta-standard-cutoff-closeness} and Lemma~\ref{B-lem:normalized-modulation} give $u_\ell=\mathsf{ret}_{\mathcal U^{(\beta)}_{x_\ell,\varepsilon_\ell}}(\widetilde w_\ell)$ with $\widetilde w_\ell\to0$ in the complement. Criticality makes $\widetilde w_\ell$ solve \eqref{II-eq:beta-projected-equation}. Uniqueness in Lemma~\ref{I-lem:LS} identifies it with $w^{(\beta)}_{x_\ell,\varepsilon_\ell}$, and the chain rule gives \eqref{II-eq:beta-criticality-relation}.
\end{proof}

\begin{proposition}[Fourth-order one-bubble expansion at umbilic points]
\label{II-prop:quartic-reduced-umbilic}
Under the hypotheses and cutoff convention of Lemma~\ref{II-lem:second-source-response-expansion}, let $\mathcal J_{1,\beta}$ be the reduced functional of Lemma~\ref{II-lem:beta-cutoff-relation}.

For $p\in\mathcal C$ with $\mathring{\mathrm{II}}_g(p)=0$, the limit
\[
\Xi_g^{0}(p)
\coloneqq 
\lim_{\varepsilon\downarrow0}
\varepsilon^{-4}
\left(
\frac{\mathcal J_g[
\mathcal U^{(\beta)}_{p,\varepsilon}]}{S_\ast}-1
\right)
\]
exists. With $\mathscr S_{2,p}$ and $\mathscr L_0$ as in Lemma~\ref{II-lem:second-source-response-expansion}, define
\begin{equation}\label{II-eq:Xi-reduced-definition}
\Xi_g^{\mathrm{red}}(p)
\coloneqq 
\Xi_g^{0}(p)
-\frac1{2S_\ast}
\left\langle
\mathscr L_0^{-1}\mathscr S_{2,p},
\mathscr S_{2,p}
\right\rangle .
\end{equation}
These coefficients are independent of $\beta$ in the stated range. Uniformly on compact subsets of the umbilic locus in $\mathcal C$,
\begin{align}
\frac{\mathcal J_{1,\beta}(p,\varepsilon)}{S_\ast}
&=
1+\Xi_g^{\mathrm{red}}(p)\varepsilon^4
+o(\varepsilon^4),
\label{II-eq:quartic-reduced-value}\\
\frac1{S_\ast}
\partial_\varepsilon
\mathcal J_{1,\beta}(p,\varepsilon)
&=
4\Xi_g^{\mathrm{red}}(p)\varepsilon^3
+o(\varepsilon^3).
\label{II-eq:quartic-reduced-scale}
\end{align}
In particular, the third-order coefficient of the reduced functional vanishes at every such point.
\end{proposition}

\begin{proof}
For $n\ge7$, the fourth gradient moment $\int|Y|^4|\nabla U_+|^2$ and fourth trace moment $\int|Y'|^4U_+^q$ are finite. The fourth-order Fermi expansion has fifth-order remainder $O(\varepsilon^5\log R)$ at $n=7$ and $O(\varepsilon^5)$ for $n\ge8$. For $m\le3$, its omitted gradient tail is
\[
O(\varepsilon^mR_\beta^{m+2-n})
=O(\varepsilon^{\beta(n-2)+m(1-\beta)})=o(\varepsilon^4).
\]
The flat and fourth-order tails are $O(R_\beta^{2-n})=o(\varepsilon^4)$ and $\varepsilon^4O(R_\beta^{6-n})=o(\varepsilon^4)$. Scalar and boundary terms have at least the same decay. These bounds persist after one weighted scale derivative.

The first three coefficients vanish at the umbilic point with zero mean curvature. Thus, with $\mathbf D=\varepsilon\partial_\varepsilon$,
\[
\frac{\mathcal J_g[\mathcal U^{(\beta)}_{p,\varepsilon}]}{S_\ast}
=1+\Xi_g^0(p)\varepsilon^4+e_0,
\qquad |e_0|+|\mathbf De_0|=o(\varepsilon^4).
\]
The differentiated pairing estimate \eqref{II-eq:second-pairing-C1} gives
\[
\langle(\mathscr L^{(\beta)})^{-1}\mathscr S^{(\beta)},\mathscr S^{(\beta)}\rangle
=\varepsilon^4\langle\mathscr L_0^{-1}\mathscr S_{2,p},\mathscr S_{2,p}\rangle+e_1,
\qquad |e_1|+|\mathbf De_1|=o(\varepsilon^4).
\]
Insert these identities into the quadratic expansion of Lemma~\ref{II-lem:beta-cutoff-relation}, whose remaining error has the same bound. This proves \eqref{II-eq:Xi-reduced-definition} and \eqref{II-eq:quartic-reduced-value}; applying $\varepsilon^{-1}\mathbf D$ proves \eqref{II-eq:quartic-reduced-scale}. The tail estimates and the limiting source equation contain no $\beta$, so both coefficients are independent of the admissible cutoff exponent.
\end{proof}

\subsection{Conformal invariance}

\begin{lemma}[Conformal invariance of the reduced fourth-order coefficient]
\label{II-lem:quartic-gauge-completion}
Assume $n\ge7$, \textup{(Y$^+_{\partial}$)}, \textup{(BG$^{5+}$)}, and that the boundary is totally umbilic near $p$. Let $g$ be the fixed representative with $H_g=0$. Choose a globally defined conformal metric $\widehat g=\varphi^{4/(n-2)}g$ with $H_{\widehat g}=0$, conformal-Fermi volume density $1+O(|x|^6)$ at $p$, and $\varphi(p)=1$, $d\varphi(p)=0$.
Construct the one-bubble reductions for $g$ and $\widehat g$ using the same admissible $\beta\in(4/(n-2),1)$ and the same cutoff profile. Then $\Xi_g^{\mathrm{red}}(p)=\Xi_{\widehat g}^{\mathrm{red}}(p)$.

More precisely, let $B,B'$ be the fourth-order coefficients in energy units of the two uncorrected families after the conformal change, let $\mathscr S,\mathscr S'\in X_+^*$ be their second-order source coefficients, and let $\mathscr L_0:X_+\to X_+^*$ be the flat complementary Hessian. If, after trace normalization and modulation, the two families differ at order $\varepsilon^2$ by $a\in X_+$, then
\begin{equation}\label{II-eq:quartic-gauge-completion}
\begin{gathered}
B'=B+\mathscr S[a]+\frac12\langle\mathscr L_0a,a\rangle,
\qquad
\mathscr S'=\mathscr S+\mathscr L_0a,\\
B'-\frac12\left\langle\mathscr L_0^{-1}\mathscr S',\mathscr S'\right\rangle
=
B-\frac12\left\langle\mathscr L_0^{-1}\mathscr S,\mathscr S\right\rangle.
\end{gathered}
\end{equation}
The reduced coefficient is invariant under this change of approximate family.
\end{lemma}

\begin{proof}
By Lemma~\ref{I-lem:graph-norm}, multiplication $\mathcal C_\varphi v=\varphi v$ is an isometry of the energy spaces and preserves the boundary constraint:
\begin{equation}\label{II-eq:quartic-conformal-slice-isometry}
\langle\mathcal C_\varphi v,\mathcal C_\varphi w\rangle_{E,g}
=\langle v,w\rangle_{E,\widehat g},\qquad
\|\mathcal C_\varphi v\|_{L^q(g)}=\|v\|_{L^q(\widehat g)},\qquad
\mathcal J_g[\mathcal C_\varphi v]=\mathcal J_{\widehat g}[v].
\end{equation}
It also commutes with trace normalization. We compare the two bubble families on the same constraint.

Align the tangential frames at $p$. Since $\varphi(p)=1$, $d\varphi(p)=0$, and both boundaries are totally geodesic, the Fermi transition has no quadratic term:
\begin{equation}\label{II-eq:quartic-Fermi-transition}
\widehat\Phi_p^{-1}\circ\Phi_p(x)=x+C_3(x)+O(|x|^4),
\qquad \varphi(\Phi_p(x))=1+q_2(x)+O(|x|^3),
\end{equation}
where $C_3$ and $q_2$ are homogeneous of degrees three and two. Before trace normalization, critical rescaling gives
\begin{equation}\label{II-eq:quartic-profile-gauge-difference}
\varepsilon^2a_p^0+o_{\mathcal D^{1,2}}(\varepsilon^2),
\qquad a_p^0=t_\ast(q_2U_++C_3\cdot\nabla U_+).
\end{equation}
Its gradient is $O((1+|Y|)^{3-n})$, and hence
\[
\|\nabla a_p^0\|_{L^2(|Y|>R)}\le CR^{-(n-6)/2}.
\]
The next Taylor terms have Dirichlet norm at most $C\varepsilon^3\ell_n(R)$, with $\ell_n$ as in Lemma~\ref{II-lem:second-source-response-expansion}. Higher products are absorbed by this bound on $\varepsilon R<1$. Trace normalization and projection onto the flat complement give
\[
a_p=P_+^\perp\left(a_p^0-U_\ast
\int_{\partial\R^n_+}U_\ast^{q-1}a_p^0\right)\in X_+.
\]
The uniform modulation estimate gives weighted parameter changes $O(\varepsilon^2)$, or changes in the center and scale $d_\partial(p_\varepsilon,p)+|\widehat\varepsilon-\varepsilon|=O(\varepsilon^3)$. Thus
\begin{equation}\label{II-eq:quartic-profile-gauge-comparison}
\mathcal C_\varphi\widehat{\mathcal U}^{(\beta)}_{p_\varepsilon,\widehat\varepsilon}
=\mathsf{ret}_{\mathcal U^{(\beta)}_{p,\varepsilon}}
\bigl(\varepsilon^2\mathfrak B^{(\beta)}_{p,\varepsilon}a_p+o_{H^1}(\varepsilon^2)\bigr).
\end{equation}

Let $\Delta\chi_{R_\beta}$ denote the difference between the two cutoff functions after the change of Fermi coordinates. On the cutoff annulus the rescaled coordinate displacement is $O(\varepsilon^2R_\beta^3)$, so its ratio to the annular radius is $O(\varepsilon^2R_\beta^2)=o(1)$. The mean-value theorem gives
\[
|\Delta\chi_{R_\beta}|\le C\varepsilon^2R_\beta^2,
\qquad |\nabla\Delta\chi_{R_\beta}|\le C\varepsilon^2R_\beta.
\]
Using $U_+=O(R_\beta^{2-n})$ and $\nabla U_+=O(R_\beta^{1-n})$ there,
\[
\|\nabla(\Delta\chi_{R_\beta}U_+)\|_2
\le C\varepsilon^2R_\beta^{-(n-6)/2}.
\]
Thus, if $\Delta_\varepsilon$ is the critically rescaled profile difference after normalization and modulation, the leading tail and the next Taylor terms give
\[
\|\varepsilon^{-2}\Delta_\varepsilon-a_p\|_{\mathcal D^{1,2}}
\le C\bigl(R_\beta^{-(n-6)/2}+\varepsilon\ell_n(R_\beta)\bigr)+o(1)=o(1).
\]
The displayed bounds account for the cutoff-transition error and the higher-order Taylor terms. The difference $e_\varepsilon=\varepsilon^{-2}\Delta_\varepsilon-\chi_{R_\beta}a_p$ tends to zero in the Dirichlet norm. Before the final projection its support lies in a half-ball of radius $C\varepsilon R_\beta$ on $M$, so Sobolev gives
\[
\|\mathcal T_{p,\varepsilon}e_\varepsilon\|_{L^2(M)}
\le C\varepsilon R_\beta\|e_\varepsilon\|_{\mathcal D^{1,2}}.
\]
The bounded energy projection and frame convergence therefore give the stated $o_{H^1}(\varepsilon^2)$ comparison.

The Taylor expansion on the common constraint is
\begin{equation}\label{II-eq:quartic-common-slice-expansion}
\mathcal J_g\!\left[\mathsf{ret}_{\mathcal U^{(\beta)}_{p,\varepsilon}}
(\varepsilon^2\mathfrak B^{(\beta)}_{p,\varepsilon}z)\right]
=S_\ast+\varepsilon^4\left(B+\mathscr S[z]+
\frac12\langle\mathscr L_0z,z\rangle\right)+o(\varepsilon^4)
\end{equation}
for each fixed $z\in X_+$, uniformly on compact subsets of $X_+$. Here the source converges in the model dual, the Hessian pairings converge on fixed model functions, and the $C^{2,\gamma}$ Taylor error is $O(\varepsilon^{4+2\gamma})$.

After conformal transport the second family has leading coordinates $a_p+z$ in this common parametrization. This follows first for compactly supported smooth model functions by the Fermi transition above, and then for fixed $z\in X_+$ by Dirichlet approximation and the uniform bounds for rescaling and projection. The limiting normal frames agree, so projection preserves this leading coordinate comparison. The full source estimate \eqref{II-eq:second-global-source-operator-bound} and the bounded Hessian give $\|D\mathcal J_g\|\le C\varepsilon^2$ throughout these $O(\varepsilon^2)$ neighborhoods. For a profile error $e=o_{H^1}(\varepsilon^2)$ the mean-value estimate is therefore
\[
|\mathcal J_g[u+e]-\mathcal J_g[u]|
\le\sup_{0\le s\le1}\|D\mathcal J_g[u+se]\|\|e\|_{H^1}
=o(\varepsilon^4).
\]

For every fixed $z\in X_+$ the common expansion gives
\[
B'+\mathscr S'[z]+\tfrac12\langle\mathscr L_0z,z\rangle
=B+\mathscr S[a_p+z]
+\tfrac12\langle\mathscr L_0(a_p+z),a_p+z\rangle.
\]
Comparison of the constant and linear terms gives the first two identities of \eqref{II-eq:quartic-gauge-completion}. Substituting $\mathscr S'=\mathscr S+\mathscr L_0a_p$ into the quadratic form cancels the added linear and quadratic terms and proves the third. The approximation \eqref{II-eq:second-response-comparison} identifies the minima of these quadratic expressions with the reduced coefficients. Division by $S_\ast$ proves invariance. The parameter changes contribute no fourth-order term: $\widehat\varepsilon/\varepsilon=1+O(\varepsilon^2)$, the quartic coefficient is continuous in the center, and the quadratic and cubic coefficients vanish throughout the locally totally umbilic neighborhood.

\smallskip\noindent \emph{The normalized representative}.
Conformal Fermi coordinates were constructed by Marques \cite[Proposition~3.1]{Marques2005}. The conformal-Fermi construction used here has $\varphi(p)=1$ and tangential $d\varphi(p)=0$; see \cite[Lemma~2.2]{KimMussoWei2021}. The mean-curvature transformation gives the normal derivative zero as well. Since the first derivatives of the original Fermi metric vanish, the corrections that normalize the volume density begin in degree two. For this fixed smooth metric, apply \cite[Proposition~2.2]{Almaraz2010} with $N=6$, giving a local metric $g_0$ with volume density $1+O(|x|^6)$ and $H_{g_0}=O(|x'|^5)$. Choose $\chi\in C_c^\infty([0,\infty))$ equal to one near zero. In these Fermi coordinates, the factor
\[
\psi=\exp(x_nf(x')\chi(x_n)),
\qquad f=(n-2)H_{g_0}/2=O(|x'|^5),
\]
has boundary value one and $\partial_{\nu_{g_0}}\psi=(n-2)H_{g_0}/2$. It makes $H$ vanish and, being $1+O(|x|^6)$, preserves the required Taylor coefficients. The resulting factor relative to $g$ has zero normal derivative on the boundary. Multiply its logarithm by a smooth cutoff equal to one near $p$ with zero boundary normal derivative, and extend by zero. The global positive factor preserves $H=0$ and the normalization at $p$. The energy isometry preserves strict coercivity.
\end{proof}

\subsection{The curvature formula and its sign}

\begin{theorem}[Exact formula for the fourth-order coefficient]
\label{II-thm:quartic-Weyl-sign}
Assume $n\ge7$, \textup{(Y$^+_{\partial}$)}, and \textup{(BG$^{5+}$)}. Let $g$ be the fixed representative with $H_g=0$, and suppose that the boundary is totally umbilic near $p$. Let $\widehat g=\varphi^{4/(n-2)}g$ be the normalized conformal-Fermi representative from Lemma~\ref{II-lem:quartic-gauge-completion}. In $\widehat g$-Fermi coordinates at $p$, put
\[
T_{ab}=R_{\widehat g,\nu a\nu b}(p),\qquad
\overline W=W_{\widehat g|_{T\partial M}}(p),\qquad
D_n=(n-1)(n-3)(n-4)(n-5)(n-6).
\]
Then $T$ is trace-free, and the standard-bubble coefficient is
\begin{equation}\label{II-eq:Xi-uncorrected-CF-explicit}
\Xi_{\widehat g}^{0}(p)=\frac{8-n}{D_n}|T|^2-b_n|\overline W|^2,
\qquad
b_n=\frac{n-2}{24(n-1)(n-3)(n-5)(n-6)}.
\end{equation}
The second-order source depends only on $T$; denote it by $\mathscr S_T$. For $T\ne0$, the coefficient of the corresponding quadratic form is
\begin{equation}\label{II-eq:kappa-n-quartic-model}
\kappa_n\coloneqq \frac1{2S_\ast}
\frac{\langle\mathscr L_0^{-1}\mathscr S_T,\mathscr S_T\rangle}{|T|^2}
=\frac{3n^3+30n^2-79n+6}
{3(n-2)(n-1)^2(n+1)(n-3)(n-4)(n-5)(n-6)}.
\end{equation}
Consequently,
\begin{equation}\label{II-eq:Xi-reduced-CF-explicit}
\Xi_g^{\mathrm{red}}(p)=\Xi_{\widehat g}^{\mathrm{red}}(p)
=-e_n|T|^2-b_n|\overline W|^2,
\end{equation}
where
\begin{equation}\label{II-eq:quartic-en-exact}
e_n=\frac{3n^3-18n^2+21n+14}
{3(n-2)(n-1)^2(n+1)(n-4)(n-5)(n-6)}>0.
\end{equation}
The optimal algebraic comparison with the ambient Weyl norm is
\begin{equation}\label{II-eq:Xi-Weyl-negative}
\Xi_g^{\mathrm{red}}(p)\le-c(n)|W_g(p)|_g^2,
\end{equation}
with
\begin{equation}\label{II-eq:quartic-cn-Weyl}
c(n)=\frac{3n^3-18n^2+21n+14}
{12(n-1)^2(n+1)(n-3)(n-4)(n-5)(n-6)}.
\end{equation}
Equality in \eqref{II-eq:Xi-Weyl-negative} holds precisely when $\overline W=0$. In particular,
\begin{equation}\label{II-eq:Xi-zero-iff-Weyl-zero}
\Xi_g^{\mathrm{red}}(p)=0\quad\Longleftrightarrow\quad W_g(p)=0.
\end{equation}
\end{theorem}

\begin{proof}
Here barred curvature tensors belong to the boundary metric, and semicolons denote covariant derivatives in the indicated directions. The conformal-Fermi identities at $p$ are
\[
\overline{\Ric}=0,
\qquad
\Ric(\nu,\nu)=0,
\qquad
\Scal=0,
\qquad
\Ric_{ab}=R_{\nu a\nu b},
\qquad
R_{\nu abc}=0,
\]
and
\[
\Scal_{;aa}
=-\frac16|\overline W|^2,
\qquad
\Ric_{\nu\nu;\nu\nu}+2|T|^2=0,
\qquad
T_{ab;ab}
=-\frac12\Scal_{;\nu\nu}-|T|^2;
\]
see \cite[Lemmas~2.4--2.5]{Almaraz2010}.

For the uncorrected expansion, use \cite[Section~3, (3-7)--(3-14)]{Almaraz2010} with correction parameter $A=0$. Fix a small chart radius $\delta>0$ and a smooth radial cutoff $\chi$ equal to one on $[0,1]$ and zero on $[2,\infty)$. In the $\widehat g$-Fermi chart put
\[
V_\varepsilon(x)=\left(\frac{\varepsilon}{|x'|^2+(\varepsilon+x_n)^2}\right)^{(n-2)/2},
\quad \psi_{\varepsilon,\delta}=\chi(|x|/\delta)V_\varepsilon,
\quad b_\varepsilon=\int_{\partial M}\psi_{\varepsilon,\delta}^q\,d\sigma_{\widehat g}.
\]
The flat bubble identity and integration by parts give
\[
\int_{B_\delta^+}|\nabla V_\varepsilon|^2-S_\ast b_\varepsilon^{2/q}
=O_\delta(\varepsilon^{n-2})+o(\varepsilon^4).
\]
The cutoff annulus contributes $O_\delta(\varepsilon^{n-2})$, and the volume density contributes $O_\delta(\varepsilon^5)$ under the order-six normalization. The Taylor errors are $O_\delta(\varepsilon^5\log(\delta/\varepsilon))$ for $n=7$ and $O_\delta(\varepsilon^5)$ for $n\ge8$. Replacing truncated fourth moments by their full limits adds $\varepsilon^4O((\delta/\varepsilon)^{6-n})$. These estimates also hold after $\mathbf D=\varepsilon\partial_\varepsilon$, so all errors are $o(\varepsilon^4)$ for $n\ge7$.
Thus, with
\[
I_n
\coloneqq 
\int_0^\infty\frac{r^n}{(1+r^2)^n}\,dr,
\qquad
J_n
\coloneqq 
\int_{\mathbb R^n_+}
\frac{|y'|^2}
{\bigl(|y'|^2+(1+y_n)^2\bigr)^{n-2}}\,dy,
\]
the curvature identities above give the following cancellation. If $D'_n=(n-3)(n-4)(n-5)(n-6)$, the three moments multiplying $T_{ab;ab}$, $|T|^2$, and $\Scal_{;\nu\nu}$, with $s=|y'|^2+(1+y_n)^2$, are
\begin{align*}
\int_{\mathbb R^n_+}\frac{y_n^2|y'|^4}{s^n}\,dy
&=\frac{2(n+1)|\mathbb S^{n-2}|I_n}{D'_n},\\
\int_{\mathbb R^n_+}\frac{y_n^4|y'|^2}{s^n}\,dy
&=\frac{24|\mathbb S^{n-2}|I_n}{(n-2)D'_n},\qquad\int_{\mathbb R^n_+}\frac{y_n^2}{s^{n-2}}\,dy
=\frac{8(n-2)|\mathbb S^{n-2}|I_n}{D'_n}.
\end{align*}
The intrinsic scalar-curvature term is
\[
\frac{n-2}{8(n-1)^2}\Scal_{;aa}J_n
=-\frac{n-2}{48(n-1)^2}J_n|\overline W|^2;
\]
see \cite[(3-14)]{Almaraz2010}. The coefficients of the three displayed moments are $(n-2)^2/((n-1)(n+1))$, $(n-2)^2/(2(n-1))$, and $(n-2)/(8(n-1))$, respectively. Substituting $T_{ab;ab}=-\Scal_{;\nu\nu}/2-|T|^2$ cancels $\Scal_{;\nu\nu}$ and yields
\begin{equation*}
\mathcal N_{\widehat g}^\circ
(\psi_{\varepsilon,\delta})
-S_\ast b_\varepsilon^{2/q}
=
\varepsilon^4\bigg[
\frac{2|\mathbb S^{n-2}|I_n(8-n)(n-2)}
{(n-1)(n-3)(n-4)(n-5)(n-6)}
|T|^2
-\frac{n-2}{48(n-1)^2}
J_n|\overline W|^2
\bigg]
+\mathscr E_{\varepsilon,\delta},
\end{equation*}
where
\[
|\mathscr E_{\varepsilon,\delta}|
+
|\mathbf D\mathscr E_{\varepsilon,\delta}|
=
o(\varepsilon^4).
\]
The flat energy is $2(n-2)|\mathbb S^{n-2}|I_n$, and the beta integrals give
\[
\frac{J_n}{|\mathbb S^{n-2}|I_n}
=\frac{4(n-1)(n-2)}{(n-3)(n-5)(n-6)}.
\]
Division by the flat energy gives
\eqref{II-eq:Xi-uncorrected-CF-explicit}. The fourth-order moments converge for $n\ge7$. Homogeneity of the quotient and the tail estimates in Proposition~\ref{II-prop:quartic-reduced-umbilic} transfer the value and weighted scale derivative to the $R_\beta(\varepsilon)=\varepsilon^{-\beta}$ family.

The second-order inverse-metric term is
\[
(g^{ab})_2
=
\frac13\overline R^{a}{}_{c}{}^{b}{}_{d}y_cy_d
+
T^{ab}y_n^2.
\]
The intrinsic term vanishes after contraction with the tangential derivative of the radial bubble. The remaining second-order terms vanish in these coordinates. Hence
\begin{equation}\label{II-eq:quartic-T-source-explicit}
\mathscr S_T[\phi]
=
2t_\ast
\int_{\mathbb R^n_+}
y_n^2T^{ab}
\partial_aU_+
\partial_b\phi\,dy
\qquad(\phi\in X_+),\qquad t_\ast=S_\ast^{1/2}.
\end{equation}
The factor $t_\ast$ comes from the trace-normalized base $U_\ast=t_\ast U_+$; the model direction $\phi$ is rescaled without this amplitude. The same integral defines $\mathscr S_T$ on the whole Dirichlet space. Its continuity follows from $\|y_n^2\nabla U_+\|_{L^2}<\infty$, which holds for $n\ge7$. The normal to the flat trace constraint is a multiple of $U_+$ in the Dirichlet pairing. Both $U_+$ and the dilation field $Z_0^\flat$ are radial in $y'$. For either of these functions, the tangential angular integral contains
\[
\int_{\mathbb S^{n-2}}T_{ab}\theta_a\theta_b\,d\theta
=\frac{|\mathbb S^{n-2}|}{n-1}\operatorname{tr}T=0.
\]
For a translation field $Z_c^\flat=-\partial_cU_+$ the integrand is odd under $y'\mapsto-y'$. Thus
\[
\mathscr S_T[U_+]=\mathscr S_T[Z_0^\flat]
=\mathscr S_T[Z_c^\flat]=0\qquad(1\le c\le n-1).
\]
These are the stated orthogonality conditions.

The equivariance follows from a change of variables. For $R\in O(n-1)$, put $(\mathcal R_R\phi)(y',y_n)=\phi(R^{-1}y',y_n)$. Then
\[
\mathscr S_{RTR^{\mathsf T}}[\mathcal R_R\phi]
=\mathscr S_T[\phi].
\]
The Dirichlet form, the radial boundary potential in $\mathscr L_0$, and the subspace $X_+$ are invariant under $\mathcal R_R$. Hence $\mathscr L_0^{-1}$ commutes with this action, and the quadratic form
\[
T\longmapsto
\langle\mathscr L_0^{-1}\mathscr S_T,\mathscr S_T\rangle
\]
is a multiple of $|T|^2$ by the invariant quadratic-form argument in Proposition~\ref{I-prop:spin2-schur-identification}. The source is nonzero for $T\ne0$, as is seen by pairing it with a compactly supported function of the form $Q_T h(|y'|,t)$. Positivity of the flat form therefore makes $\kappa_n$ positive. To evaluate it, it suffices by homogeneity to assume $|T|=1$. Set $a=1+t$, $s=|y'|^2+a^2$, and
\[
f_2(v)=n(n-2)A_n\int_{\R^n_+}t^2Q_Ts^{-(n+2)/2}v\,dy.
\]
Tangential integration by parts in \eqref{II-eq:quartic-T-source-explicit} gives $\mathscr S_T=-2t_\ast f_2$. Since $\mathscr L_0=2\mathcal B_n$ and $t_\ast^2=S_\ast$,
\[
\kappa_n=f_2(\mathcal B_n^{-1}f_2).
\]
Let $h_0,b_0,R_n$ be as in Lemma~\ref{I-lem:Robin-pairing-exact}, with $Q=Q_T$, and put
\[
p_2=A_nQ_Ts^{-n/2}
\left[\frac{n-2}{6}(a^2-3a+3)+\frac{s}{12}\right],
\qquad K_2=\frac{(n-3)A_n}{6}.
\]
Direct differentiation shows that $-\Delta p_2$ is the density defining $f_2$ and $\mathcal Rp_2=K_2b_0$. Hence
\[
\mathcal B_n^{-1}f_2=p_2-K_2h_0.
\]
The function $p_2$ has finite Dirichlet energy for $n\ge7$. Define
\[
I_0=\frac1{n-6}-\frac5{n-5}+\frac{10}{n-4}-\frac9{n-3}+\frac3{n-2},\quad I_1=\frac1{n-6}-\frac2{n-5}+\frac1{n-4}.
\]
The angular identity from Lemma~\ref{I-lem:Robin-pairing-exact} and radial beta integrals give
\begin{align*}
f_2(p_2)=&\frac{n-2}{12(n-1)}I_0+\frac{n}{12(n-3)(n-1)}I_1,\\
\int_{\R^{n-1}}K_2b_0p_2=\frac1{36(n-1)}+&\frac1{36(n-2)(n-5)},\quad K_2^2R_n=\frac{n(n-3)(n^2-2n-5)}
{36(n-2)(n-1)^2(n+1)(n-4)}.
\end{align*}
Substitution in \eqref{I-eq:Robin-Green-pairing} proves \eqref{II-eq:kappa-n-quartic-model}. Subtracting $\kappa_n$ from $(8-n)/D_n$ gives $-e_n$ and proves \eqref{II-eq:Xi-reduced-CF-explicit}. The numerator of $e_n$, written with $x=n-7\ge0$, is $3x^3+45x^2+210x+308$, so $e_n>0$.

Finally, the conformal-Fermi identities give
\[
W_{\nu a\nu b}=\frac{n-3}{n-2}T_{ab},\qquad
W_{\nu abc}=0,\qquad
|W_{\widehat g}|^2=|\overline W|^2+\frac{4(n-3)}{n-2}|T|^2.
\]
Thus the best algebraic constant is the smaller of $b_n$ and $e_n(n-2)/(4(n-3))$. The latter is strictly smaller because
\[
b_n-\frac{n-2}{4(n-3)}e_n
=\frac{n^2-3n-2}{24(n-5)(n-4)(n-1)^2(n+1)}>0
\qquad(n\ge7).
\]
This proves \eqref{II-eq:quartic-cn-Weyl} and the equality statement. Since $\varphi(p)=1$, the Weyl norm at $p$ is unchanged. Lemma~\ref{II-lem:quartic-gauge-completion} identifies the two reduced coefficients, completing the proof.
\end{proof}

\begin{corollary}[Intrinsic formula and conformal weight]
\label{II-cor:quartic-conformal-weight}
Under the hypotheses of Theorem~\ref{II-thm:quartic-Weyl-sign}, let $E_g=W_g(\nu,\cdot,\nu,\cdot)$ on $T_p\partial M$, and let $\overline W_g$ be the intrinsic Weyl tensor of the boundary metric. Then
\begin{equation}\label{II-eq:quartic-intrinsic}
\Xi_g^{\mathrm{red}}(p)
=-e_n\frac{(n-2)^2}{(n-3)^2}|E_g(p)|_g^2
-b_n|\overline W_g(p)|_g^2.
\end{equation}
If $\widetilde g=\psi^{4/(n-2)}g$ is another representative with minimal boundary, then
\begin{equation}\label{II-eq:quartic-conformal-weight}
\Xi_{\widetilde g}^{\mathrm{red}}(p)
=\psi(p)^{-8/(n-2)}\Xi_g^{\mathrm{red}}(p).
\end{equation}
\end{corollary}

\begin{proof}
In the normalized conformal-Fermi representative, $E=(n-3)T/(n-2)$. The normalizing factor has value one at $p$, so conformal covariance identifies the Weyl tensors and their norms there with those of $g$. Substitution in \eqref{II-eq:Xi-reduced-CF-explicit} gives \eqref{II-eq:quartic-intrinsic}. Under $\widetilde g=\psi^{4/(n-2)}g$, both squared norms in this formula acquire the factor $\psi(p)^{-8/(n-2)}$. This proves \eqref{II-eq:quartic-conformal-weight}.
\end{proof}

\begin{corollary}[Fourth-order subcriticality and compactness]
\label{II-cor:quartic-subcriticality-compactness}
Assume the hypotheses of Proposition~\ref{II-prop:quartic-reduced-umbilic} and work in a global representative satisfying $H_g=0$.
\begin{enumerate}[label=\textup{(\alph*)},leftmargin=2em]
\item If $p\in\mathcal C$ is umbilic and $\Xi_g^{\mathrm{red}}(p)<0$, then $C^*_{\Esc}(M,g)<S_\ast$.

\item This conclusion holds in particular if the boundary is totally umbilic near $p$ and $W_g(p)\ne0$.

\item Suppose that the boundary is totally umbilic. If a normalized sequence of constrained critical points satisfies
\[
\mathcal J_g[u_\ell]\longrightarrow S_\ast
\]
and is not precompact, then every limiting blow-up center $p$ satisfies $W_g(p)=0$. If $\inf_{\partial M}|W_g|_g>0$, then, for every $c_0<c_{\mathrm{mix}}$, the normalized constrained critical points with values at most $c_0$ form a compact subset of $H^1(M)$.
\end{enumerate}
\end{corollary}

\begin{proof}
For (a), \eqref{II-eq:quartic-reduced-value} gives $\mathcal J_{1,\beta}(p,\varepsilon)<S_\ast$ for small $\varepsilon$. Part (b) follows from Theorem~\ref{II-thm:quartic-Weyl-sign}.

For (c), Lemma~\ref{II-lem:beta-cutoff-relation} gives centers $x_\ell\to p$ and scales $\varepsilon_\ell\downarrow0$ with
\[
0=S_\ast^{-1}\partial_\varepsilon\mathcal J_{1,\beta}(x_\ell,\varepsilon_\ell)
=4\Xi_g^{\mathrm{red}}(x_\ell)\varepsilon_\ell^3+o(\varepsilon_\ell^3).
\]
Total umbilicity makes the expansion uniform in $x_\ell$. Dividing by $\varepsilon_\ell^3$ gives $\Xi_g^{\mathrm{red}}(p)=0$, hence $W_g(p)=0$ by \eqref{II-eq:Xi-zero-iff-Weyl-zero}. This excludes blow-up at $S_\ast$ if the Weyl norm is bounded away from zero. The compactness assertion below $c_{\mathrm{mix}}$ then follows from Theorems~\ref{I-thm:compactness-escobar} and~\ref{I-thm:signfree-first-gap-compactness}, as in Corollary~\ref{II-cor:threshold-compactness}.
\end{proof}

\begin{remark}[Boundary geometry and dimension]
At a general umbilic point, tangential derivatives of $\mathring{\mathrm{II}}$ enter the second-order source. They vanish under local total umbilicity, where $W_g(p)=0$ also gives $\Xi_g^{\mathrm{red}}(p)=0$. The radial integral $\int_1^R r^{5-n}\,dr$ converges as $R\to\infty$ for $n\ge7$ and grows logarithmically in dimension six.
\end{remark}

\appendix

\section{Supplementary stability estimates}
\label{A-sec:supplementary-stability}

The coercivity and modulation results of Section~\ref{I-sec:coercivity-first} also yield deficit estimates.

\subsection{The optimizer orbit and flat trace stability}
Retain $F_\ast$, $\mathcal M_S$, and the energy inner product from Section~\ref{I-sec:coercivity-first}. On the hemisphere, $\mathcal S$ denotes the unit boundary $L^q$-norm constraint.

\begin{lemma}[Modulation near the optimizer manifold on $S^n_+$]\label{I-lem:modulation-sphere}
There are $0<\delta_*<\delta_{\rm tube}$ and a $C^1$ map assigning to each $u\in\mathcal S$ with $\operatorname{dist}_{H^1}(u,\mathcal M_S)\le\delta_*$ an optimizer $v^\sharp\in\mathcal M_S$ satisfying $u-v^\sharp\perp T_{v^\sharp}\mathcal M_S$ in the energy inner product. It is the unique optimizer satisfying this orthogonality and $\|u-v^\sharp\|_{H^1}<\delta_{\rm tube}$. Moreover, the parameter map is $C^1$ and $\|u-v^\sharp\|_{H^1}\lesssim \mathrm{dist}_{H^1}(u,\mathcal M_S)$.
\end{lemma}

\begin{proof}
Use the width convention
\[
U^*_{\lambda,\xi}(y',t)=t_\ast\lambda^{-\alpha}
U_+\bigl((y'-\xi)/\lambda,t/\lambda\bigr),
\qquad \alpha=(n-2)/2,
\]
and let $F_{\lambda,\xi}$ be its Cayley pullback. Impose
$G_j(u,p)=\langle u-F_p,Z_j(p)\rangle_E=0$, where
$Z_0=\lambda\partial_\lambda F_p$ and $Z_a=\lambda\partial_{\xi^a}F_p$.
At $u=F_p$, the normalized parameter Jacobian is the negative energy Gram matrix. For the Gram calculation, identify $F_p$ and $Z_j$ with their conformal images on the half-space. At $(\lambda,\xi)=(1,0)$ the ratios $Z_j/U_\ast$ are $\alpha$ times the stereographic coordinate functions on $S^{n-1}$, up to signs. Since the normalized boundary measure is the uniform spherical measure and $\int x_ix_j=\delta_{ij}/n$, the Jacobi equation gives
\[
\langle Z_i,Z_j\rangle_E
=(q-1)S_\ast\int_{\partial\R^n_+} U_\ast^{q-2}Z_iZ_j\,dy'
=\frac{S_\ast(n-2)}4\delta_{ij}.
\]
The normalized first two parameter derivatives have uniform energy bounds. The quantitative implicit-function theorem therefore gives a fixed local radius, a unique $C^1$ parameter choice, and the distance estimate. Conformal energy isometries make the constants uniform along the orbit; fixed energy--$H^1$ equivalence transfers the result to the stated norm.

For uniqueness throughout the stated neighborhood, send one candidate to $F_\ast$ by a conformal map. When the relative translation and dilation parameters leave every compact set, the corresponding profiles converge weakly to zero in $\mathcal D^{1,2}$, so their energy pairing with $F_\ast$ tends to zero and their distance from it stays positive. On a compact relative parameter set, continuity and injectivity give the same separation outside the implicit-function neighborhood. Two sufficiently close candidates must therefore lie in that neighborhood and coincide. Applying the local estimate at approximate closest optimizers gives $\|u-v^\sharp\|_{H^1}\le C\operatorname{dist}_{H^1}(u,\mathcal M_S)$. Choose $\delta_*$ small enough to place the result inside the neighborhood where uniqueness holds.
\end{proof}

\begin{proposition}[Model trace stability]
\label{I-prop:deficit-to-proximity}\label{I-thm:local-esc-sphere}
\label{I-cor:local-esc-halfspace}\label{I-lem:local-H1-from-Dirichlet}
Let $n\ge3$. There are $\delta_0,c_*>0$, depending only on $n$, such that every positive $u\in\mathcal S$ on the round hemisphere with $\Def^{\Esc}(u)\le\delta_0$ admits the modulated optimizer $v^\sharp$ of Lemma~\ref{I-lem:modulation-sphere}, and
\[
u-v^\sharp\perp T_{v^\sharp}\mathcal M_S,
\qquad
\Def^{\Esc}(u)\ge c_*\|u-v^\sharp\|_{E,g_{\mathrm{round}}}^2
\asymp\|u-v^\sharp\|_{H^1(S^n_+)}^2.
\]
In particular, arbitrarily small deficit forces arbitrarily small distance to $\mathcal M_S$.

Write $\Def^{\Esc}_{\mathrm{flat}}(f)=\|\nabla f\|_2^2/\|f\|_{L^q(\partial\R^n_+)}^2-S_\ast$. If $f\in\mathcal D^{1,2}(\R^n_+)$ is positive, has the same boundary norm as $U_+$, and has sufficiently small deficit, there is $c_{\mathrm{tr}}=c_{\mathrm{tr}}(n)>0$ such that the corresponding optimizer $V^\sharp$ with this boundary norm satisfies
\[
\Def^{\Esc}_{\mathrm{flat}}(f)\ge c_{\mathrm{tr}}
\|\nabla(f-V^\sharp)\|_{L^2(\R^n_+)}^2,
\qquad
\|f-V^\sharp\|_{H^1(B_R^+)}^2
\le C(n)(1+R^2)\Def^{\Esc}_{\mathrm{flat}}(f)
\quad(R\ge1).
\]
\end{proposition}

\begin{proof}
Ho's trace stability theorem \cite[Theorem~1.2]{Ho2022} gives
\[
\inf_{V\in\mathcal M_{\mathrm{flat}}}\|\nabla(f-V)\|_2^2
\le C\bigl(\|\nabla f\|_2^2-S_\ast\|f\|_q^2\bigr)
=C\|f\|_q^2\Def^{\Esc}_{\mathrm{flat}}(f),
\]
where $\mathcal M_{\mathrm{flat}}$ is the full real-amplitude family of half-space optimizers.
If $V$ is an approximate closest optimizer, the trace inequality gives $|\|V\|_q-\|f\|_q|\le C\|\nabla(f-V)\|_2$. Since $\|\nabla V\|_2=\sqrt{S_\ast}\|V\|_q$, rescaling $V$ to the prescribed boundary norm changes the distance by at most a fixed multiple. For positive $f$ and small deficit, $V$ is positive. Lemma~\ref{I-lem:graph-norm} transfers the estimate to the hemisphere, and Lemma~\ref{I-lem:modulation-sphere} gives the orthogonal optimizer with comparable distance. Finally,
\[
\|w\|_{L^2(B_R^+)}\le CR\|\nabla w\|_{L^2(\R^n_+)}
\qquad(w\in\mathcal D^{1,2}(\R^n_+))
\]
by Sobolev and H\"older, proving the local $H^1$ estimate.
\end{proof}

\subsection{A deficit estimate near a boundary bubble}
We return to the fixed manifold $(M,g)$ and its boundary constraint $\mathcal S$.

\begin{lemma}[Taylor expansion under the boundary normalization]\label{I-lem:quotient-cubic-remainder}
Let $v\in\mathcal S$ be nonnegative and satisfy $\|v\|_{E,g}\le C_0$. For $h\in T_v\mathcal S$ with sufficiently small $\|h\|_{H^1}$, write $\mathcal R_v(h)=\mathsf{ret}_v(h)$. Uniformly on this set of base functions, and in particular for the normalized boundary bubbles $v=U_{+,x,\varepsilon}/\|U_{+,x,\varepsilon}\|_{L^q(\partial M)}$, $0<\varepsilon\le\varepsilon_0$, one has
\[
\mathcal J_g[\mathcal R_v(h)]
=\mathcal J_g[v]+L_{\mathcal S}[v](h)+\frac12Q_{\mathcal S,v}[h]+\mathcal R_v^{(3)}(h),
\]
where $L_{\mathcal S}[v]=D(\mathcal J_g\circ\mathsf{ret}_v)(0)$ and $Q_{\mathcal S,v}=D^2(\mathcal J_g\circ\mathsf{ret}_v)(0)$ on $T_v\mathcal S$, and
\begin{equation}\label{I-eq:quotient-cubic-remainder-chart}
|\mathcal R_v^{(3)}(h)|\le C\|h\|_{E,g}^{\min\{3,q\}}.
\end{equation}
For an arbitrary $u\in\mathcal S$ sufficiently close to $v$, the inverse chart is
\[
h=\frac{u}{\int_{\partial M}v^{q-1}u\,d\sigma_g}-v\in T_v\mathcal S,
\qquad \mathsf{ret}_v(h)=u.
\]
If $w=u-v$, then $h-w=O_{H^1}(\|w\|_{H^1}^2)$.
\end{lemma}

\begin{proof}
For $2<q\le3$, the trace Taylor estimate used in the proof of Lemma~\ref{II-lem:holder-quadratic-expansion} gives a remainder $O(\|h\|_E^q)$. For $q\ge3$, the pointwise remainder is bounded by $C(|v|+|h|)^{q-3}|h|^3$; H\"older's inequality gives $O(\|h\|_E^3)$ on the stated bounded set. Tangency removes the linear term in $\int|v+h|^q$. The quadratic energy numerator and the smooth reciprocal power of this boundary integral prove \eqref{I-eq:quotient-cubic-remainder-chart}. The inverse formula is the one used in Lemma~\ref{B-lem:normalized-modulation}; applying the same boundary Taylor estimate to $u=v+w\in\mathcal S$ gives $\int v^{q-1}u=1+O(\|w\|_{H^1}^2)$ and the claimed error.
\end{proof}

\begin{theorem}[Quantitative estimate near a boundary bubble]\label{I-thm:quant-escobar-M-cov}
Assume $n\ge5$, \textup{(Y$^+_{\partial}$)}, and \textup{(BG$^{2+}$)}. There exist $\delta_0>0$, $\varepsilon_0>0$, $C>0$, and $C_{\mathrm{geo}}>0$, depending only on $(M,g,n)$, with the following property.

Suppose $u\in\mathcal S$ and $\Def^{\mathrm{Esc}}_g(u)\le\delta_0$. Assume also that $u$ lies in the neighbourhood considered in Lemma~\ref{I-lem:modulation} and that the corresponding boundary bubble $v^\sharp=a^\sharp U_{+,x^\sharp,\varepsilon^\sharp}$ has $0<\varepsilon^\sharp<\varepsilon_0$. Then
\begin{equation}\label{I-eq:quant-escobar-M-cov}
\|u-v^\sharp\|_{E,g}^2\le C\Def^{\mathrm{Esc}}_g(u)+C_{\mathrm{geo}}(\varepsilon^\sharp)^2,
\end{equation}
and consequently
\[
\|u-v^\sharp\|_{H^1(M)}^2\lesssim\Def^{\mathrm{Esc}}_g(u)+C_{\mathrm{geo}}(\varepsilon^\sharp)^2.
\]
\end{theorem}

\begin{proof}
Since $C^*_{\Esc}\le S_\ast$, the small deficit bounds $\|u\|_E^2\le S_\ast+\delta_0$. Write $U=U_{+,x^\sharp,\varepsilon^\sharp}$, $e=\varepsilon^\sharp$, $r=u-v^\sharp$, and $\lambda=\|v^\sharp\|_q$. By shrinking the modulation neighborhood, we may assume that $r$ has small $H^1$ norm; $|\lambda-1|\le\|r\|_q$ then bounds $\lambda$ above and below. Put $v=v^\sharp/\lambda=tU\in\mathcal S$. Boundary orthogonality gives
\[
h=\frac{u}{\lambda}-v=\frac{r}{\lambda}\in T_v\mathcal S,
\quad u=\mathsf{ret}_v(h),\quad \|h\|_E\asymp\|r\|_E,
\]
and preserves the energy orthogonality to $\partial_eU$ and $\partial_{x^a}U$. Homogeneity and Theorem~\ref{I-thm:intrinsic-single-gap} give
\[
D^2(\mathcal J_g\circ\mathsf{ret}_v)(0)[h,h]
=t^{-2}D^2\mathcal J_g[U][h,h]\ge c\|h\|_E^2.
\]
Lemma~\ref{I-lem:first-order-drift} bounds the linear term by $Ce\|h\|_E$. Lemma~\ref{I-lem:quotient-cubic-remainder} and Young's inequality, with $p=\min\{3,q\}>2$, therefore yield
\[
\mathcal J_g[u]-\mathcal J_g[v]
\ge c'\|h\|_E^2-Ce\|h\|_E-C\|h\|_E^p
\ge\frac{c'}2\|h\|_E^2-C'e^2.
\]
Now $\mathcal J_g[v]\ge C^*_{\Esc}(M,g)$ proves \eqref{I-eq:quant-escobar-M-cov}; energy--$H^1$ equivalence proves the final estimate.
\end{proof}

\begin{remark}[Zero-trace comparison families]\label{I-lem:no-interior-channel}
Let $K_0=\ker(\operatorname{Tr}:H^1(M)\to L^q(\partial M))=H^1_0(M)$, and let $\mathcal C_{\mathrm{int}}$ be any nonempty subset of $K_0$. Interior-supported cutoff profiles and their scalar multiples are examples. If $\|u\|_{L^q(\partial M)}=1$, then for every $W\in K_0$ the trace inequality gives
$1=\|u-W\|_{L^q(\partial M)}\le C_{\mathrm{tr}}\|u-W\|_{H^1}$. Hence
\begin{equation}\label{I-eq:int-trace-bound}
\dist_{H^1}(u,\mathcal C_{\mathrm{int}})
\ge\dist_{H^1}(u,K_0)\ge C_{\mathrm{tr}}^{-1}.
\end{equation}
Interior harmonic profiles of Palais--Smale sequences vanish by Proposition~\ref{I-prop:signfree-boundary-profile-decomposition}.
\end{remark}

\begin{remark}[Minimizing sequences and profile normalization]\label{I-rem:interior-exclusion-threshold}\label{B-recall-I-par:moments-first}
For a non-precompact minimizing sequence with $Y=S_\ast$, Ekeland's principle and Lemma~\ref{I-lem:extremal-profile-onesigned} give, after passing to a subsequence, one boundary bubble. Below $S_\ast$, minimizing sequences are precompact by Corollary~\ref{I-cor:escobar-precompact-cov}. For the profile moments of our Dirichlet-normalized $U_+$, including their relations, first-order cancellation, and the logarithmic term in dimension three, see \cite[Lemma~3.9, Equation~(3.16), Remark~3.10, and Lemma~4.2]{EscobarGeometry}.
\end{remark}

\section{Hemisphere comparison and coefficient consequences}\label{I-subsec:closure}

\begin{proposition}[Hemisphere comparison and rigidity]
\label{I-prop:GC-known}\label{I-def:GC}
Let $(M^n,g)$ be a smooth compact connected Riemannian manifold with nonempty boundary, $n\ge3$, satisfying \textup{(Y$^+_\partial$)}. Then either $(M,g)$ is conformally diffeomorphic to the round hemisphere or
\[
C^*_{\Esc}(M,g)<S_\ast.
\]
\end{proposition}

\begin{proof}
The cases $n=3$, and $n=4,5$ with umbilic boundary, follow from Escobar~\cite{EscobarAnnals92}. A non-umbilic point gives strict inequality by the Escobar--Marques comparison; see \cite[pp.~383--384]{Marques2007} and \cite[pp.~1--2]{Almaraz2010}. We may therefore assume $n\ge6$ and total umbilicity.

Put $d=\lfloor(n-2)/2\rfloor$. If $\nabla^jW_g(p)\ne0$ for some boundary point $p$ and $0\le j\le d-2$, Chen's local comparison gives the result \cite[Theorem~2]{Chen2010}. Otherwise fix $p\in\partial M$. Pass to the conformal-Fermi representative with totally geodesic boundary used in \cite[Sections~2 and~4]{Chen2010}. In these coordinates $g_{in}=\delta_{in}$ and $\partial_ng_{ab}=0$ on the boundary. The vanishing Weyl jets imply that the Taylor coefficients of $\log g$ through degree $d$ vanish, as in the proof of \cite[Theorem~3]{Chen2010}. Hence
\[
g-\delta=O_2(r^{d+1}),\qquad \Scal_g=O(r^{d-1}),\qquad r=|x|,
\]
where $T=O_2(r^a)$ means $|\partial^jT|\le Cr^{a-j}$ for $0\le j\le2$.

Let $G$ be the positive conformal Green function with singular coefficient one. It has homogeneous Neumann data away from $p$. Equation~(12) of \cite{Chen2010} gives $F\coloneqq G-r^{2-n}=O(r^{d+3-n})$. On a half-annulus of radius $\rho$,
\[
L_g^\circ F=-L_g^\circ r^{2-n}=O(r^{d+1-n}),
\qquad \partial_nF=0\quad\text{on the flat boundary}.
\]
After rescaling, $F_\rho(X)=\rho^{n-d-3}F(\rho X)$ has bounded supremum norm on fixed half-annuli. Smoothness and the vanishing metric jets give uniformly elliptic coefficients and a uniformly bounded $C^{0,\alpha}$ right-hand side there, for any fixed $0<\alpha<1$. Local interior and Neumann boundary Schauder estimates on smaller half-annuli give $|\partial^jF|\le Cr^{d+3-n-j}$ for $j=1,2$. Thus $G=r^{2-n}(1+O_2(r^{d+1}))$. Inversion $y=x/|x|^2$ now gives
\[
\widehat g\coloneqq G^{4/(n-2)}g=\delta+O_2(|y|^{-d-1}),
\qquad d+1>(n-2)/2.
\]
This metric is complete and scalar-flat, with minimal boundary and one asymptotically half-space end. Its mass is nonnegative, and zero mass forces the Euclidean half-space, by \cite[Theorem~2.10]{GongKimMussoWei2026}. See also \cite[Theorem~1.1]{KhuriWangWang2026} and \cite[Section~4]{AlmarazBarbosaDeLima2016}. In the nonorientable case, the connected orientation cover has two half-space ends of the same mass, so rigidity excludes zero mass. Thus either conformal compactification gives the round hemisphere or the mass is positive. In the latter case Chen's Green-function test gives strict inequality \cite[Theorem~3 and Proposition~10]{Chen2010}.\end{proof}

\begin{lemma}[Equality at the hemisphere value]
\label{I-lem:attained-equality-rigidity}
Suppose that the hypotheses of Proposition~\ref{I-prop:GC-known} hold. Equality $C^*_{\Esc}(M,g)=S_\ast$ forces conformal equivalence with the round hemisphere. A positive smooth extremal $u$ with multiplier $\lambda$ gives $\widehat g=u^{4/(n-2)}g$ with
\[
\Scal_{\widehat g}=0,\qquad H_{\widehat g}=\frac{2\lambda}{n-2},
\qquad
\mathring{\mathrm{II}}_{\widehat g}=u^{2/(n-2)}\mathring{\mathrm{II}}_g.
\]
Thus total umbilicity is preserved; $\lambda=C^*_{\Esc}(M,g)$ if $\|u\|_q=1$.
\end{lemma}
\begin{proof}
Use Proposition~\ref{I-prop:GC-known} and conformal covariance of the equations $L_g^\circ u=0$, $B_g^\circ u=\lambda u^{q-1}$.
\end{proof}

\begin{remark}[Disjoint hemispheres]\label{I-rem:GC-status}
The disjoint union of two round hemispheres has $C^*_{\Esc}=S_\ast$ and is not conformally equivalent to one hemisphere.
\end{remark}

\begin{theorem}[Coefficient consequences and the Escobar threshold]
\label{I-thm:closure-escobar}\label{I-cor:model-LCF}\label{I-thm:zmc-second}
Assume $n\ge6$, \textup{(Y$^+_{\partial}$)}, \textup{(BG$^{4+}$)}, and $H_g=0$.
\begin{enumerate}[label=\textup{(\roman*)},leftmargin=2em]
\item A non-umbilic point gives $C^*_{\Esc}(M,g)<S_\ast$, attainment, and compactness of all normalized minimizing sequences and all normalized Palais--Smale sequences at the infimum. Either criterion $\mathfrak R_g^{\mathrm{bare}}(x)<0$ or $\mathfrak R_g^{\mathrm{bare}}(x)=0>\Theta_g(x)$ implies these conclusions, as do the corresponding criteria with $\mathfrak R_g^{\mathrm{red}}$. For the standard bubble, a vanishing quadratic coefficient and a negative cubic coefficient can occur only when $n=6$.
\item At every umbilic point, Lemma~\ref{I-lem:umbilic-lower-coefficients} gives the vanishing of the lower-order coefficients and the indicated derivative, and Proposition~\ref{I-prop:standard-bubble-inputs}\textup{(ii)} gives, with the stated cutoffs, the energy expansions with remainder $o(\varepsilon^3)$. At a non-umbilic point in dimension six, the quadratic coefficient of the standard bubble vanishes but the cubic coefficient need not vanish.
\item If $C^*_{\Esc}(M,g)=S_\ast$, the boundary is totally geodesic and the quantities in \textup{(ii)} vanish everywhere. At every limiting concentration point $p$ of a sequence of positive constrained critical points at $S_\ast$, one has $\mathrm{II}(p)=0$ and $\mathfrak R_g^{\mathrm{bare}}(p)=\mathfrak R_g^{\mathrm{red}}(p)=\nabla_\partial\mathfrak R_g^{\mathrm{red}}(p)=\Theta_g(p)=0$. Non-umbilic concentration is excluded for all $n\ge5$ by Theorem~\ref{II-thm:threshold-Sstar-onebubble}.
\item If $M$ is connected, $C^*_{\Esc}(M,g)=S_\ast$ implies conformal equivalence with the round hemisphere.
\end{enumerate}
\end{theorem}

\begin{proof}
The quadratic formulas in Propositions~\ref{I-prop:standard-bubble-inputs} and~\ref{I-rem:reduced-vs-bare}, and $\Theta_g=0$ at umbilic points, show that every criterion in \textup{(i)} forces non-umbilicity; the case in which the standard bubble has vanishing quadratic coefficient at a non-umbilic point requires $n=6$. Apply Theorem~\ref{I-thm:subcriticality-nonumbilic} and Corollary~\ref{I-cor:escobar-precompact-cov}.

The results cited in \textup{(ii)} use \cite[Proposition~4.7, Lemma~4.11, Equations~(4.25) and~(B.21), and Corollary~B.10]{EscobarGeometry}; for the non-umbilic cubic terms see \cite[Theorem~B.9 and Remark~B.11]{EscobarGeometry}. Equality forces $\mathring{\mathrm{II}}\equiv0$ by \textup{(i)}, proving \textup{(iii)} since $H_g=0$; \textup{(iv)} is Proposition~\ref{I-prop:GC-known}.
\end{proof}

\section{Genericity at finite regularity}
\label{app:genericity-proof}

Throughout this appendix, the space of metrics and the notation are those of Theorem~\ref{II-thm:generic-compactness}.

\begin{lemma}[Finite-regularity conformal and coordinate estimates]\label{I-lem:finite-regularity-preparation}
Let $g$ be a fixed $C^m$ metric on the compact manifold, with $m\ge3$. There is a positive $C^m$ factor $\phi$ with boundary value one for which $H_{\phi^{4/(n-2)}g}=0$. At $m=3$, the radial bubbles used above have two uniformly bounded normalized parameter derivatives in $H^1$ for $n\ge5$. The first-order metric expansion and its radial derivative satisfy the estimates used in the quadratic reduction.
\end{lemma}
\begin{proof}
In a fixed smooth collar write $\nu=a\partial_t+V$, where $V$ is tangential and $a>0$. Prescribe
\[
F|_{\partial M}=0,\qquad
\partial_tF=f\coloneqq \frac{n-2}{2a}H_g\in C^{m-1}(\partial M).
\]
In each boundary chart the jets $F=0$,
$\partial_\top^\alpha\partial_tF=\partial_\top^\alpha f$ for $|\alpha|\le m-1$, and all jets with at least two transverse derivatives equal to zero satisfy the Whitney compatibility conditions. The extension theorem \cite[Theorem~I, p.~65]{Whitney1934}, a boundary partition of unity, and a collar cutoff give a global $C^m$ function $F$. With $\phi=e^F$, the prescribed normal derivative is $(n-2)H_g/2$, so conformal covariance gives minimal boundary and preserves coercivity.

For the coordinate estimates write the principal Fermi coefficient as
$A(y)=I+A_1(y)+E(y)$, with $A_1$ the first jet determined by $\mathrm{II}(p)$. The Jacobi and Riccati integral equations, using $\operatorname{Rm}\in C^1$ and $\mathrm{II}\in C^2$, give
\begin{equation}\label{I-eq:finite-C3-Fermi-remainder}
|E|+|y\cdot\nabla E|\le C|y|^2,
\qquad |\nabla E|\le C|y|.
\end{equation}
These estimates also follow by smooth approximation in $C^3$ with uniform integral bounds. Scalar-curvature terms and their radial derivatives use only $\operatorname{Rm}\in C^1$.

For modulation, write the radial profile in a fixed collar using the normal distance $t(z)$, the boundary footpoint $\pi_\partial(z)$, and $d_{\bar g}(x,\pi_\partial(z))^2$. The first-variation formula
$\partial_xd_{\bar g}(x,z')^2=-2\exp_x^{-1}(z')$ gives joint $C^3$ regularity of squared distance near the diagonal; the footpoint is $C^2$. Two center derivatives and one spatial derivative therefore give the required rescaled $H^1$ bounds. For $R=\varepsilon^{-3/4}$ the cutoff contribution to an ordinary first parameter derivative has $H^1$ norm $O(\varepsilon^{(3n-14)/8})=o(1)$ for $n\ge5$; normalized derivatives through order two remain bounded. These bounds give the modulation, source, and first scale-derivative estimates for each fixed $C^3$ metric.
\end{proof}

\begin{proof}[Proof of Theorem~\ref{II-thm:generic-compactness}]
\smallskip\noindent \emph{Boundary transversality.}
For smooth metrics the transversality argument is due to Marques; see \cite[p.~343]{Almaraz2011CVPDE}.

Put $\mathcal M=\operatorname{Met}_{Y,BG}^m(M)$. Consider the section $(g,x)\mapsto\mathring{\mathrm{II}}_g(x)$ over $\mathcal M\times\partial M$. At a zero $(g,x_0)$, choose a smooth boundary chart whose transverse coordinate vector equals $\nu_g$ at $x_0$, with $g_{an}(x_0)=0$ and $g_{nn}(x_0)=1$. For $A\in\operatorname{Sym}_0^2(T_{x_0}^*\partial M,\bar g)$, choose cutoffs with $\chi(0)=1$, $\eta(0)=0$, $\eta'(0)=1$, and set
\[
h_{ab}\coloneqq -2\chi(x')\eta(x_n)A_{ab},
\qquad
h_{an}=h_{nn}\coloneqq 0.
\]
Since $h|_{\partial M}=0$, the boundary metric and normal have zero variation there. At $x_0$, the variation of the second fundamental form is therefore $-\tfrac12\partial_{x_n}h_{ab}=A_{ab}$. Because $A$ is trace-free, one has
\[
D_g\mathring{\mathrm{II}}_g(x_0)[h]=A.
\]
Put $d=n-1$ and $N=n(n-1)/2-1$. In local trace-free fiber trivializations, $F(g,x)=\mathring{\mathrm{II}}_g(x)$ is $C^{m-1}$, and the preceding variation makes $D_gF$ onto at every zero. Thus $\mathcal Z=F^{-1}(0)$ is a $C^{m-1}$ Banach submanifold, and the projection $\pi:\mathcal Z\to\mathcal M$ is Fredholm of index $d-N<0$. The space of $C^m$ metrics is separable, and these manifolds have countable bases. Since $m-1\ge2>\max\{d-N,0\}$, Smale's theorem \cite[Theorem~1.3]{Smale1965} gives a residual set of regular values. Since the Fredholm index is negative, every regular fiber is empty. Therefore $\mathcal G_1=\{g\in\mathcal M:\mathring{\mathrm{II}}_g\ne0\text{ on }\partial M\}$ is dense; compactness of $\partial M$ makes it open.

\smallskip\noindent \emph{An interior perturbation.}
Fix $y_0\in M^\circ$. In a smooth chart with $g_{ij}(y_0)=\delta_{ij}$, let $\mathbf W\ne0$ be an algebraic Weyl tensor and let $\chi\in C_c^\infty(M^\circ)$ be equal to one near $y_0$. Write $j^1_{y_0}h$ for the value and first derivatives of $h$ at $y_0$. The variation
\[
h^{\mathbf W}_{ij}(x)
\coloneqq 
-\tfrac13\mathbf W_{ikjl}x^kx^l\chi(x)
\]
satisfies
\[
j^1_{y_0}h^{\mathbf W}=0,
\qquad
D_gW_g(y_0)[h^{\mathbf W}]
=
\mathbf W.
\]
Only the second derivatives enter the curvature variation, because the value and first derivatives vanish. Hence $\mathcal G_2=\{g\in\mathcal M:W_g(y_0)\ne0\}$ is open and dense. Taking $\mathcal G=\mathcal G_1\cap\mathcal G_2$ gives \eqref{II-eq:generic-geometric-properties}. Since conformal flatness is preserved by conformal changes and diffeomorphisms, no metric in $\mathcal G$ is conformally diffeomorphic to the round hemisphere.

\smallskip\noindent \emph{Analytic consequences.}
Apply Lemma~\ref{I-lem:finite-regularity-preparation} to pass to a representative with minimal boundary while preserving the stated regularity and coercivity. Its coordinate bounds justify the quadratic reduction, including modulation and the first scale derivative at $m=3$. Theorem~\ref{I-thm:subcriticality-nonumbilic}, Corollary~\ref{I-cor:escobar-precompact-cov}, and Corollary~\ref{II-cor:threshold-compactness} give strict inequality, attainment, minimizing-sequence compactness, and compactness of the sets of constrained critical points with the stated energy bounds.
\end{proof}

\section{Supplementary one-bubble consequences}\label{B-sec:consequences}

\subsection{Recovery of the quadratic coefficient from two scales}

\begin{theorem}[Two-scale coefficient recovery]\label{I-thm:inverse-precise}
Assume $n\ge7$, \textup{(Y$^+_{\partial}$)}, \textup{(BG$^{4+}$)}, and $H_g=0$ on $\partial M$. For the standard normalized bubbles with the cutoff $R(\varepsilon)=\varepsilon^{-3/4}$ put
\[
\Delta_x(e)=\frac{\mathcal J_g[\mathcal U_{x,e}]-S_\ast}{S_\ast},
\qquad D_x(e)=\frac{\Delta_x(2e)-2\Delta_x(e)}{2e^2}.
\]
For all sufficiently small $e$, uniformly in $x\in\partial M$,
\[
|D_x(e)-\mathfrak R_g^{\mathrm{bare}}(x)|\le Ce,
\qquad
\left|\frac{\Delta_x(e)}e-eD_x(e)\right|\le Ce^2.
\]
In particular, these estimates recover the quadratic coefficient of the standard bubble at the modulated center of a concentrating sequence by evaluating the standard bubbles at its scale and twice its scale.
\end{theorem}

\begin{proof}
Put $\zeta_n=\min\{2,3(n-2)/4-2\}\ge7/4$. The expansion in \cite[Proposition~4.7 and Lemma~4.11]{EscobarGeometry}, with the cutoff tails of Proposition~\ref{I-prop:standard-bubble-inputs}, gives
\[
\Delta_x(e)=\mathfrak R_g^{\mathrm{bare}}(x)e^2+\Theta_g(x)e^3+O(e^{2+\zeta_n}).
\]
Substituting $e$ and $2e$ yields
\[
D_x(e)-\mathfrak R_g^{\mathrm{bare}}(x)=3\Theta_g(x)e+O(e^{\zeta_n}),
\qquad
\frac{\Delta_x(e)}e-eD_x(e)=-2\Theta_g(x)e^2+O(e^{1+\zeta_n}).
\]
The coefficients are uniformly bounded and $\zeta_n>1$.
\end{proof}

\subsection{Strictly mean-convex Euclidean domains}

\begin{theorem}[Compactness on strictly mean-convex Euclidean domains]\label{I-thm:convex-compactness}
Let $\Omega\subset\mathbb R^n$ be a bounded smooth domain, $n\ge4$, and let $g=g_{\mathrm{flat}}$. Assume that the boundary mean curvature with respect to the
\emph{inner} unit normal satisfies
\[
H_g(x)\ge h_0>0 \qquad \text{for all }x\in\partial\Omega.
\]
Assume moreover that $(\Omega,g_{\mathrm{flat}})$ is not conformally equivalent to the unit ball $B^n$ (equivalently, not conformally diffeomorphic to $(S^n_+,g_{\mathrm{round}})$). If $n\ge15$, assume in addition that the trace-free second fundamental form $\mathring{\mathrm{II}}$ of $\partial\Omega$ vanishes nowhere on $\partial\Omega$.

Let $u_\ell>0$ be normalized positive critical points of the Escobar quotient, i.e.
\[
L_g^\circ u_\ell=0 \ \text{ in }\Omega,\qquad
B_g^\circ u_\ell=\lambda_\ell u_\ell^{2^*_\partial-1}\ \text{ on }\partial\Omega,
\qquad
\int_{\partial\Omega}u_\ell^{2^*_\partial}d\sigma=1.
\]
Then $(u_\ell)$ is precompact in $C^2(\overline\Omega)$, hence in $H^1(\Omega)$. In particular, no boundary bubbling occurs.
\end{theorem}

\begin{proof}
The trace--Poincar\'e inequality and $H_g\ge h_0>0$ imply \textup{(Y$^+_\partial$)}. Hence $\mu\coloneqq C^*_{\Esc}(\Omega,g)>0$. Testing the equation gives $\lambda_\ell=\mathcal J_g[u_\ell]\ge\mu$. Put
\[
c_\ell=(\lambda_\ell/\mu)^{(n-2)/2},\qquad v_\ell=c_\ell u_\ell.
\]
Then $L_g^\circ v_\ell=0$ and $B_g^\circ v_\ell=\mu v_\ell^{q-1}$, where $q=2^*_\partial$. The fixed-coefficient compactness theorem \cite[Theorem~1]{GongKimMussoWei2026} applies for $4\le n\le14$, and \cite[Theorem~1.2]{Almaraz2011CVPDE} applies for $n\ge15$ under the nowhere-umbilic hypothesis. Thus a subsequence of $(v_\ell)$ converges in $C^2(\overline\Omega)$ to a positive solution. Since $c_\ell=\|v_\ell\|_{L^q(\partial\Omega)}$, these factors converge to a positive number, and $u_\ell=v_\ell/c_\ell$ converges in $C^2(\overline\Omega)$.
\end{proof}

\section*{Acknowledgments}
The authors thank IIT Bombay for providing ideal working conditions, and Saikat Mazumdar and Souptik Chakraborty for helpful discussions.

\bibliographystyle{alpha}
\bibliography{Escobar_II}

@article {Beckner1993,
    AUTHOR = {Beckner, William},
     TITLE = {Sharp {S}obolev inequalities on the sphere and the
              {M}oser-{T}rudinger inequality},
   JOURNAL = {Ann. of Math. (2)},
  FJOURNAL = {Annals of Mathematics. Second Series},
    VOLUME = {138},
      YEAR = {1993},
    NUMBER = {1},
     PAGES = {213--242},
      ISSN = {0003-486X,1939-8980},
   MRCLASS = {58G30 (46E35 53C21 58G26)},
  MRNUMBER = {1230930},
MRREVIEWER = {Paul\ C.\ Yang},
       DOI = {10.2307/2946638},
       URL = {https://doi.org/10.2307/2946638},
}

@book {AmbrosettiMalchiodi2006,
    AUTHOR = {Ambrosetti, Antonio and Malchiodi, Andrea},
     TITLE = {Perturbation methods and semilinear elliptic problems on
              {${\bf R}^n$}},
    SERIES = {Progress in Mathematics},
    VOLUME = {240},
 PUBLISHER = {Birkh\"{a}user Verlag, Basel},
      YEAR = {2006},
     PAGES = {xii+183},
      ISBN = {978-3-7643-7321-4; 3-7643-7321-0},
   MRCLASS = {35-02 (35B20 35B25 35B33 35J60 47J30 58E05)},
  MRNUMBER = {2186962},
MRREVIEWER = {Jes\'{u}s\ Hern\'{a}ndez},
       DOI = {10.1007/3-7643-7396-2},
       URL = {https://doi.org/10.1007/3-7643-7396-2},
}

@article {ChenFrankWeth2013,
    AUTHOR = {Chen, Shibing and Frank, Rupert L. and Weth, Tobias},
     TITLE = {Remainder terms in the fractional {S}obolev inequality},
   JOURNAL = {Indiana Univ. Math. J.},
  FJOURNAL = {Indiana University Mathematics Journal},
    VOLUME = {62},
      YEAR = {2013},
    NUMBER = {4},
     PAGES = {1381--1397},
      ISSN = {0022-2518,1943-5258},
   MRCLASS = {46E35 (26D15)},
  MRNUMBER = {3179693},
MRREVIEWER = {Jan\ Vyb\'{\i}ral},
       DOI = {10.1512/iumj.2013.62.5065},
       URL = {https://doi.org/10.1512/iumj.2013.62.5065},
}

@article {DengSunWei2025,
    AUTHOR = {Deng, Bin and Sun, Liming and Wei, Jun-Cheng},
     TITLE = {Sharp quantitative estimates of {S}truwe's decomposition},
   JOURNAL = {Duke Math. J.},
  FJOURNAL = {Duke Mathematical Journal},
    VOLUME = {174},
      YEAR = {2025},
    NUMBER = {1},
     PAGES = {159--228},
      ISSN = {0012-7094,1547-7398},
   MRCLASS = {35A23 (26D10 35B35 35J20)},
  MRNUMBER = {4867256},
       DOI = {10.1215/00127094-2024-0026},
       URL = {https://doi.org/10.1215/00127094-2024-0026},
}

@article {Gerard1998,
    AUTHOR = {G\'{e}rard, Patrick},
     TITLE = {Description du d\'{e}faut de compacit\'{e} de l'injection de
              {S}obolev},
   JOURNAL = {ESAIM Control Optim. Calc. Var.},
  FJOURNAL = {ESAIM. Control, Optimisation and Calculus of Variations.
              European Series in Applied and Industrial Mathematics},
    VOLUME = {3},
      YEAR = {1998},
     PAGES = {213--233},
      ISSN = {1292-8119,1262-3377},
   MRCLASS = {46E35 (49K10)},
  MRNUMBER = {1632171},
MRREVIEWER = {W.\ P.\ Ziemer},
       DOI = {10.1051/cocv:1998107},
       URL = {https://doi.org/10.1051/cocv:1998107},
}

@article {FigalliGlaudo2020,
    AUTHOR = {Figalli, Alessio and Glaudo, Federico},
     TITLE = {On the sharp stability of critical points of the {S}obolev
              inequality},
   JOURNAL = {Arch. Ration. Mech. Anal.},
  FJOURNAL = {Archive for Rational Mechanics and Analysis},
    VOLUME = {237},
      YEAR = {2020},
    NUMBER = {1},
     PAGES = {201--258},
      ISSN = {0003-9527,1432-0673},
   MRCLASS = {35A23 (26D10 35B35 35J20)},
  MRNUMBER = {4090466},
MRREVIEWER = {Stefan\ Steinerberger},
       DOI = {10.1007/s00205-020-01506-6},
       URL = {https://doi.org/10.1007/s00205-020-01506-6},
}

@article {GhimentiMichelettiPistoia2018,
    AUTHOR = {Ghimenti, Marco and Micheletti, Anna Maria and Pistoia,
              Angela},
     TITLE = {Linear perturbation of the {Y}amabe problem on manifolds with
              boundary},
   JOURNAL = {J. Geom. Anal.},
  FJOURNAL = {Journal of Geometric Analysis},
    VOLUME = {28},
      YEAR = {2018},
    NUMBER = {2},
     PAGES = {1315--1340},
      ISSN = {1050-6926,1559-002X},
   MRCLASS = {58J05 (35J60 35R01 53C21)},
  MRNUMBER = {3790502},
MRREVIEWER = {Qu\cfac{o}c\ Anh\ Ng\^{o}},
       DOI = {10.1007/s12220-017-9864-6},
       URL = {https://doi.org/10.1007/s12220-017-9864-6},
}

@article {Marques2005,
    AUTHOR = {Marques, Fernando C.},
     TITLE = {Existence results for the {Y}amabe problem on manifolds with
              boundary},
   JOURNAL = {Indiana Univ. Math. J.},
  FJOURNAL = {Indiana University Mathematics Journal},
    VOLUME = {54},
      YEAR = {2005},
    NUMBER = {6},
     PAGES = {1599--1620},
      ISSN = {0022-2518,1943-5258},
   MRCLASS = {53C21},
  MRNUMBER = {2189679},
MRREVIEWER = {Mohameden\ Ahmedou},
       DOI = {10.1512/iumj.2005.54.2590},
       URL = {https://doi.org/10.1512/iumj.2005.54.2590},
}

@article {Lions85a,
    AUTHOR = {Lions, P.-L.},
     TITLE = {The concentration-compactness principle in the calculus of
              variations. {T}he limit case. {I}},
   JOURNAL = {Rev. Mat. Iberoamericana},
  FJOURNAL = {Revista Matem\'{a}tica Iberoamericana},
    VOLUME = {1},
      YEAR = {1985},
    NUMBER = {1},
     PAGES = {145--201},
      ISSN = {0213-2230},
   MRCLASS = {49A22 (58E30)},
  MRNUMBER = {834360},
MRREVIEWER = {Ll.\ G.\ Chambers},
       DOI = {10.4171/RMI/6},
       URL = {https://doi.org/10.4171/RMI/6},
}

@article {LiZhu1997,
    AUTHOR = {Li, Yanyan and Zhu, Meijun},
     TITLE = {Sharp {S}obolev trace inequalities on {R}iemannian manifolds
              with boundaries},
   JOURNAL = {Comm. Pure Appl. Math.},
  FJOURNAL = {Communications on Pure and Applied Mathematics},
    VOLUME = {50},
      YEAR = {1997},
    NUMBER = {5},
     PAGES = {449--487},
      ISSN = {0010-3640,1097-0312},
   MRCLASS = {46E35 (58G30)},
  MRNUMBER = {1443055},
MRREVIEWER = {H.\ Triebel},
       DOI = {10.1002/(SICI)1097-0312(199705)50:5<449::AID-CPA2>3.3.CO;2-5},
       URL =
              {https://doi.org/10.1002/(SICI)1097-0312(199705)50:5<449::AID-CPA2>3.3.CO;2-5},
}

@article {Whitney1934,
    AUTHOR = {Whitney, Hassler},
     TITLE = {Analytic extensions of differentiable functions defined in
              closed sets},
   JOURNAL = {Trans. Amer. Math. Soc.},
  FJOURNAL = {Transactions of the American Mathematical Society},
    VOLUME = {36},
      YEAR = {1934},
    NUMBER = {1},
     PAGES = {63--89},
      ISSN = {0002-9947,1088-6850},
   MRCLASS = {26B05},
  MRNUMBER = {1501735},
       DOI = {10.2307/1989708},
       URL = {https://doi.org/10.2307/1989708},
}

@article{Mazumdar2017,
    AUTHOR = {Mazumdar, Saikat},
     TITLE = {Struwe's decomposition for a polyharmonic operator on a
              compact {R}iemannian manifold with or without boundary},
   JOURNAL = {Commun. Pure Appl. Anal.},
  FJOURNAL = {Communications on Pure and Applied Analysis},
    VOLUME = {16},
      YEAR = {2017},
    NUMBER = {1},
     PAGES = {311--330},
      ISSN = {1534-0392,1553-5258},
   MRCLASS = {58J60 (35B33 35J30 35J35 35R01)},
  MRNUMBER = {3583528},
       DOI = {10.3934/cpaa.2017015},
       URL = {https://doi.org/10.3934/cpaa.2017015}
}

@article{Almaraz2010,
 AUTHOR = {de Moura Almaraz, S\'ergio},
     TITLE = {An existence theorem of conformal scalar-flat metrics on
              manifolds with boundary},
   JOURNAL = {Pacific J. Math.},
  FJOURNAL = {Pacific Journal of Mathematics},
    VOLUME = {248},
      YEAR = {2010},
    NUMBER = {1},
     PAGES = {1--22},
      ISSN = {0030-8730,1945-5844},
   MRCLASS = {53C21 (35J65 35R01)},
  MRNUMBER = {2734161},
MRREVIEWER = {Gabjin\ Yun},
       DOI = {10.2140/pjm.2010.248.1},
       URL = {https://doi.org/10.2140/pjm.2010.248.1},
}

@article{Almaraz2011CVPDE,
 AUTHOR = {de Moura Almaraz, S\'ergio},
     TITLE = {A compactness theorem for scalar-flat metrics on manifolds
              with boundary},
   JOURNAL = {Calc. Var. Partial Differential Equations},
  FJOURNAL = {Calculus of Variations and Partial Differential Equations},
    VOLUME = {41},
      YEAR = {2011},
    NUMBER = {3-4},
     PAGES = {341--386},
      ISSN = {0944-2669,1432-0835},
   MRCLASS = {53C21 (35J65 35R01)},
  MRNUMBER = {2796235},
MRREVIEWER = {Seongtag\ Kim},
       DOI = {10.1007/s00526-010-0365-8},
       URL = {https://doi.org/10.1007/s00526-010-0365-8},
}

@article {Almaraz2014,
    AUTHOR = {de Moura Almaraz, S\'ergio},
     TITLE = {The asymptotic behavior of {P}alais-{S}male sequences on
              manifolds with boundary},
   JOURNAL = {Pacific J. Math.},
  FJOURNAL = {Pacific Journal of Mathematics},
    VOLUME = {269},
      YEAR = {2014},
    NUMBER = {1},
     PAGES = {1--17},
      ISSN = {0030-8730,1945-5844},
   MRCLASS = {35R01 (35J20 35J65 53C21)},
  MRNUMBER = {3233907},
MRREVIEWER = {Anna\ Maria\ Candela},
       DOI = {10.2140/pjm.2014.269.1},
       URL = {https://doi.org/10.2140/pjm.2014.269.1},
}

@article {AlmarazBarbosaDeLima2016,
    AUTHOR = {de Moura Almaraz, S\'ergio and Barbosa, Ezequiel and de Lima, Levi
              Lopes},
     TITLE = {A positive mass theorem for asymptotically flat manifolds with
              a non-compact boundary},
   JOURNAL = {Comm. Anal. Geom.},
  FJOURNAL = {Communications in Analysis and Geometry},
    VOLUME = {24},
      YEAR = {2016},
    NUMBER = {4},
     PAGES = {673--715},
      ISSN = {1019-8385,1944-9992},
   MRCLASS = {53C20 (53C21)},
  MRNUMBER = {3570413},
MRREVIEWER = {Mohameden\ Ahmedou},
       DOI = {10.4310/CAG.2016.v24.n4.a1},
       URL = {https://doi.org/10.4310/CAG.2016.v24.n4.a1},
}

@article{AlmarazQueirozWang2019,
   AUTHOR = {de Moura Almaraz, S\'ergio and de Queiroz, Olivaine S. and Wang,
              Shaodong},
     TITLE = {A compactness theorem for scalar-flat metrics on 3-manifolds
              with boundary},
   JOURNAL = {J. Funct. Anal.},
  FJOURNAL = {Journal of Functional Analysis},
    VOLUME = {277},
      YEAR = {2019},
    NUMBER = {7},
     PAGES = {2092--2116},
      ISSN = {0022-1236,1096-0783},
   MRCLASS = {53C21 (35J66 35R01)},
  MRNUMBER = {3989139},
MRREVIEWER = {Hong\ Zhang},
       DOI = {10.1016/j.jfa.2019.01.001},
       URL = {https://doi.org/10.1016/j.jfa.2019.01.001},
}

@article {AlmarazWang2022,
    AUTHOR = {de Moura Almaraz, S\'ergio and Wang, Shaodong},
     TITLE = {Energy bounds of sign-changing solutions to {Y}amabe equations
              on manifolds with boundary},
   JOURNAL = {Nonlinear Anal.},
  FJOURNAL = {Nonlinear Analysis. Theory, Methods \& Applications. An
              International Multidisciplinary Journal},
    VOLUME = {225},
      YEAR = {2022},
     PAGES = {Paper No. 113131, 12},
      ISSN = {0362-546X,1873-5215},
   MRCLASS = {35M30},
  MRNUMBER = {4475441},
       DOI = {10.1016/j.na.2022.113131},
       URL = {https://doi.org/10.1016/j.na.2022.113131},
}

@article{BianchiEgnell1991,
  AUTHOR = {Bianchi, Gabriele and Egnell, Henrik},
     TITLE = {A note on the {S}obolev inequality},
   JOURNAL = {J. Funct. Anal.},
  FJOURNAL = {Journal of Functional Analysis},
    VOLUME = {100},
      YEAR = {1991},
    NUMBER = {1},
     PAGES = {18--24},
      ISSN = {0022-1236,1096-0783},
   MRCLASS = {46E35},
  MRNUMBER = {1124290},
MRREVIEWER = {Petr\ Gurka},
       DOI = {10.1016/0022-1236(91)90099-Q},
       URL = {https://doi.org/10.1016/0022-1236(91)90099-Q},
}

@article {BorquezCajuBosch2026,
    AUTHOR = {Borquez, Benjam\'in and Caju, Rayssa and Van Den Bosch, Hanne},
     TITLE = {Quantitative stability for {Y}amabe minimizers on manifolds
              with boundary},
   JOURNAL = {J. Lond. Math. Soc. (2)},
  FJOURNAL = {Journal of the London Mathematical Society. Second Series},
    VOLUME = {113},
      YEAR = {2026},
    NUMBER = {2},
     PAGES = {Paper No. e70464, 25},
      ISSN = {0024-6107,1469-7750},
   MRCLASS = {58E15 (26D10 35J20 53C18)},
  MRNUMBER = {5033816},
       DOI = {10.1112/jlms.70464},
       URL = {https://doi.org/10.1112/jlms.70464},
}

@article {BrezisLieb1983,
    AUTHOR = {Br\'ezis, Ha\"im and Lieb, Elliott},
     TITLE = {A relation between pointwise convergence of functions and
              convergence of functionals},
   JOURNAL = {Proc. Amer. Math. Soc.},
  FJOURNAL = {Proceedings of the American Mathematical Society},
    VOLUME = {88},
      YEAR = {1983},
    NUMBER = {3},
     PAGES = {486--490},
      ISSN = {0002-9939,1088-6826},
   MRCLASS = {28A20 (46E30 49A99)},
  MRNUMBER = {699419},
MRREVIEWER = {R.\ G.\ Bartle},
       DOI = {10.2307/2044999},
       URL = {https://doi.org/10.2307/2044999},
}

@article {Cherrier1984,
    AUTHOR = {Cherrier, Pascal},
     TITLE = {Probl\`emes de {N}eumann non lin\'eaires sur les vari\'et\'es
              riemanniennes},
   JOURNAL = {J. Funct. Anal.},
  FJOURNAL = {Journal of Functional Analysis},
    VOLUME = {57},
      YEAR = {1984},
    NUMBER = {2},
     PAGES = {154--206},
      ISSN = {0022-1236},
   MRCLASS = {58G30 (35J65 53C21)},
  MRNUMBER = {749522},
MRREVIEWER = {Dennis\ M.\ DeTurck},
       DOI = {10.1016/0022-1236(84)90094-6},
       URL = {https://doi.org/10.1016/0022-1236(84)90094-6},
}

@article {Ekeland1979,
    AUTHOR = {Ekeland, Ivar},
     TITLE = {Nonconvex minimization problems},
   JOURNAL = {Bull. Amer. Math. Soc. (N.S.)},
  FJOURNAL = {American Mathematical Society. Bulletin. New Series},
    VOLUME = {1},
      YEAR = {1979},
    NUMBER = {3},
     PAGES = {443--474},
      ISSN = {0273-0979,1088-9485},
   MRCLASS = {49A50 (46N05 47H99 58B20)},
  MRNUMBER = {526967},
MRREVIEWER = {Alexander\ D.\ Ioffe},
       DOI = {10.1090/S0273-0979-1979-14595-6},
       URL = {https://doi.org/10.1090/S0273-0979-1979-14595-6},
}

@article{Escobar1988,
    AUTHOR = {Escobar, Jos\'e{} F.},
     TITLE = {Sharp constant in a {S}obolev trace inequality},
   JOURNAL = {Indiana Univ. Math. J.},
  FJOURNAL = {Indiana University Mathematics Journal},
    VOLUME = {37},
      YEAR = {1988},
    NUMBER = {3},
     PAGES = {687--698},
             DOI = {10.1512/iumj.1988.37.37033},
       URL = {https://doi.org/10.1512/iumj.1988.37.37033},
}

@article{EscobarAnnals92,
AUTHOR = {Escobar, Jos\'e{} F.},
     TITLE = {Conformal deformation of a {R}iemannian metric to a scalar
              flat metric with constant mean curvature on the boundary},
   JOURNAL = {Ann. of Math. (2)},
  FJOURNAL = {Annals of Mathematics. Second Series},
    VOLUME = {136},
      YEAR = {1992},
    NUMBER = {1},
     PAGES = {1--50},
      ISSN = {0003-486X,1939-8980},
   MRCLASS = {53C21 (53C25 58G30)},
  MRNUMBER = {1173925},
MRREVIEWER = {Dennis\ M.\ DeTurck},
       DOI = {10.2307/2946545},
       URL = {https://doi.org/10.2307/2946545},
}

@article{EscobarJDG92,
AUTHOR = {Escobar, Jos\'e{} F.},
     TITLE = {The {Y}amabe problem on manifolds with boundary},
   JOURNAL = {J. Differential Geom.},
  FJOURNAL = {Journal of Differential Geometry},
    VOLUME = {35},
      YEAR = {1992},
    NUMBER = {1},
     PAGES = {21--84},
      ISSN = {0022-040X,1945-743X},
   MRCLASS = {53C21 (35J20 58E11 58G30)},
  MRNUMBER = {1152225},
MRREVIEWER = {Dennis\ M.\ DeTurck},
       URL = {http://projecteuclid.org/euclid.jdg/1214447805},
}

@article {EscobarAddendum94,
    AUTHOR = {Escobar, Jos\'{e} F.},
     TITLE = {Addendum: ``{C}onformal deformation of a {R}iemannian metric
              to a scalar flat metric with constant mean curvature on the
              boundary'' [{A}nn. of {M}ath. (2) {\bf 136} (1992), no. 1,
              1--50; {MR}1173925 (93e:53046)]},
   JOURNAL = {Ann. of Math. (2)},
  FJOURNAL = {Annals of Mathematics. Second Series},
    VOLUME = {139},
      YEAR = {1994},
    NUMBER = {3},
     PAGES = {749--750},
      ISSN = {0003-486X,1939-8980},
   MRCLASS = {53C21 (53C25 58G30)},
  MRNUMBER = {1283876},
       DOI = {10.2307/2118578},
       URL = {https://doi.org/10.2307/2118578},
}

@article {Escobar1996a,
    AUTHOR = {Escobar, Jos\'e{} F.},
     TITLE = {Conformal metrics with prescribed mean curvature on the
              boundary},
   JOURNAL = {Calc. Var. Partial Differential Equations},
  FJOURNAL = {Calculus of Variations and Partial Differential Equations},
    VOLUME = {4},
      YEAR = {1996},
    NUMBER = {6},
     PAGES = {559--592},
      ISSN = {0944-2669,1432-0835},
   MRCLASS = {53C21 (35J60 53A30 58E15)},
  MRNUMBER = {1416000},
MRREVIEWER = {Robert\ McOwen},
       DOI = {10.1007/BF01261763},
       URL = {https://doi.org/10.1007/BF01261763},
}

@article {GhimentiMicheletti2020NA,
    AUTHOR = {Ghimenti, Marco and Micheletti, Anna Maria},
     TITLE = {Compactness for conformal scalar-flat metrics on umbilic
              boundary manifolds},
   JOURNAL = {Nonlinear Anal.},
  FJOURNAL = {Nonlinear Analysis. Theory, Methods \& Applications. An
              International Multidisciplinary Journal},
    VOLUME = {200},
      YEAR = {2020},
     PAGES = {111992, 30},
      ISSN = {0362-546X,1873-5215},
   MRCLASS = {58J32 (35R01 53C18 53C21)},
  MRNUMBER = {4103358},
MRREVIEWER = {Yunyan\ Yang},
       DOI = {10.1016/j.na.2020.111992},
       URL = {https://doi.org/10.1016/j.na.2020.111992},
}

@article {GhimentiMicheletti2021CVPDE,
    AUTHOR = {Ghimenti, Marco G. and Micheletti, Anna Maria},
     TITLE = {A compactness result for scalar-flat metrics on low
              dimensional manifolds with umbilic boundary},
   JOURNAL = {Calc. Var. Partial Differential Equations},
  FJOURNAL = {Calculus of Variations and Partial Differential Equations},
    VOLUME = {60},
      YEAR = {2021},
    NUMBER = {3},
     PAGES = {Paper No. 119, 24},
      ISSN = {0944-2669,1432-0835},
   MRCLASS = {53C21 (35J65)},
  MRNUMBER = {4270815},
MRREVIEWER = {S\'{e}rgio\ Almaraz},
       DOI = {10.1007/s00526-021-01983-x},
       URL = {https://doi.org/10.1007/s00526-021-01983-x},
}

@article{HanLi1999,
    AUTHOR = {Han, Zheng-Chao and Li, Yanyan},
     TITLE = {The {Y}amabe problem on manifolds with boundary: existence and
              compactness results},
   JOURNAL = {Duke Math. J.},
  FJOURNAL = {Duke Mathematical Journal},
    VOLUME = {99},
      YEAR = {1999},
    NUMBER = {3},
     PAGES = {489--542},
      ISSN = {0012-7094,1547-7398},
   MRCLASS = {53C21 (35J60 58J60)},
  MRNUMBER = {1712631},
MRREVIEWER = {David\ L.\ Finn},
       DOI = {10.1215/S0012-7094-99-09916-7},
       URL = {https://doi.org/10.1215/S0012-7094-99-09916-7},
}

@article {Ho2022,
    AUTHOR = {Ho, Pak Tung},
     TITLE = {A note on the {S}obolev trace inequality},
   JOURNAL = {Proc. Amer. Math. Soc.},
  FJOURNAL = {Proceedings of the American Mathematical Society},
    VOLUME = {150},
      YEAR = {2022},
    NUMBER = {3},
     PAGES = {1257--1267},
      ISSN = {0002-9939,1088-6826},
   MRCLASS = {46E35 (35A23)},
  MRNUMBER = {4375719},
       DOI = {10.1090/proc/15751},
       URL = {https://doi.org/10.1090/proc/15751},
}

@article {KimMussoWei2021,
    AUTHOR = {Kim, Seunghyeok and Musso, Monica and Wei, Juncheng},
     TITLE = {Compactness of scalar-flat conformal metrics on
              low-dimensional manifolds with constant mean curvature on
              boundary},
   JOURNAL = {Ann. Inst. H. Poincar\'e{} C Anal. Non Lin\'eaire},
  FJOURNAL = {Annales de l'Institut Henri Poincar\'e{} C. Analyse Non
              Lin\'eaire},
    VOLUME = {38},
      YEAR = {2021},
    NUMBER = {6},
     PAGES = {1763--1793},
      ISSN = {0294-1449,1873-1430},
   MRCLASS = {53C18 (35B40 35J65 35R01 53C21)},
  MRNUMBER = {4327897},
MRREVIEWER = {Gabjin\ Yun},
       DOI = {10.1016/j.anihpc.2021.01.005},
       URL = {https://doi.org/10.1016/j.anihpc.2021.01.005},
}

@article{LiZhu1995,
    AUTHOR = {Li, Yan Yan and Zhu, Meijun},
     TITLE = {Uniqueness theorems through the method of moving spheres},
   JOURNAL = {Duke Math. J.},
  FJOURNAL = {Duke Mathematical Journal},
    VOLUME = {80},
      YEAR = {1995},
    NUMBER = {2},
     PAGES = {383--417},
      ISSN = {0012-7094},
       DOI = {10.1215/S0012-7094-95-08016-8},
   MRCLASS = {35J60 (35B05 53C21)},
  MRNUMBER = {1369398},
}

@article{Lions85b,
  AUTHOR = {Lions, P.-L.},
     TITLE = {The concentration-compactness principle in the calculus of
              variations. {T}he limit case. {II}},
   JOURNAL = {Rev. Mat. Iberoamericana},
  FJOURNAL = {Revista Matem\'atica Iberoamericana},
    VOLUME = {1},
      YEAR = {1985},
    NUMBER = {2},
     PAGES = {45--121},
      ISSN = {0213-2230},
   MRCLASS = {49A22 (46E35 49F15 58E30)},
  MRNUMBER = {850686},
MRREVIEWER = {Ll.\ G.\ Chambers},
       DOI = {10.4171/RMI/12},
       URL = {https://doi.org/10.4171/RMI/12},
}

@article {Marques2007,
    AUTHOR = {Marques, Fernando C.},
     TITLE = {Conformal deformations to scalar-flat metrics with constant
              mean curvature on the boundary},
   JOURNAL = {Comm. Anal. Geom.},
  FJOURNAL = {Communications in Analysis and Geometry},
    VOLUME = {15},
      YEAR = {2007},
    NUMBER = {2},
     PAGES = {381--405},
      ISSN = {1019-8385,1944-9992},
   MRCLASS = {53C21 (53C44)},
  MRNUMBER = {2344328},
MRREVIEWER = {Mohameden\ Ahmedou},
       DOI = {10.4310/cag.2007.v15.n2.a7},
       URL = {https://doi.org/10.4310/cag.2007.v15.n2.a7},
}

@article{MayerNdiaye2017,
  AUTHOR = {Mayer, Martin and Ndiaye, Cheikh Birahim},
     TITLE = {Barycenter technique and the {R}iemann mapping problem of
              {C}herrier-{E}scobar},
   JOURNAL = {J. Differential Geom.},
  FJOURNAL = {Journal of Differential Geometry},
    VOLUME = {107},
      YEAR = {2017},
    NUMBER = {3},
     PAGES = {519--560},
      ISSN = {0022-040X,1945-743X},
   MRCLASS = {53A30 (53C21)},
  MRNUMBER = {3715348},
MRREVIEWER = {Yueh-Ju\ Lin},
       DOI = {10.4310/jdg/1508551224},
       URL = {https://doi.org/10.4310/jdg/1508551224},
}

@article {Smale1965,
    AUTHOR = {Smale, S.},
     TITLE = {An infinite dimensional version of {S}ard's theorem},
   JOURNAL = {Amer. J. Math.},
  FJOURNAL = {American Journal of Mathematics},
    VOLUME = {87},
      YEAR = {1965},
     PAGES = {861--866},
      ISSN = {0002-9327,1080-6377},
   MRCLASS = {57.55 (57.50)},
  MRNUMBER = {185604},
MRREVIEWER = {Richard\ Beals},
       DOI = {10.2307/2373250},
       URL = {https://doi.org/10.2307/2373250},
}

@article{Struwe1984,   
  AUTHOR = {Struwe, Michael},
     TITLE = {A global compactness result for elliptic boundary value
              problems involving limiting nonlinearities},
   JOURNAL = {Math. Z.},
  FJOURNAL = {Mathematische Zeitschrift},
    VOLUME = {187},
      YEAR = {1984},
    NUMBER = {4},
     PAGES = {511--517},
      ISSN = {0025-5874,1432-1823},
   MRCLASS = {35J60 (35A35 35J20)},
  MRNUMBER = {760051},
MRREVIEWER = {Bernhard\ Kawohl},
       DOI = {10.1007/BF01174186},
       URL = {https://doi.org/10.1007/BF01174186},
}

@article {ZhangZhouZou2025,
    AUTHOR = {Zhang, Yingfang and Zhou, Yuxuan and Zou, Wenming},
     TITLE = {Sharp quantitative stability for the fractional {S}obolev
              trace inequality},
   JOURNAL = {Math. Z.},
  FJOURNAL = {Mathematische Zeitschrift},
    VOLUME = {310},
      YEAR = {2025},
    NUMBER = {4},
     PAGES = {Paper No. 90, 45},
      ISSN = {0025-5874,1432-1823},
   MRCLASS = {35J25 (26A33)},
  MRNUMBER = {4920022},
       DOI = {10.1007/s00209-025-03788-6},
       URL = {https://doi.org/10.1007/s00209-025-03788-6},
}

@article{Chen2010,
      title={Conformal Deformation to Scalar Flat Metrics with Constant Mean Curvature on the Boundary in Higher Dimensions},
      journal={arXiv e-prints},
      author={Szu-yu Sophie Chen},
      year=2010,
      month=apr,
      eprint={0912.1302},
      archivePrefix={arXiv},
      note={arXiv:0912.1302},
}

@article{ChenKim2024,
      title={Sharp quantitative stability of the {Y}amabe problem},
      journal={arXiv e-prints},
      author={Haixia Chen and Seunghyeok Kim},
      year=2024,
      month=may,
      eprint={2404.13961},
      archivePrefix={arXiv},
      note={arXiv:2404.13961},
}

@article{GongKimMussoWei2026,
      title={On Critical Dimensions for Compactness in the Boundary {Y}amabe Problem, {II}},
      journal={arXiv e-prints},
      author={Liuwei Gong and Seunghyeok Kim and Monica Musso and Juncheng Wei},
      year=2026,
      month=oct,
      eprint={2610.10485},
      archivePrefix={arXiv},
      note={arXiv:2610.10485},
}

@article{KhuriWangWang2026,
      title={Riemannian Positive Mass Theorem in All Dimensions in the Presence of Low-Codimension Singularities}, 
      journal={arXiv e-prints},
      author={Marcus Khuri and Jian Wang and Jinmin Wang},
      year=2026,
      month=jul,
      eprint={2606.23529},
      archivePrefix={arXiv},
      note={arXiv:2606.23529},
}

@article{EscobarGeometry,
      title = {Quantitative stability and strict inequalities for the {E}scobar functional},
      journal = {arXiv e-prints},
      author = {Mayukh Mukherjee and Utsab Sarkar},
      year = 2026,
      month = oct,
      eprint = {2601.22665},
      archivePrefix = {arXiv},
      note= {arXiv:2601.22665},
}

@article {YuZhou2026,
      title={Sharp One-bubble Critical-Point Stability and Global Compactness for the {S}obolev Trace Inequality},
      journal={arXiv e-prints},
      author={Bao Yu and Yang Zhou},
      year=2026,
      month=aug,
      eprint={2607.21429},
      archivePrefix={arXiv},
      note={arXiv:2607.21429},
}

@article{MukherjeeSarkarGN2026,
  author = {Mukherjee, Mayukh and Sarkar, Utsab},
  title = {{G}agliardo--{N}irenberg extremals on manifolds with small holes},
  journal = {arXiv e-prints},
  year = {2026},
  month = oct,
  eprint = {2610.32173},
  archivePrefix = {arXiv},
  primaryClass = {math.AP},
  note = {arXiv:2610.32173},
  url = {https://arxiv.org/abs/2610.32173}
}
\end{document}